\documentclass[english]{article}

\usepackage{geometry}
\usepackage{amsmath,bbm,bm,amssymb}
\usepackage{amsfonts}
\usepackage{amsthm}
\usepackage{mathtools}
\usepackage{authblk,natbib}
\usepackage[colorlinks=true, linkcolor=blue, citecolor=blue, urlcolor=blue]{hyperref}

\usepackage{tabularx}

\usepackage{amsmath,bbm,bm}
\usepackage{amsfonts}
\usepackage{amsthm}
\usepackage{mathtools}

\usepackage{algorithm}
\usepackage{algorithmic}

\newtheorem{thm}{Theorem}
\newtheorem{lem}{Lemma}
\newtheorem{prop}{Proposition}
\newtheorem{cor}{Corollary}

\newtheorem{aspt}{Assumption}

\newtheorem{defn}{Definition}

\usepackage{xcolor}
\newcount\comments  
\newcommand{\genComment}[2]{\ifnum\comments=1{\textcolor{#1}{\textsf{\footnotesize #2}}}\fi}

\def\eqref#1{equation~\ref{#1}}

\def\1{\bm{1}}

\DeclareMathAlphabet{\mathsfit}{\encodingdefault}{\sfdefault}{m}{sl}
\SetMathAlphabet{\mathsfit}{bold}{\encodingdefault}{\sfdefault}{bx}{n}

\def\gA{{\mathcal{A}}}

\def\gC{{\mathcal{C}}}

\def\gX{{\mathcal{X}}}

\def\R{{\mathbb{R}}}

\DeclareMathOperator*{\argmin}{arg\,min}

\def\bx{{\bm{x}}}

\def\bI{{\bm{I}}}

\def\ba{{\bm a}}
\def\bA{{\bm{A}}}

\def\btheta{{\bm \theta}}
\def\bSigma{{\bm \Sigma}}
\def\bphi{{\bm \phi}}

\DeclareSymbolFont{matha}{OML}{txmi}{m}{it}
\DeclareMathSymbol{\varv}{\mathord}{matha}{118}

\newcommand\bLambda{\boldsymbol{\Lambda}}

\newcommand{\thetastar}{{\bm \theta}^{\star}}

\usepackage{color}
\definecolor{yxc}{RGB}{255,0,0}
\definecolor{yjc}{RGB}{125,0,0}
\definecolor{ytw}{RGB}{255,69,0}
\definecolor{gen}{RGB}{0,0,200}
\definecolor{wei}{RGB}{0,125,125}

\allowdisplaybreaks

\theoremstyle{definition}

\title{Statistical Inference under Adaptive Sampling with LinUCB}
\author{Wei Fan$^{*}$, Kevin Tan\thanks{Equal contribution.}, Yuting Wei}
\affil{
  Department of Statistics and Data Science\\
  The Wharton School, University of Pennsylvania}
\date{\today}

\begin{document}

\maketitle

\begin{abstract}
    Adaptively collected data has become ubiquitous within modern practice. However, even seemingly benign adaptive sampling schemes can introduce severe biases, rendering traditional statistical inference tools inapplicable. This can be mitigated by a property called stability, which states that if the rate at which an algorithm takes actions converges to a deterministic limit, one can expect that certain parameters are asymptotically normal. Building on a recent line of work for the multi-armed bandit setting, we show that the linear upper confidence bound (LinUCB) algorithm for linear bandits satisfies this property. In doing so, we painstakingly characterize the behavior of the eigenvalues and eigenvectors of the random design feature covariance matrix in the setting where the action set is the unit ball, showing that it decomposes into a rank-one direction that locks onto the true parameter and an almost-isotropic bulk that grows at a predictable $\sqrt{T}$ rate. This allows us to establish a central limit theorem for the LinUCB algorithm, establishing asymptotic normality for the limiting distribution of the estimation error where the convergence occurs at a $T^{-1/4}$ rate. The resulting Wald-type confidence sets and hypothesis tests do not depend on the feature covariance matrix and are asymptotically tighter than existing nonasymptotic confidence sets. Numerical simulations corroborate our findings.
\end{abstract}

\tableofcontents

\section{Introduction}

Statistical inference for adaptively collected data is essential for providing rigorous justification and interpretability in modern data analysis, with applications ranging from scientific discovery to social decision-making. 
In sharp contrast to classical i.i.d.\ settings, adaptive data collection induces intricate dependencies across samples, often rendering traditional inferential tools unreliable. 
It is now well recognized that even seemingly benign adaptive sampling schemes can introduce severe biases and complicate the asymptotic distribution of estimators 
(see, e.g.~\cite{dickey1979distribution,lai1982least,deshpande2023online}), which in turn complicates the task of uncertainty quantification (\cite{deshpande2018accurate,khamaru2024inference,lin2023statistical,zhang2020inferenceforbatchedbandits,zhang2021statisticalinferencewithmestimatorsonadaptivelycollecteddata}).

A central theme emerging from this growing body of work is that the very process of learning dynamically reshapes the statistical properties of the data. Although this complicates the analysis, recent work within the multi-armed bandit setting \citep{kalvit2021closerlookworstcasebehavior, khamaru2024inference, han2024ucb, halder2025stable} has shown that certain algorithms exhibit a notion of ``stability'' that allows for asymptotic normality of the arm mean reward estimates. In other words, if the rate at which a multi-armed bandit algorithm pulls each arm is asymptotically deterministic, then under suitable conditions this alone can ensure that the estimated mean rewards are asymptotically normal. This property is satisfied for the UCB algorithm \citep{khamaru2024inference, kalvit2021closerlookworstcasebehavior, han2024ucb}, but not Thompson sampling \citep{zhang2021statisticalinferencewithmestimatorsonadaptivelycollecteddata} unless the posterior variance is inflated by a logarithmic factor \citep{halder2025stable}. 


We explore whether this phenomena of stability for the UCB algorithm \citep{khamaru2024inference, kalvit2021closerlookworstcasebehavior, han2024ucb} also extends to the linear bandit problem, a classical and influential model in reinforcement learning and the bandit literature. The linear bandit formalizes sequential decision-making in which the expected reward is a linear function of an action’s features, making it both a natural and useful abstraction. 
Mathematically, given a content $\bx_t \in \gX$, and an action $\ba_t \in \gA$ at time $t$, the learner receives the reward 
\begin{align}
    \label{equ:linear}
    r_t = \langle \bphi(\bx_t,\ba_t),{\bm \theta}^{\star}\rangle+{ \epsilon}_t \in \mathbb{R}, 
\end{align}
according to an unknown parameter ${\bm \theta}^{\star} \in \R^d$.
Here, $\bphi: \gX \times \gA \to \R^d$ is a feature map, and the action ${\ba}_t \in \mathcal{F}_{t-1}$ with $\mathcal{F}_{t-1}$ being the $\sigma$-field generated by the history including previous actions and rewards $\{{\ba}_1, r_1,\ldots, {\ba}_{t-1}, r_{t-1}\}$, collected up to time $t-1$. 
The noise satisfies $\mathbb{E}[{ \epsilon}_t|\mathcal{F}_{t-1}] = 0$.

Within this linear bandit framework, we focus on the linear upper confidence bound algorithm (\emph{LinUCB}) of \citet{Li_2010,abbas2011improved}, a canonical UCB-type method tailored to linear bandits. LinUCB stands out for its principled balance between exploration and exploitation and for its strong practical performance. Yet while LinUCB’s regret guarantees are well understood, the distributional behavior of its estimators—and, consequently, tools for valid statistical inference under LinUCB—remain underdeveloped. A common fallback is the familiar non-asymptotic confidence set for $\btheta^\star$ (also used for action selection), but because it depends on the empirical feature covariance $\bLambda_T$, it provides neither a limiting distribution for the estimator nor a deterministic characterization of the set’s width. Beyond this stopgap, only a handful of results from stochastic approximation or stochastic gradient descent \citep{polyak1992acceleration,su2023higrad,wu2025uncertainty,chen2020statistical} and a covariance characterization within LinUCB \citep{banerjee2023exploration} speak to inference, and these remain insufficient for conducting valid statistical inference with LinUCB. This gap motivates the following question:
\\

{\centering \emph{Can LinUCB be used not only as a learning algorithm, but also as a vehicle for valid statistical inference?}}\\

\subsection{Prior art}
Existing results do not provide a complete answer to this question. 
It has long been known that the confidence sets constructed by LinUCB contain $\bm{\theta}^*$ with high probability at each timestep. More precisely, when $\lVert\bphi(\bx_t, \ba_t)\rVert_2 \leq L$ for all $t=1,...,T$, $\Vert\bm{\theta}^{\star}\Vert_2\leq S$ and $\bLambda_t$ is the feature covariance matrix at time $t$ (defined in (\ref{equ:cov})), Theorem 2 of \cite{abbas2011improved} states that, given  some regularization parameter $\lambda > 0$, for all $t=1,...,T$, with probability at least $1-\delta$:
\begin{align}
     \bm{\theta}^*\in \gC_t, \qquad \gC_t = \left\{\btheta \in \R^d \;\;\bigg|\;\; \lVert \widehat{\btheta}_t - \btheta \rVert_{\bLambda_t} \leq \sigma\sqrt{d\log\left(\frac{1+TL^2/\lambda}{\delta}\right)} + \lambda^{1/2} S\right\},
     \label{eq:linucb-confidence-sets}
\end{align}
  where $\Vert\ba\Vert_{\bm{\Lambda}_t}\coloneqq\sqrt{\ba^\top\bm{\Lambda}_t\ba}$. This yields a simultaneous confidence set for all coefficients within the parameter $\btheta^*$ with non-asymptotic $1-\delta$ coverage. However, this simply states that $\btheta^*$ is contained within some ellipsoid centered at $\widehat\btheta_t$ and beyond that, we have no information about the distribution of $\widehat{\btheta}_t - \btheta^*$. Further, this non-asymptotic result provides no information on how the eigenvalues of $\bLambda_t$ scale with $t$, leaving the dependence of the size of the confidence set on $t$ unclear. 
  We fill this gap within this paper, both by characterizing the asymptotic limiting distribution and by providing non-asymptotic guarantees for convergence to said limiting distribution. Unlike their confidence set, ours explicitly utilizes the quantile of the chi-squared distribution, entailing a stronger distributional result.

Regarding other known results, Theorem 3 in \cite{lai1982least} establishes the asymptotic normality of the least squares estimator for the adaptive linear regression model. Here, in contrast to the multi-armed bandit setting, the ``stability'' condition requires that the sample covariance matrix (and not the sequence of arm pull rates) stabilizes to a deterministic sequence. They also demonstrate that asymptotic normality can fail in the absence of this stability property.
In the multi-armed bandit (MAB) setting, a special case of linear bandits considered here, \cite{khamaru2024inference} verifies the stability property and proves asymptotic normality for the UCB algorithm. In a broader context, \cite{banerjee2023exploration} controls the minimum eigenvalues of the design matrix generated by any linear bandit algorithm with sublinear regret. As we shall see shortly, this quantity also plays a central role in our analysis.
However, it remains unclear whether the output of the LinUCB algorithm is asymptotically normal without additional assumptions, and if so, what the limiting covariance structure would be. This question is appealing both theoretically and practically: if asymptotic normality holds, one can construct substantially tighter confidence sets than those derived from concentration inequalities.

\begin{figure}[t]
    \centering
    \includegraphics[width=0.7\linewidth]{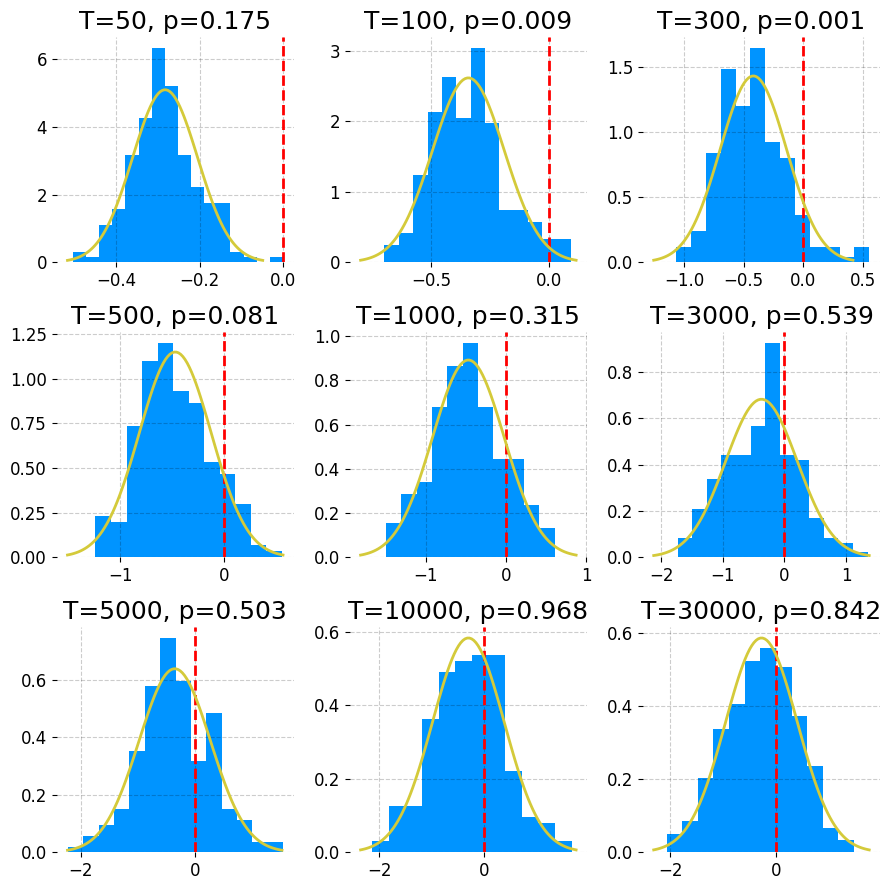}
    \caption{Asymptotic normality of the LinUCB algorithm in case where the action set is the unit ball. For some random vector $u$ on the unit ball, we plot $\widehat{\sigma}^{-1} \big(\frac{2\beta^2 T}{d+1}\big)^{1/4} u^\top ( \widehat{\btheta}_T - \btheta^\star)$ over 1000 independent trials, with KDE estimate overlaid as well as Shapiro-Wilk $p$-values provided as a test for non-normality. Asymptotic normality is indeed demonstrated, but the rate of convergence to the true parameter is certainly empirically slower than the $1/\sqrt{T}$ parametric rate, corroborating our theory.}
    \label{fig:normality-ball}
\end{figure}

\subsection{Our contributions}

In this paper, we study the asymptotic behavior of the LinUCB algorithm in the setting where the true parameter lies on the unit sphere and the action set is the unit ball (Assumption~\ref{aspt:unconstrained}), under sub-Gaussian noise (Assumption~\ref{aspt:subgaussian}).
We show that the algorithm is asymptotically normal, aligning with empirical results in Figure \ref{fig:normality-ball}.
\begin{itemize}
    \item Our main result, Theorem~\ref{thm:ucb-unconstrained} shows that, when projected onto  $(\bm{\theta}^{\star})^{\perp}$ (the orthogonal complement of $\boldsymbol{\theta}^\star$), the asymptotic covariance of $\widehat{\boldsymbol{\theta}}_{T}$ is isotropic—of magnitude $\sqrt{\tfrac{d+1}{2\beta^2 T}}\,\boldsymbol{I}_{d-1}$ for a broad class of exploration schedules $\beta=\beta(T,d)$. Equivalently, after an explicit rescaling proportional to $\beta$, the projected error obeys a central limit theorem: for any $\boldsymbol{U}\in\mathbb{R}^{d\times(d-1)}$ with orthonormal columns orthogonal to $\boldsymbol{\theta}^\star$,
\[
\left(\frac{2 \beta^2 T}{d+1}\right)^{\!1/4}\,\boldsymbol{U}^{\top}\!\left(\widehat{\boldsymbol{\theta}}_{T}-\boldsymbol{\theta}^{\star}\right) \xrightarrow{d} \mathcal{N}\!\left(0, \sigma^2 \boldsymbol{I}_{d-1}\right).
\]
Since $\widehat{\boldsymbol{\theta}}_{T}$ lies on the unit sphere, this essentially pins down the full asymptotic law of $\widehat{\boldsymbol{\theta}}_{T}$. To our knowledge, these results give the first asymptotic convergence guarantee for the parameter estimate based on a dependent, adaptively collected data sequence generated by the LinUCB policy.


    \item The above result allows us to provide an asymptotic $(1-\delta)$ Wald-type confidence set for $\bm{\theta}^{\star}$. In contrast to the confidence set from \cite{abbas2011improved} in  (\ref{eq:linucb-confidence-sets}), we require only an asymptotically spherical confidence set
  $$
\mathcal{C}_\delta=\left\{\boldsymbol{\theta} \in \mathcal{S}^{d-1}:\left\|\widehat{\boldsymbol{\theta}}_T-\boldsymbol{\theta}\right\|_2^2 \leq \widehat{\sigma}^2 \sqrt{\frac{d+1}{2 \beta^2 T}} \cdot \chi_{d-1,1-\delta}^2\right\}, 
$$
    for which we can provide a precise and deterministic (modulo randomness in the variance estimate $\widehat{\sigma}^2$) characterization of its diameter.  
    If the user desires an ellipsoid confidence set for better finite-sample performance, they can utilize
    $$\mathcal{C}_{\delta}^{\mathsf{ellipsoid}} = \left\{\bm{\theta}\in\mathcal{S}^{d-1}: \left\Vert \widehat{\bm{\theta}}_{T} - \bm{\theta}\right\Vert^2_{\bm{\Lambda}_T}\leq \widehat{\sigma}^2\chi_{d-1,1-\delta}^2\right\},$$
    which is asymptotically equivalent to $\mathcal{C}_\delta$ while maintaining the same distributional guarantees.
   Whereas the LinUCB confidence set in (\ref{eq:linucb-confidence-sets}) is derived via martingale concentration, our sets are standard Wald-type confidence sets used for asymptotically normal models under i.i.d.\ sampling. A key caveat is the rate: with \(T\) observations, the estimator concentrates at \(O_p(T^{-1/4})\), in contrast to the \(O_p(T^{-1/2})\) rate (and corresponding confidence-set width) in the i.i.d.\ setting. This perspective allows us to conduct statistical inference on data collected by LinUCB, paralleling the classical UCB algorithm for multi-armed bandits \citep{khamaru2024inference}.

    \item In the process of proving this result, we ended up proving several results that may be of independent interest.  These include:
    \begin{itemize}
        \item A tighter uniform control over the error of the estimated parameter at each timestep within Theorem \ref{thm:noise-bound}. Instead of obtaining a high-probability uniform bound for $\lVert \widehat{\btheta}_t - \btheta^\star \rVert_{\bLambda_t}$ ($t\in[T]$) scaling in $O(\sqrt{\log T})$ as shown in~(\ref{eq:linucb-confidence-sets}), we obtain an improved bound of $O( \sqrt{ \log \log T})$ when the failure probability is $1/\log T$. 
        \item A complete characterization of the eigenvalues of the feature covariance $\boldsymbol{\Lambda}_t$ under LinUCB. 
In particular, we describe how these eigenvalues evolve throughout the learning process, summarized in Propositions~\ref{prop:first-stage}--\ref{prop:fourth-stage}.  First, in Propositions~\ref{prop:first-stage} and \ref{prop:second-stage}, we establish the early-stage pattern: all non-leading eigenvalues are of the same order -- $\Theta(t)$ in the first phase and $\Theta(\sqrt{t})$ in the second. Next, Proposition~\ref{prop:third-stage} identifies the key transition, wherein the top eigenvector concentrates more tightly around $\boldsymbol{\theta}^\star$; building on this transition, Propositions~\ref{prop:fourth-stage} develop a fine-grained analysis of the non-leading eigenvalues and  show that they converge to a deterministic limit.

        
    \end{itemize} 
    Although analogous results are available for the multi-armed bandit setting with the UCB algorithm~\citep{khamaru2024inference}, that setting features a finite, fixed action set. In contrast, our linear bandit setting allows arbitrary adaptive exploration over the unit ball, yielding a continuum of actions. This richer action space makes the asymptotic analysis substantially more delicate and requires techniques beyond those used in the finite-arm case.
\end{itemize}


\subsection{Other related works}

\paragraph{Bandit algorithms and inference for bandits. } Dating back to the seminal works \cite{robbins1952some,thompson1933likelihood}, bandit algorithms have attracted tremendous attention for their simplicity and flexibility in modeling adaptive data collection. The UCB algorithm was first proposed in \cite{lai1985asymptotically,lai1987adaptive}, and its linear extension, LinUCB, was introduced in \cite{Li_2010}. UCB and its many variants have since been widely applied to problems with dynamic data, leaving a profound impact across statistics, operations research, and reinforcement learning. Classical research has primarily focused on regret analysis, often establishing sublinear bounds for UCB and its extensions (see \cite{bubeck2012regret,lattimore2020bandit} and references therein). More recent work has refined these results, providing precise regret bounds (\cite{han2024ucb,fan2024precise}), and in some cases explicitly incorporating stability to enable tractable inference alongside learning (\cite{sengupta2024stablebatchedbanditoptimalregretwithfreeinference}). While the inferential properties of UCB have been recently studied in the multi-armed bandit setting \cite{khamaru2024inference}, much less is known about the inferential properties of LinUCB, beyond the non-asymptotic guarantees provided by confidence sets \cite{abbas2011improved}.

\paragraph{Inference with adaptively collected data. } 
As discussed above, when data are collected adaptively, as in multi-armed and contextual bandits, standard i.i.d. asymptotics fail, since sampling depends on past observations, introduces complex dependencies (\cite{lai1982least,dwork2015reusable}). 
To address this challenge, a growing body of literature develops new statistical inference procedures that come with theoretical guarantees. 
Broadly, depending on the flavor of these results, they can be classified into two categories: those providing finite-sample guarantees and those establishing asymptotic characterizations. High probability bounds that hold for any finite samples often rely on tools such as concentration of self-normalized martingales (e.g.~\cite{abbas2011improved,shi2023statistical,waudbysmith2024anytimevalidoffpolicyinferenceforcontextualbandits,nair2023randomizationtestsforadaptivelycollecteddata,dimakopoulou2021online,wu2024statistical}). 
In contrast, the asymptotic line of work focuses on characterizing the limiting distribution of estimators of interest, enabling confidence intervals and hypothesis testing in large-sample regimes (e.g.~\cite{hadad2021confidence,halder2025stable,deshpande2018accurate,zhang2020inferenceforbatchedbandits,niu2025assumption,wu2025uncertainty,guo2025statistical}).


\subsection{Organization and notation}
\paragraph{Paper organization.} The remainder of the paper is organized as follows. Section~\ref{sec:preliminary} reviews basics for linear bandits and the LinUCB algorithm, and introduces the stability condition necessary for establishing asymptotic results under LinUCB and other adaptive data-collection schemes.
Section~\ref{sec:main-result} presents our main results—an asymptotic normality theorem and an associated confidence set—along with further remarks and a proof sketch. Section~\ref{sec:covariance} analyzes the evolution of the design covariance under LinUCB, a key ingredient in establishing our main result. Section~\ref{sec:conclusion} concludes our paper with a discussion and outlines several directions for future work.

\paragraph{Notation.} Throughout this paper, \(\|\bx\|_2\) (or simply \(\|\bx\|\)) denotes the Euclidean norm for a vector \(\bx\), and \(\|\bA\|_2\) (or \(\|\bA\|\)) denotes the spectral (matrix 2-) norm for a matrix \(\bA\). For a positive definite matrix \(\bLambda\) and a vector \(\bx\), write the weighted norm as \(\|\bx\|_{\bLambda} := \sqrt{\bx^\top \bLambda \bx}\). Let \(\mathcal{S}^{d-1} := \{\bx \in \mathbb{R}^d : \|\bx\|_2 = 1\}\) be the unit sphere in \(\mathbb{R}^d\) and \(\mathcal{B}^d := \{\bx \in \mathbb{R}^d : \|\bx\|_2 \le 1\}\) be the unit ball in \(\mathbb{R}^d\). Define the projection onto the unit sphere by \(\mathcal{P}(\bx) := \bx/\|\bx\|_2\) for \(\bx \neq \mathbf{0}\). For any vector $\bm{x}\in\mathbb{R}^{d}$, define the orthogonal complement of $\bm{x}$ as $\bm{x}^{\perp} = \{\bm{y}\in\mathbb{R}^d: \langle\bm{x},\bm{y}\rangle = 0\}$, a $(d-1)$ dimensional subspace of $\mathbb{R}^d$.

In addition, for any two functions $f(T)$ and $g(T)$, we write $f(T)\lesssim g(T)$ (equivalently, $f(T)=O(g(T))$) if there exists a constant $c_1>0$ such that $|f(T)|\le c_1|g(T)|$. Conversely, we denote $f(T)\gtrsim g(T)$ (equivalently, $f(T)=\Omega(g(T))$) if there exists a constant $c_2>0$ such that $|f(T)|\ge c_2|g(T)|$. We also adopt $f(T)\asymp g(T)$ (equivalently, $f(T)=\Theta(g(T))$) to indicate that both $f(T)\lesssim g(T)$ and $f(T)\gtrsim g(T)$ hold simultaneously. We write $f(T) = O_p(g(T))$ if $f(T)/g(T)$ bounded in probability as $T\to\infty$. We write $f(T)=\widetilde{O}(g(T))$ to indicate that the $f(T) = O(g(T))$ holds up to logarithmic factors. Moreover, we denote $f(T)=o(g(T))$ if $f(T)/g(T)\to 0$ as $T\to\infty$, and we denote $f(T)\gg g(T)$ if $f(T)/g(T)\to\infty$ as $T\to\infty$. Finally, $c$ and $C$ denote universal constants that do not depend on $T$.














\section{Preliminaries}
\label{sec:preliminary}

\subsection{Stochastic linear bandits}

We formally define the stochastic linear bandits in $d$ dimensions. This framework models sequential decision-making where the expected reward of each action is a linear function of its associated features. The learning procedure unfolds over a time horizon of $T$ rounds. For clarity in what follows, we use the subscript $T$ to denote terminal quantities (i.e., those evaluated after $T$ rounds), and the subscript $t\in[T]:={1,\ldots,T}$ to denote per-round quantities.

In each round $t=1,2,\ldots,T$, the learner observes a context $\bx_t \in \gX$, and is presented a finite or infinite action set $\mathcal{A}_t$. We denote by $\gX$ the context space and by $\gA \supseteq \bigcup_{t=1}^T \gA_t$ the overall action space, i.e., the union of all actions available across rounds.
Given context $\bx \in \gX$, and action $\ba \in \gA$, we assume access to a feature mapping $\bphi: \gX \times \gA \to \R^d$. We define
\begin{align}
\label{eqn:feature-map}
\bm{\Phi}_t = \Big\{\bm{\phi}(\bm{x}_t,\bm{a}_t): \bm{a}_t\in\mathcal{A}_t\Big\}
\end{align}
be the set of features available to the learner at time $t$. 

As briefly introduced in (\ref{equ:linear}), the learner plays an action ${\ba}_t\in\mathcal{A}_t$ based on the trajectory of previous actions and rewards $\tau^{(t)} \coloneqq \{{\ba}_1, r_1,\ldots, {\ba}_{t-1}, r_{t-1}\}$ and receives a reward
\begin{align*}
    r_t = \langle \bphi(\bx_t,\ba_t),{\bm \theta}^{\star}\rangle+{ \epsilon}_t,
\end{align*}
where ${\bm \theta}^{\star}$ is an unknown parameter that defines the expected reward function. 
We assume without loss of generality that \(\|\bm{\theta}^{\star}\|_2 = 1\).
The expected reward $\langle  \bphi(\bx_t,\ba_t),{\bm \theta}^{\star}\rangle$ is a linear function with respect to ${\bm \theta}^{\star}$ and ${ \epsilon}_t$ is noise that satisfies $\mathbb{E}[{ \epsilon}_t|\mathcal{F}_{t-1}] = 0$. 
The feature map can be a neural embedding, random Fourier feature map, polynomial embedding, kernel embedding, or other similar feature maps. All we shall require is that it is known that the features are bounded at every timestep, and that the rewards are a linear function of the features. As such, and especially when pretrained embeddings are readily available, the linear bandit model can be surprisingly expressive.

The learner's goal is to minimize the cumulative regret, which characterizes the difference between the total expected reward and the best possible reward that the learner could possibly obtained. Formally speaking, the regret after round $T$ is defined as
\begin{align}
    R_T := \sum_{t=1}^{T} \langle \bphi(\bx_t, {\ba}_t^{\star}), {\bm \theta}^{\star} \rangle - \sum_{t=1}^{T} \langle \bphi(\bx_t, {\ba}_t), {\bm \theta}^{\star}\rangle, 
\end{align}
where ${\ba}_t^{\star} = \mathrm{argmax}_{{\ba}\in\mathcal{A}_t} \langle \bphi(\bx_t,{\ba}),{\bm \theta}^{\star}\rangle$, is defined as the best possible action that the learner could take at time $t$.
When the action sets $\gA_1 =...=\gA_T = \gA$ are the same for all $t=1,...,T$ and the context is the same for all $t$ as well,  the expression then simplifies:
\begin{align}
    R_T := \sum_{t=1}^{T} \langle \bphi(\bx, {\ba}^{\star}) - \bphi(\bx, \ba_t), {\bm \theta}^{\star} \rangle, \qquad {\ba}^{\star} = \mathrm{argmax}_{{\ba}\in\mathcal{A}} \langle \bphi(\bx, {\ba}),{\bm \theta}^{\star}\rangle. 
\end{align}

\subsection{The LinUCB algorithm}


A widely used and conceptually elegant strategy for minimizing regret in linear bandits is the Upper Confidence Bound (UCB) principle, adapted to the linear setting as \emph{LinUCB}.  The method embodies \emph{optimism in the face of uncertainty}: the learner maintains a high-probability confidence region for the unknown parameter \(\bm\theta^\star\in\mathbb{R}^d\) and acts as if the most favorable parameter in this region were the truth.  Concretely, at each round \(t\) the learner assigns to every candidate action \(\ba\in\mathcal{A}_t\) a UCB score that trades off predicted reward and an uncertainty bonus:
\begin{align}
\label{equ:ucb}
\mathrm{UCB}_t(\ba)
\;=\;
\big\langle \bphi(\bx_t,\ba),\, \widehat{\bm{\theta}}_{t-1}\big\rangle
\;+\;
\beta\,\sqrt{\bphi(\bx_t,\ba)^\top \bm{\Lambda}_{t-1}^{-1}\bphi(\bx_t,\ba)}\,,
\end{align}
where the first term is the estimated reward and the second is a data-dependent exploration bonus measuring uncertainty along the direction \(\bphi(\bx_t,\ba)\). To rationalize this principle, we  detail these two ingredients of the UCB score in  (\ref{equ:ucb}).

\paragraph{Estimated reward.} 
At time \(t{-}1\), define the cumulative covariance matrix
\begin{align}
\label{equ:cov}
\bLambda_{t-1}
= \lambda \bI_d + \sum_{s=1}^{t-1} \bphi(\bx_s,\ba_s)\bphi(\bx_s,\ba_s)^{\top}.
\end{align}
The (ridge) regularized least-squares estimator is then
\begin{align}
\label{equ:overline-theta}
\overline{\bm{\theta}}_{t-1}
\;\in\;
\argmin_{\bm\theta\in\mathbb{R}^d}
\Big\{\sum_{s=1}^{t-1}\!\big(r_s-\langle\bphi(\bx_s,\ba_s),\bm\theta\rangle\big)^2
+\lambda\|\bm\theta\|_2^2\Big\}
\;=\;
\bLambda_{t-1}^{-1}\!\sum_{s=1}^{t-1}\bphi(\bx_s,\ba_s)\,r_s,
\end{align}
Because the ground truth \(\bm\theta^\star\) has unit norm, we project the ridge estimate onto the unit sphere \(\mathcal{S}^{d-1}\):
\begin{align}
\label{equ:hat-theta}
\widehat{\bm{\theta}}_{t-1}
\;=\;\mathcal{P}\big(\overline{\bm{\theta}}_{t-1}\big),
\end{align}
With \(\widehat{\bm{\theta}}_{t-1}\) in hand, we score any candidate action \(\ba\) at time \(t\) by the estimated reward $\big\langle \bphi(\bx_t,\ba),\,\widehat{\bm{\theta}}_{t-1}\big\rangle$. In practice, it is common to set \(\lambda=1\) and \(\bLambda_0=\bI_d\), which keeps the cumulative covariance \(\bLambda_t\) invertible in the early stages while balancing bias and variance.

\paragraph{Exploration bonus.} The bonus term is set to be the conventional choice $\beta\sqrt{\bphi(\bx_t, \ba)^\top \bm{\Lambda}_{t-1}^{-1} \bphi(\bx_t, \ba)}$. 
The scalar \(\beta\) is chosen so that, with probability at least \(1-\delta\), the expected reward of every action is upper-bounded by \(\mathrm{UCB}_t(\ba)\) uniformly over \(t\) and \(\ba\) (the standard optimism property). 
For example, it is sufficient to choose $\beta$ as \citep{abbas2011improved}
\begin{align}
\label{eqn:beta-choice}
     \beta = \sigma \sqrt{d \log\bigl(1 + TL^2/d\bigr) + 2 \log(1/\delta)} + 1,
\end{align}
where \(\sigma\) is the sub-Gaussian parameter of the noise, and \(L\) is an upper bound for \(\|\bphi(\bx_t,\ba_t)\|_2\). 
With such a schedule, LinUCB achieves high-probability regret \(\widetilde{O}(d\sqrt{T})\) when \(L\) is constant \citep{abbas2011improved}.

Given the UCB score in  (\ref{equ:ucb}), LinUCB selects at time \(t\)
\[
\ba_t \in \arg\max_{\ba \in \mathcal{A}_t} \mathrm{UCB}_t(\ba),
\]
 the action whose upper confidence bound on the expected reward is largest, thereby optimistically balancing exploitation of high estimated rewards with exploration of uncertain actions. The learner then observes the immediate reward \(r_t\) and update UCB score for the next round. The procedure is summarized in Algorithm~\ref{alg:linucb}.

\subsection{Asymptotic normality of LinUCB with stability condition}
\label{sec:asymptotics-ucb} 


Beyond minimizing cumulative regret, LinUCB also serves as a natural candidate for statistical inference purposes. LinUCB outputs an estimator $\widehat{\bm\theta}_T$ at time $T$ via~(\ref{equ:hat-theta}) with the sequence of feature vectors $\bphi(\bx_1,\ba_1),\ldots,\bphi(\bx_T,\ba_T)$ and rewards $r_1,\ldots,r_T$ collected adaptively. 
If one can characterize the limiting distribution of \(\widehat{\bm\theta}_T\), Wald-type confidence sets can be constructed for $\thetastar.$
In the fixed-design case, where the action sequence is chosen a priori, if the terminal design $\bm{\Lambda}_T$ satisfies $\lambda_{\min}(\bm{\Lambda}_T)\to\infty$ and the noise is i.i.d. with mean zero and variance $\sigma^2$, then the ridge estimator $\overline{\bm\theta}_T$ obeys $\bm{\Lambda}_T^{1/2}(\bm{\overline{\theta}}_T - {\bm\theta}^{\star})\to \mathcal{N}(0,\sigma^2\bm{I}_d)$ (see Example 2.28 in \cite{van2000asymptotic}).
The same asymptotic normality holds in the i.i.d. random-design setting, provided the design covariance matrix is full rank.

When samples are collected adaptively (e.g., using LinUCB), the situation is more subtle. The terminal design $\bm{\Lambda}_T$ is random and history-dependent; its spectrum can vary across runs, so a single deterministic normalization under which $\overline{\bm{\theta}}_{T}-\bm{\theta}^{\star}$ has a normal limit may not exist; in particular, the usual Lindeberg–Feller CLT for deterministic designs is not directly applicable. 

In prior work, additional regularity conditions are made to ensure asymptotic normality. One such condition is the so-called \emph{stability} condition (\citet{lai1982least}), which assumes that the cumulative covariance matrix $\{\bm{\Lambda}_T\}$ admits a deterministic limit:
\begin{defn}[Stability]\label{defn:stable}
    The sequence of sample covariance matrices $\{\bm{\Lambda}_T\}$ is stable if there exists a sequence of deterministic positive definite matrices $\{{\bm \Sigma}_T\}$ such that 
    ${\bm \Sigma}_T^{-1} \bm{\Lambda}_T\longrightarrow {\bI}_d.$
\end{defn}

 
It is unclear that without the stability assumption or other similar assumptions, whether the asymptotic normality and, therefore statistical inference can be achieved. 

\begin{algorithm}[t]
\caption{Linear UCB Algorithm}
\begin{algorithmic}[1]
    \STATE \textbf{Input:} Horizon \( T \), action set $\mathcal{A}_t$ and feature map $\bm{\phi}$, and exploration bonus $\beta$.
    \STATE \textbf{Initialize:} \( \bm{\Lambda}_0 =\bm{I}_d \).
    \FOR{each round \( t = 1, 2, \dots, T \)}
        \STATE Compute 
        $\overline{\bm{\theta}}_{t-1}$ as in~(\ref{equ:overline-theta}) and $\bm{\widehat{\theta}}_{t-1}$ as in~(\ref{equ:hat-theta}).
        \STATE For each \( \ba \in \mathcal{A}_t \), compute UCB score $
        \mathrm{UCB}_t(\ba)$ as in~(\ref{equ:ucb}).
        \STATE Select and play action:
        $\ba_t = \arg\max_{\ba \in \mathcal{A}_t} \mathrm{UCB}_t(\ba)$ and observe reward \( r_t \).
        \STATE 
        Update $\bm{\Lambda}_t = \bm{\Lambda}_{t-1} + \bphi(\bx_t, \ba_t) \bphi(\bx_t, \ba_t)^\top$.
    \ENDFOR
\end{algorithmic}
\label{alg:linucb}
\end{algorithm}


\section{Main results}
\label{sec:main-result}

In this section, we develop an asymptotic theory for LinUCB without imposing the aforementioned stability assumption. We begin by stating several mild assumptions about our model and a few notation in Section~\ref{sec:ucb-unconstrained}. Our main results are provided in Section~\ref{sec:asymptotic-noramlity}, 
followed by a few remarks and implications. We present the main proof strategies of Theorem \ref{thm:ucb-unconstrained} in Section \ref{sec:technical-overview}. 


\subsection{Assumptions}
\label{sec:ucb-unconstrained}



In this work, we analyze a canonical and broadly applicable regime where at every round \(t\), the learner may choose any vector inside the unit ball of \(\mathbb{R}^d\) ($d\geq 2$) as the action. Formally, 

\begin{aspt}[Unconstrained action set on unit ball]\label{aspt:unconstrained}
    We assume that the set of feature maps is given by \( \bm{\Phi}_t = \mathcal{B}^{d} \) ($d\geq 2$) for all \( t \in [T] \), where \( \mathcal{B}^{d} \) denotes the unit ball in \( \mathbb{R}^d \).
\end{aspt}

The unit-ball action set allows the richest possible exploration directions. This stands in sharp contrast to the multi-armed bandit setting~\citep{khamaru2024inference}, where the learner is limited to exploring a finite number of pre-specified directions. 
Under Assumption~\ref{aspt:unconstrained}, the feature map is time-invariant with image \(\bm{\Phi}_t \equiv \mathcal{B}^d\), the unit ball in \(\mathbb{R}^d\). Consequently, choosing \(\bm a_t\) given \(\bm x_t\) is equivalent to selecting a point in \(\mathcal{B}^d\). We therefore take \(\mathcal{B}^d\) as the action set and, by slight abuse of notation, we identify each action with its feature vector, and write 
\(\bm a_t \equiv \bphi(\bm x_t,\bm a_t)\) throughout.

\medskip
We further assume the noise sequence $\{\epsilon_t\}_{t=1}^T$ (as in~(\ref{equ:linear})) is sub-Gaussian with parameter $\sigma.$ More concretely, we assume: 


\begin{aspt}[Sub-Gaussian noise]\label{aspt:subgaussian}
Let \((\mathcal F_t)_{t=0}^{T}\) denote the natural filtration generated by Algorithm \ref{alg:linucb}; that is,
\(\mathcal F_t := \sigma(\{\bm a_1,r_1,\ldots,\bm a_t,r_t\})\). The noise \((\epsilon_t)_{t=1}^T\) is an \((\mathcal F_t)\)-adapted martingale difference sequence with conditional mean zero and conditional variance \(\sigma^2\), and is conditionally sub-Gaussian with variance proxy \(\sigma^2\): for each \(t\ge1\) and all \(\lambda\in\mathbb{R}\),
\[
\mathbb{E}[\epsilon_t|\mathcal{F}_{t-1}] = 0,\quad \mathrm{Var}[\epsilon_t|\mathcal{F}_{t-1}] = \sigma^2,\quad 
\mathbb{E}\!\left[\exp(\lambda \epsilon_t)\mid \mathcal F_{t-1}\right]\le \exp\!\left(\frac{1}{2}\sigma^2\lambda^2\right).
\]
\end{aspt}

The sub-Gaussian property ensures the concentration of self-normalized martingales, which we shall leverage to control the cumulative effect of noise under adaptive sampling (see Theorem \ref{thm:noise-bound}). 
In fact, this assumption is crucial for 
establishing the asymptotic normality under adaptive sampling, whereas  for non-adaptive algorithms, asymptotic normality often follows from much weaker conditions (e.g., finite moment conditions, etc.). 

\subsection{Main theorems: asymptotic normality and confidence set}
\label{sec:asymptotic-noramlity}


We now state our main result on the asymptotic behavior of the LinUCB estimator \(\widehat{\bm\theta}_{T}\). Since \(\bm\theta^\star\) has unit norm, we project and obtain \(\widehat{\bm\theta}_{T}\) on the unit sphere \(\mathcal{S}^{d-1}\). Consequently, it cannot admit a nondegenerate \(d\)-dimensional limit in \(\mathbb{R}^d\). Instead, its first-order fluctuations are confined to the orthogonal linear subspace of \(\bm\theta^\star\), denoted as \((\bm\theta^\star)^{\perp}\). Accordingly, we characterize its limit distribution after projecting onto \((\bm\theta^\star)^{\perp}\).

\begin{thm}[Asymptotic normality for LinUCB]\label{thm:ucb-unconstrained}
Under Assumptions \ref{aspt:unconstrained}--\ref{aspt:subgaussian}, fix any matrix \(\bm U\in\mathbb R^{d\times(d-1)}\) with orthonormal columns orthogonal to \(\bm\theta^{\star}\); equivalently, let $\bm{Q} = (\bm{\theta}^{\star}, \bm{U})$, then $\bm{Q}^{\top}\bm{Q} = \bm{I}_d$. With $\beta\gg d^2(\sigma\sqrt{d+\log\log T}+1)$ and $\beta = O(\mathrm{poly}\log T)$,
the estimator \(\widehat{\bm\theta}_{T}\) in (\ref{equ:hat-theta}) satisfies Central Limit Theorem
\begin{align}
     \left(\frac{2\beta^{2}T}{d+1}\right)^{1/4}{\bm U}^{\top}\left(\widehat{\bm{\theta}}_{T} - \bm{\theta}^{\star}\right)\to\mathcal{N}(0,\sigma^2\bm{I}_{d-1}),
     \qquad \text{as } T\to \infty.
\end{align}
\end{thm}

 In Theorem \ref{thm:ucb-unconstrained}, the columns of matrix 
\(\bm U=(\bm u_1,\ldots,\bm u_{d-1})\) form an orthonormal basis of \((\bm\theta^\star)^{\perp}\); hence
\[
\bm U^{\top}\widehat{\bm\theta}_T
=\bigl(\bm u_1^{\top}\widehat{\bm\theta}_T,\ldots,\bm u_{d-1}^{\top}\widehat{\bm\theta}_T\bigr)
\]
gives the coordinates of the orthogonal projection of $\widehat{\bm\theta}_T$ onto \((\bm\theta^\star)^{\perp}\), namely \(\bm U\bm U^{\top}\widehat{\bm\theta}_T\) . 
Since the limiting covariance is \(\sigma^2\bm I_{d-1}\), this result is invariant to any choice of \(\bm U\): the estimator converges at the same rate in every direction of \((\bm\theta^\star)^{\perp}\). While we state the result in terms of \(\bm U^{\top}(\widehat{\bm\theta}_T-\bm\theta^\star)\), it is effectively an asymptotic normality statement for \(\widehat{\bm\theta}_T\) as well, since in a neighborhood of \(\bm\theta^\star\) on \(\mathcal{S}^{d-1}\), the mapping \(\bm\theta \mapsto \bm U^{\top}\bm\theta\) gives a one-to-one (indeed, near isometric, see (\ref{equ:near-isometric}) below) local reparameterization.

Building on this result, we construct an asymptotic \((1-\delta)\) confidence set for \(\bm{\theta}^{\star}\).
A naïve attempt would invert the limit law of \(\bm U^{\top}(\widehat{\bm\theta}_T-\bm\theta^\star)\), leading to a set defined by \(\|\bm U^{\top}(\widehat{\bm\theta}_T-\bm\theta)\|_2^2\); this is infeasible because \(\bm U\) depends on the unknown \(\bm\theta^\star\). Instead, we work directly with the Euclidean distance \(\|\widehat{\bm{\theta}}_{T}-\bm{\theta}\|_2^2\).
Using a local expansion around \(\bm{\theta}^{\star}\), we have
\begin{align}
\label{equ:near-isometric}
\big\|\widehat{\bm{\theta}}_{T}-\bm{\theta}^{\star}\big\|_2^2
= \Big[1 + O\Big(\big\|\bm{U}^{\top}(\widehat{\bm{\theta}}_{T}-\bm{\theta}^{\star})\big\|_2^2\Big)\Big]
\cdot \big\|\bm{U}^{\top}(\widehat{\bm{\theta}}_{T}-\bm{\theta}^{\star})\big\|_2^2
= \big[1+o_p(1)\big]\cdot \big\|\bm{U}^{\top}(\widehat{\bm{\theta}}_{T}-\bm{\theta}^{\star})\big\|_2^2,
\end{align}
where the last equality holds as  \(\|\bm{U}^{\top}(\widehat{\bm{\theta}}_{T}-\bm{\theta}^{\star})\|_2^2\to 0\), indicating these two criteria share the same asymptotic distribution. This allows us to construct the following confidence set of $\bm{\theta}^{\star}$.
\begin{cor}[Confidence set of LinUCB]
\label{cor:confidence}
Under Assumptions \ref{aspt:unconstrained}--\ref{aspt:subgaussian}, when $\beta\gg d^2(\sigma\sqrt{d+\log\log T}+1)$ and $\beta = O(\mathrm{poly}\log T)$, an asymptotic $(1-\delta)$ confidence set for $\bm{\theta}^{\star}$, based on $\widehat{\bm\theta}_T$, is given as
\begin{align}
    \mathcal{C}_{\delta} = \left\{\bm{\theta}\in\mathcal{S}^{d-1}: \left\Vert \widehat{\bm{\theta}}_{T} - \bm{\theta}\right\Vert_2^2\leq \widehat{\sigma}^2\sqrt{\frac{d+1}{2\beta^2 T}}\cdot \chi_{d-1,1-\delta}^2\right\},
\end{align}
where $\widehat{\sigma}^2$ is the estimated noise variance,  $\chi^{2}_{\,d-1,\,1-\delta}$ denotes the
$(1-\delta)$-quantile of the $\chi^{2}$ distribution with $d-1$ degrees of freedom.
\end{cor}
Here, we note that the noise variance $\sigma^2$ can be consistently estimated due to the consistency of $\bm{\widehat{\theta}}_T$. Furthermore, the confidence set is spherical rather than not an ellipsoid, because, as shall be seen in Theorem \ref{thm:covariance}, all non-leading eigenvalues converge to the same value asymptotically as $T\to\infty$. For moderate or small $T$, however, this convergence may be far from complete, so the empirical eigenvalues can still exhibit noticeable anisotropy. In such cases, the practitioner may prefer to capture the resulting ellipsoidal structure of the sampling distribution, as in the confidence set in~(\ref{eq:linucb-confidence-sets}). The following ellipsoidal confidence set provides an analogous asymptotic guarantee while accounting for this behavior:
\begin{align}
    \mathcal{C}_{\delta}^{\mathsf{ellipsoid}} = \left\{\bm{\theta}\in\mathcal{S}^{d-1}: \left\Vert \widehat{\bm{\theta}}_{T} - \bm{\theta}\right\Vert^2_{\bm{\Lambda}_T}\leq \widehat{\sigma}^2\chi_{d-1,1-\delta}^2\right\}.
\end{align}

Next, we highlight several implications and consequences of our main results.


\paragraph{Convergence rate slowdown and effective sample size.}
 With suitably chosen $\beta$, Theorem~\ref{thm:ucb-unconstrained} together with Corollary~\ref{cor:confidence} gives 
\[
\big\|\widehat{\bm{\theta}}_{T}-\bm{\theta}^\star\big\|_{2}
= \widetilde{\Theta}_p\!\left(T^{-1/4}\right).
\]
In contrast, the standard i.i.d.\ (parametric) rate is $\Theta_p\!\left(T^{-1/2}\right)$. 
The slowdown arises because in the settings when the action set is fixed over time, regret-driven adaptivity concentrates actions near the optimum, making the terminal design covariance $\boldsymbol{\Lambda}_T$ ill-conditioned (its non-leading eigenvalues grow sublinearly in $T$), thereby yielding a slower convergence rate of the estimator $\widehat{\bm\theta}_T$.


To quantify this effect across adaptive data-collection regimes, we define the \emph{effective sample size} as follows. If, for the terminal design covariance sequence $\{\boldsymbol{\Lambda}_T\}$, there exists a deterministic sequence $\{n_{\mathrm{eff},T}\}$ such that $n_{\mathrm{eff},T}^{-1}\,\lambda_{\min}(\boldsymbol{\Lambda}_T) \to 1$ with probability, then we refer to $n_{\mathrm{eff},T}$ as the effective sample size.
The examples below show that \(n_{\mathrm{eff},T}\) governs the estimator’s convergence rate, yielding \(\Theta_p\!\left(n_{\mathrm{eff},T}^{-1/2}\right)\).
\begin{itemize}
    \item \textbf{Linear regression with i.i.d.\ design.} Consider the classical linear regression setting with i.i.d.  design, in which the feature vector is sampled i.i.d. from some action set. 
    Let \(\bm{\Sigma}=\mathbb{E}[\ba_1\ba_1^{\top}]\) and assume \(\bm{\Sigma}\) is full rank. It is easily seen that $\{\bm{\Lambda}_T\}$ is stable: with  $\bm{\Sigma}_T^{-1}\bm{\Lambda}_T\to\bm{I}_d$ for $\bm{\Sigma}_T = T\bm{\Sigma}$. Hence the usual multivariate CLT applies: under standard regularity condition (e.g., mean-zero homoskedastic noise with variance \(\sigma^2\)),
$$
\sqrt{T}\bm{\Sigma}^{1/2}\bigl(\widehat{\bm{\theta}}_T-\bm{\theta}^\star\bigr)
\;\xrightarrow{\;d\;}\; \mathcal{N}\bigl(\bm{0},\sigma^2\bm{I}_d\bigr).
$$
In this case, the effective sample size is \( n_{\mathrm{eff},T} = T\lambda_{\min}(\bm{\Sigma}) \asymp T \),  admitting the classical \( T^{-1/2} \) convergence rate. 

    \item \textbf{Multi-armed bandit with UCB.} Consider the classical \(K\)-armed bandit over \(T\) rounds. In our notation this corresponds to features \(\bphi(\bx_t, \ba)=\bm{e}_a\in\mathbb{R}^K\) for arm \(a\). At round \(t\), using data up to time \(t{-}1\), the UCB score for arm \(a\) is
    \begin{align*}
        \mathrm{UCB}_{a,t} = \overline{X}_{a,t-1} + \frac{\beta}{\sqrt{n_{a,t-1}}},
    \end{align*}
    where \(\overline{X}_{a,t-1}\) is the empirical mean reward of arm \(a\) and \(n_{a,t-1}\) is the number of times arm \(a\) has been pulled up to \(t{-}1\); \(\beta>0\) controls exploration. The algorithm then selects the arm with the largest score.
   Under mild regularity conditions, the allocation vector is asymptotically deterministic \citep{khamaru2024inference,han2024ucb}: there exist deterministic counts \(\{n_{a,T}^\star\}_{a=1}^K\) such that \(n_{a,T}/n_{a,T}^\star \to 1\) as \(T\to\infty\). The limits are characterized by the “balancing’’ equations\footnote{The form in~(\ref{equ:multibandit-stable}) differs from \citet{khamaru2024inference} but is equivalent.}
\begin{align}\label{equ:multibandit-stable}
\mu_1 + \frac{\beta}{\sqrt{n_{1,T}^\star}}
= \mu_2 + \frac{\beta}{\sqrt{n_{2,T}^\star}}
= \cdots
= \mu_K + \frac{\beta}{\sqrt{n_{K,T}^\star}},
\qquad 
\sum_{a=1}^{K} n_{a,T}^\star = T,
\end{align}
where \(\mu_a\) is the mean reward of arm \(a\). 
Then, for each arm,
\[
\sqrt{n_{a,T}^\star}\,\bigl(\overline{X}_{a,T}-\mu_a\bigr)
\;\xrightarrow{\;d\;}\; \mathcal{N}\bigl(0,\sigma_a^2\bigr),
\]
where \(\sigma_a^2\) is the reward variance for arm \(a\). When all suboptimal arms have a fixed reward gap from the best arm, \citet{khamaru2024inference} show that
\(n_{K,T}^\star=\Theta(\beta^2)=\Theta(\log T)\) in the limit. Equivalently, with the terminal design matrix
\(\bLambda_T=\sum_{a=1}^K n_{a,T}\,\bm{e}_a\bm{e}_a^\top\) stabilizing to
\(\bSigma_T=\sum_{a=1}^K n_{a,T}^\star\,\bm{e}_a\bm{e}_a^\top\), the effective sample size is \(n_{\mathrm{eff},T}:=\min_a n_{a,T}^\star=n_{K,T}^\star=\Theta(\log T)\). Consequently, the estimation error for suboptimal arms decays at rate
\(\Theta\bigl(1/\sqrt{\log T}\bigr)=\Theta\bigl(1/\sqrt{n_{\mathrm{eff},T}}\bigr)\).

\end{itemize}
Returning to our setting where the action set is the unit ball---a rich action set that covers the entire feature space---the effective sample size is
\[
n_{\mathrm{eff},T}
= \sqrt{\frac{2\beta^{2} T}{d+1}}
= \widetilde{\Theta}(\sqrt{T}),
\]
Consequently, the estimator converges at a rate of  inverse square-root of \(n_{\mathrm{eff},T}\), i.e.,
\(\widetilde{\Theta}_p(T^{-1/4})\).
Relative to an i.i.d.\ design, the effective sample size is smaller for a fixed number of observations, leading to slower convergence. By contrast, it is substantially larger than in the multi-armed bandit UCB setting: with a continuous action set near the optimum, no single action dominates, so the learner continues to explore multiple directions, which increases regret ($\widetilde{\Theta}(\sqrt{T})$ for LinUCB versus $O(\sqrt{\log T})$ for UCB) yet yields better coverage (faster growth of $\lambda_{\min}(\boldsymbol{\Lambda}_T)$) and thus faster estimator convergence.

\paragraph{Comparison to the confidence set in~\citet{abbas2011improved}.} We compare the confidence set constructed in Corollary \ref{cor:confidence} with the confidence set in (\ref{eq:linucb-confidence-sets}) from \cite{abbas2011improved}. 
Our confidence set utilizes the standard Wald-type construction, while that of \cite{abbas2011improved} is based on martingale concentrations. Unlike their confidence set, ours explicitly utilizes the quantile of the chi-squared distribution instead of a $\log(1/\delta)$ concentration inequality-like term, thereby encompassing a stronger distributional result instead of only tail concentration.

Another key difference that we would like to point out is that our construction does not rely on the random, round-by-round empirical feature covariance accumulated by the algorithm -- but can do so if the user desires an ellipsoid confidence set. As such, we can provide a precise and deterministic (modulo randomness in $\widehat{\sigma}^2$) characterization of its diameter. Our confidence set is asymptotically tighter by a factor of $\sqrt{\log T}$. 
    

Most importantly, our confidence set is simply the Wald-type confidence set commonly employed within statistics for asymptotically normal models with i.i.d. data. As in \cite{khamaru2024inference}, this amounts to saying that we can treat the data collected by LinUCB when performing statistical inference as if it was i.i.d., just as with the UCB algorithm for multi-armed bandits.

\subsection{Technical overview}
\label{sec:technical-overview}
We next point out several key steps in the proof of Theorem~\ref{thm:ucb-unconstrained}. Although the final result is of asymptotic flavor, the argument rests on a sequence of non-asymptotic results. In adaptive data collection, where each action depends on past observations, such finite-sample controls are crucial: they stabilize the (effective) design covariance, which is the key ingredient for the asymptotics. We now outline the main steps and the technical challenges of the argument.

\paragraph{Uniform bound on ridge estimation error.} The first step of our proof is to establish a refined uniform bound for the ridge estimator when actions are selected by LinUCB. 
We first obtain a uniform \(\bm{\Lambda}_t\)-norm bound on the estimation error.

\begin{thm}[Uniform control of estimation error]
\label{thm:noise-bound}
    With probability $1-\frac{1}{\log T}$, the estimation error of ridge estimator $\overline{\bm \theta}_t$ scaled by cumulative covariance matrix $\bm{\Lambda}_t$ satisfies
    \begin{align}
        \max_{1\leq t\leq T}\left\Vert\overline{\bm \theta}_t-\bm{\theta}^{\star}\right\Vert_{\bm{\Lambda}_t}\lesssim\sigma\sqrt{d+\log\log T}+1.
    \end{align}
\end{thm}

\noindent With the \(\bm{\Lambda}_t\)-norm bound in Theorem~\ref{thm:noise-bound}, 
we immediately obtain a Euclidean norm control:
\begin{cor}
\label{cor:uniform-para-error}
    With probability $1-\frac{1}{\log T}$, the estimation error of ridge estimator $\overline{\bm \theta}_t$ can be uniformly upper bounded as
    \begin{align}
        \Vert \overline{\bm \theta}_t-\bm{\theta}^{\star}\Vert_2\lesssim\frac{\sigma\sqrt{d+\log\log T}+1}{\sqrt{\lambda_{t,d}}}.
    \end{align}
    Furthermore, one has the following upper bound for the projected estimator $\widehat{\bm\theta}_T$:
    \begin{align}
        \Vert \widehat{\bm \theta}_t-\bm{\theta}^{\star}\Vert_2\lesssim \frac{\sigma\sqrt{d+\log\log T}+1}{\sqrt{\lambda_{t,d}}}.
    \end{align}
\end{cor}

While a direct application of maximal concentration inequalities yields an upper bound of the estimation error of order \(O(\sqrt{\log T})\), Theorem~\ref{thm:noise-bound} and Corollary~\ref{cor:uniform-para-error} establish a tighter \(O(\sqrt{\log\log T})\) growth rate for the maximal normalized estimation error. This is achieved by leveraging the temporal correlation structure among the estimation errors.
When the exploration parameter \(\beta\) is chosen to grow faster than \(\log\log T\) (for instance, of order \(\log T\) as in \cite{abbas2011improved}), this result implies that the noise-induced estimation error is asymptotically dominated by the exploration bonus. This property constitutes a crucial component of our analysis, enabling a precise characterization of the asymptotic behavior of the UCB algorithm. 

A related result was established in Lemma 5.1 of \citet{khamaru2024inference},\footnote{That result was derived using the argument presented in \url{http://blog.wouterkoolen.info/QnD_LIL/post.html}.}
 though the analysis in our setting is considerably more involved due to the presence of a multi-dimensional noise term with a nonstationary sample covariance matrix.  To control this term, we introduce an exponential supermartingale and control a weighted aggregation of it, before establishing uniform concentration over the aggregated process. This step requires constructing a net---chosen so that the maximum over the full region is effectively captured by its maximum on the net, while accommodating the process's nonstationary nature.
In particular, we construct a global net that jointly covers all covariance matrices and noise terms up to time $T$, yet has only $O(\log T)$ cardinality. The construction leverages the rare-switching technique of \citet{abbas2011improved}, which is commonly used to control the growth of function classes in online reinforcement learning \citep{he2023nearly, sherman2024rateoptimalpolicyoptimizationlinear, tan2025actor}. However, we employ this idea differently to build an $O(\log T)$-sized collection of representative time indices which forms a covering net such that the associated covariance matrices collectively approximate all covariance matrices. The detailed proof is provided in Appendix \ref{app:noise-bound}.

\begin{figure}[t]
    \centering
    \includegraphics[width=0.6\linewidth]{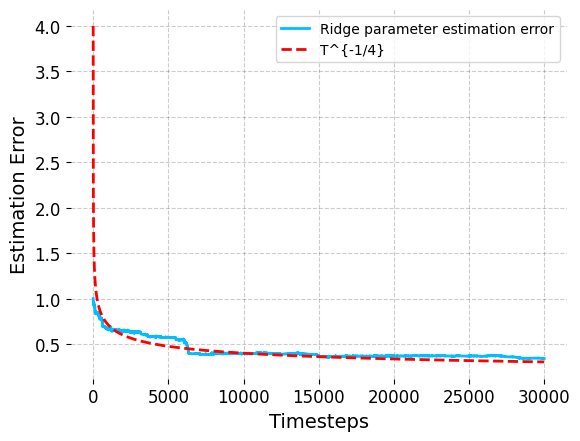}
    \caption{Estimation error of the parameter estimate obtained by ridge regression within the LinUCB algorithm. That is, we plot $\lVert \overline\btheta_t - \btheta^\star\rVert_2$ against timesteps. In this simulation, the action set is the unit ball and the optimal parameter is the first standard basis vector. We see that the estimation error decreases according to the $T^{-1/4}$ rate as predicted within Theorem \ref{thm:ucb-unconstrained}.}
    \label{fig:param-error}
\end{figure}

\paragraph{Characterization of $\bm{\Lambda}_T$.}
In the next step, we provide a careful characterization of the eigenstructure of the  terminal design \(\bm{\Lambda}_T\), which plays a critical role in developing Theorem~\ref{thm:ucb-unconstrained}. The analysis hinges on Theorem~\ref{thm:noise-bound}, which shows that, as \(T\) grows, the noise has minimal influence compared with the exploration bonus, rendering the former asymptotically negligible. 
Theorem~\ref{thm:covariance} formalizes the result, establishing (i) convergence of the top eigenvector and (ii) concentration of the non-leading eigenvalues around a deterministic limit.

\begin{thm}[Eigenstructure concentration of LinUCB]
\label{thm:covariance}
Under Assumptions \ref{aspt:unconstrained}--\ref{aspt:subgaussian}, let 
\(\{\lambda_{T,i}\}_{i=1}^{d}\) be the eigenvalues of \(\bm{\Lambda}_T\) ordered non-increasingly 
\(\lambda_{T,1}\ge \cdots \ge \lambda_{T,d}\), and let
\(\{\bm{v}_{T,i}\}_{i=1}^{d}\) be the corresponding eigenvectors. If $\beta\gg d^2(\sigma\sqrt{d+\log\log T}+1)$, then with probability $1-\frac{1}{\log T}$,

\begin{itemize}
\item \textbf{Alignment of the top eigenvector.}
The leading eigenvector $\bm{v}_{T,1}$ concentrates to the signal $\bm{\theta}^{\star}$,
\begin{align}
\label{equ:top-eig}
\|\bm{v}_{T,1}-\bm{\theta}^{\star}\|_{2}
\;\lesssim\;
\frac{\sigma\sqrt{d+\log\log T}+1}{\sqrt{\lambda_{T,d}}}\,.
\end{align}

\item \textbf{Concentration of non-leading eigenvalues.}
The non-leading eigenvalues concentrate uniformly to a deterministic limit: for any $i\geq 2$
\begin{align}
\label{equ:non-leading}
\lambda_{T,i}
\;=\;
\left[
1
+ O\!\left(
d\!\left(\frac{\beta^{8}}{T\sigma^{6}}\right)^{\!\frac{d+1}{d-1}}
+\left(\frac{\sigma\sqrt{d+\log\log T}+1}{\beta}\right)^{\!1/2}
\right)
\right]\sqrt{\frac{2\beta^2 T}{d+1}}.
\end{align}
Specifically, when $\beta = O(\mathrm{poly}\log T)$, we have 
\begin{align}
    \lambda_{T,i} = (1+o(1))\sqrt{\frac{2\beta^2 T}{d+1}}.
\end{align}
\end{itemize}   
\end{thm}
The proof of Theorem~\ref{thm:covariance} is the most technically involved part of our analysis. Our approach provides a precise account of how the cumulative covariance matrix \(\bm{\Lambda}_t\) $(t\in[T])$ evolves through four distinct phases; each phase calls for a different analysis. We outline the key arguments by phase in Section~\ref{sec:covariance} and present the full proof in Appendix~\ref{sec:proof-thm-2}.

\paragraph{Stabilizing projected covariance and establishing a CLT.} Lastly, we derive sharp characterizations for the \emph{projected} design covariance. Since the asymptotic variance of the projected estimator $\widehat{\bm\theta}_T$ is controlled by $\widetilde{\bm\Lambda}_T^{-1} = (\bm{U}^{\top}\bm{\Lambda}_{T}\bm{U})^{-1}$, it is sufficient to analyze the stability of projected sequence \(\{\widetilde{\bm\Lambda}_T\}\) rather than the full sequence \(\{\bm\Lambda_T\}\). For the diagonal matrix sequence
    \begin{align*}
        \widetilde{\bm{\Sigma}}_{T}\coloneqq \sqrt{\frac{2\beta^{2}T}{d+1}}\;\bm{I}_{d-1},
    \end{align*}
    we prove in Section~\ref{sec:proof-thm-1} that with probability $1-\frac{1}{\log T}$,
\begin{align}
\label{equ:convergence-covariance}
\bigl\Vert \widetilde{\bm{\Sigma}}_{T}^{-1}\widetilde{\bm{\Lambda}}_T-\bm{I}_{d-1}\bigr\Vert_{2}
\;\lesssim\;
d\!\left(\frac{\beta^{8}}{T\sigma^{6}}\right)^{\frac{d+1}{d-1}}
\;+\;
\left(\frac{\sigma\sqrt{\,d+\log\log T\,}+1}{\beta}\right)^{1/2}.
\end{align}

When $\beta\gg d^2(\sigma\sqrt{d+\log\log T}+1)$ and $\beta = O(\mathrm{poly}\log T)$, the right hand side of (\ref{equ:convergence-covariance}) vanishes, which suggests that \(\{\widetilde{\boldsymbol{\Lambda}}_{T}\}\) and \(\{\widetilde{\boldsymbol{\Sigma}}_{T}\}\) are asymptotically equivalent.
Intuitively, \(\beta\) trades off exploration and stability: 
if \(\beta\) is too small, the policy chases noise, the design fails to stabilize, and the estimator cannot achieve clean asymptotic normality; 
on the other hand, if \(\beta\) is too large, the policy over-explores, dispersing samples and slowing the concentration of information toward the desired limit.
The asymptotic normality of \(\widehat{\bm\theta}_T\) then follows directly from CLT with Lyapunov condition. A complete derivation of this part appears in Appendix~\ref{sec:proof-thm-1}.

\section{Non-asymptotic evolution of cumulative covariance}
\label{sec:covariance}

We now present a precise, non-asymptotic characterization of the cumulative covariance matrix $\bm{\Lambda}_t$. Its evolution over $t$ unfolds in four qualitatively distinct phases, governed by the changing influence of the bonus term; this structure, in turn, yields progressively sharper conclusions on the concentration of eigenvalues and eigenvectors as the cumulative covariance matrix $\bm{\Lambda}_t$ transitions across phases.
We formalize each phase in Propositions~\ref{prop:first-stage}--\ref{prop:fourth-stage}, which together provide a phase-wise characterization of $\bm{\Lambda}_t$. Specializing to $t=T$, Theorem~\ref{thm:covariance} follows directly from Propositions~\ref{prop:third-stage} and~\ref{prop:fourth-stage}.

\begin{figure}[t]
    \centering
    \includegraphics[width=0.7\linewidth]{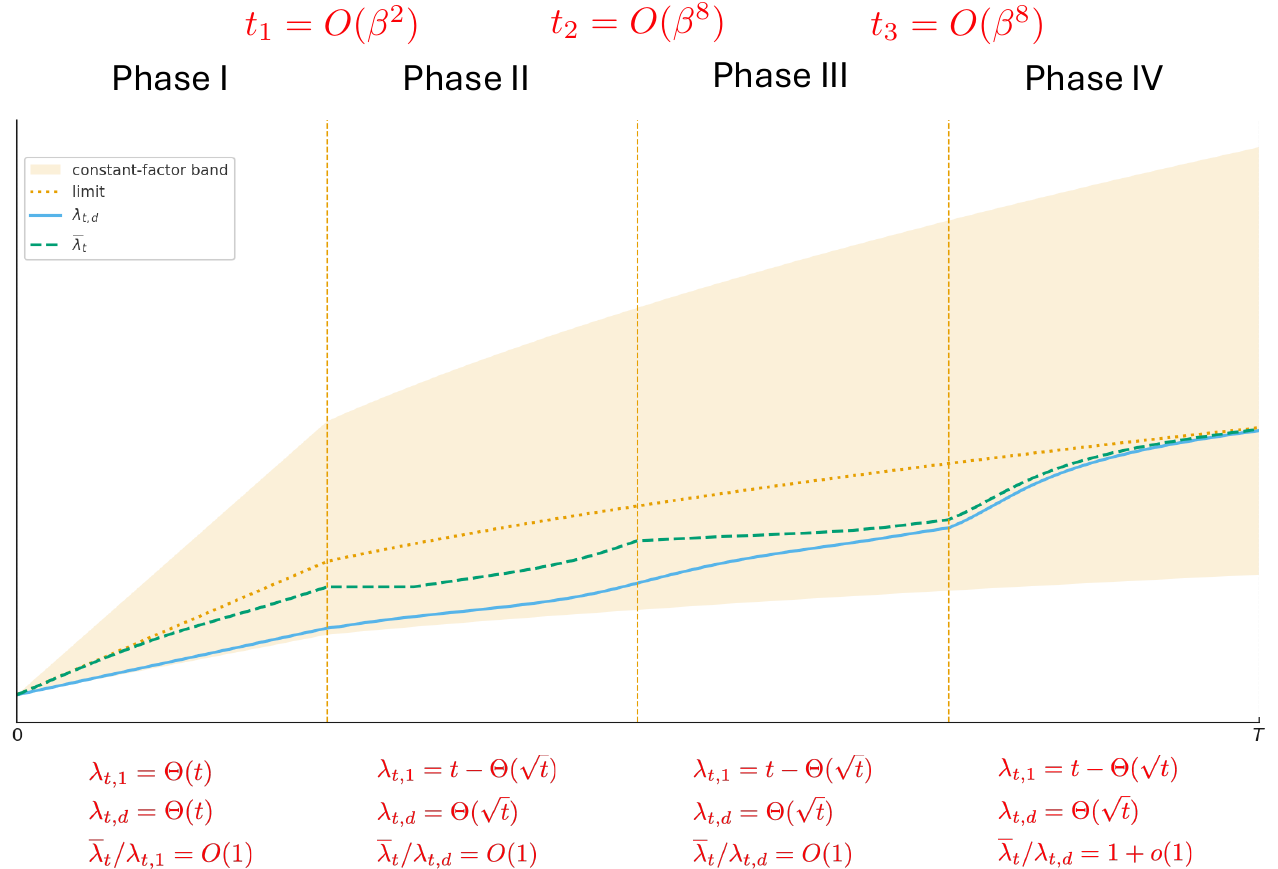}
    \caption{Growth of $\lambda_{t,d}$ and $\overline{\lambda}_t$. Throughout the entire process, these two quantities grow on the same order, falling in a constant-factor band of a deterministic growth benchmark $\lambda_t^{\star}$. When $t\geq t_2$, the minimum eigenvalue $\lambda_{t,d}$ concentrates close to the non-leading mean $\overline{\lambda}_t$, and both $\lambda_{t,d}$ and $\overline{\lambda}_t$ concentrates to a deterministic limit $\lambda_t^{\star}$ when $t=T$.}
    \label{fig:evolution-covariance}
\end{figure}

First, let us recall several notation. We use \(\lambda_{t,1}\ge \cdots \ge \lambda_{t,d}\) to denote the eigenvalues ranked non-increasingly, and use \(\overline{\lambda}_t\) as the average of the non-leading eigenvalues \(\sum_{i=2}^d \lambda_{t,i}/(d-1)\). Denote the leading eigenvector as \(\bm v_{t,1}\).
\begin{itemize}
    \item \textbf{Phase I (Initial exploration of all directions).} From the start, the minimum eigenvalue \(\lambda_{t,d}\) and  \(\overline{\lambda}_t\) grow at the same deterministic rate: there exists a deterministic sequence \(\{\lambda_t^\star\}_{t\ge 0}\) such that \(\lambda_{t,d} \asymp \overline{\lambda}_t \asymp \lambda_t^\star\) for all \(t\ge 0\). In particular, during Phase~I (i.e., for \(t \le t_1\)), we have \(\lambda_t^\star = \Theta(t)\).


    \item \textbf{Phase II (\(\Theta(\sqrt{t})\) growth of non-leading eigenvalues).} For \(t\ge t_1\), the comparability \(\lambda_{t,d} \asymp \overline{\lambda}_t \asymp \lambda_t^\star\) persists, now with \(\lambda_t^\star = \Theta(\sqrt{t})\). Meanwhile, the leading direction concentrates around \(\bm{\theta}^\star\): uniformly over \(t\), $\|\bm{v}_{t,1} - \bm{\theta}^\star\|_{\bm{\Lambda}_t} = O(\beta)$. Phase~II continues until \(t_2\), where a sharper bound for the leading eigenvector takes effect.


    
    \item \textbf{Phase III (Refined concentration of leading direction).} For \(t\ge t_2\), the leading eigenvector \(\bm v_{t,1}\) concentrates further around \(\bm\theta^\star\): the weighted error satisfies $\|\bm v_{t,1}-\bm\theta^\star\|_{\bm\Lambda_t} = O(\sqrt{\log\log T})$, uniformly in $t$, improving upon the \(O(\beta)\) bound from Phase~II. Phase~III lasts until \(t_3\), when stronger control of the non-leading eigenvalues becomes available.

    \item \textbf{Phase IV (Concentration of non-leading eigenvalues).} For all \(t\ge t_3\), we obtain the sharper concentration of non-leading eigenvalues: $\overline{\lambda}_t/\lambda_{t,d} = 1 + o(1)$, strengthening the earlier \(O(1)\) comparability. In addition, both \(\overline{\lambda}_t\) and \(\lambda_{t,d}\) start to concentrate towards the deterministic sequence \(\lambda_t^\star\); in particular, when \(t=T\), $\lambda_{T,d} = (1 + o(1))\lambda_T^\star$
    and $\overline{\lambda}_T = (1+o(1))\lambda_T^\star$ .

\end{itemize}
An illustration of the growth of the non-leading eigenvalues appears in Figure~\ref{fig:evolution-covariance}.

\begin{figure}[t]
    \centering
    \includegraphics[width=0.7\linewidth]{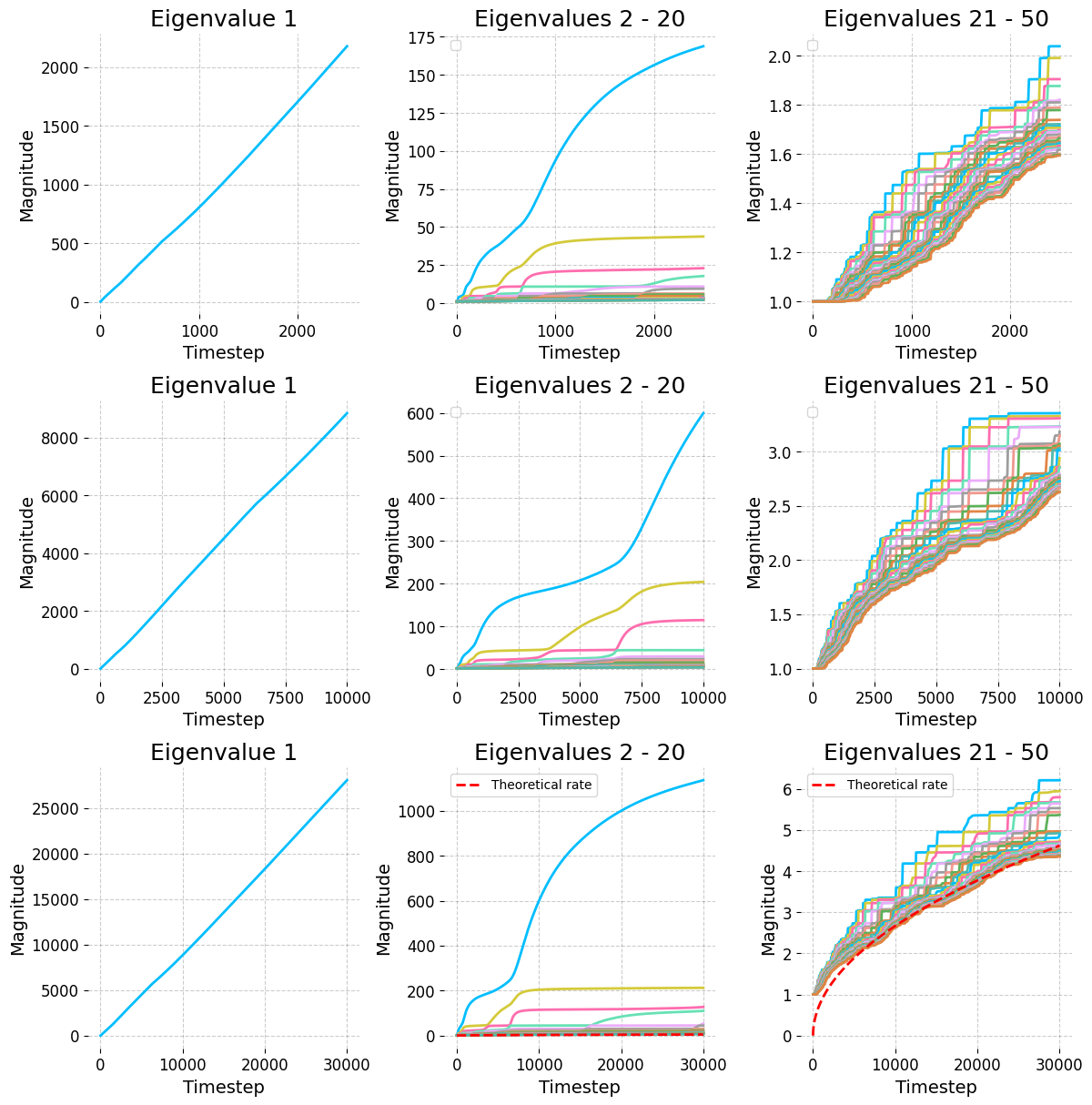}
    \caption{Rate of growth of the eigenvalues of the covariance matrix within a simulation of running LinUCB where the action set is the unit ball and the optimal parameter is the first standard basis vector. We see that the non-leading eigenvalues increase linearly at first as predicted by Proposition \ref{prop:first-stage}, before increasing on the order of $\sqrt{t}$ as predicted within Proposition \ref{prop:second-stage}. The dashed red line in the bottom row denotes the theoretical rate within Theorem \ref{thm:ucb-unconstrained}. Overall, the above simulation aligns well with our theory.}
    \label{fig:eig-growth}
\end{figure}



\paragraph{Phase I: Initial exploration of all directions.}
In the first phase, each direction is sampled only occasionally, so all eigenvalues of \(\bm{\Lambda}_t\) are small. Because the action \(\bm{a}_t\) is chosen by maximizing a UCB score based on \(\bm{\Lambda}_{t-1}\), the confidence bonus \(\beta\cdot (\bm{a}_t^\top \bm{\Lambda}_{t-1}^{-1}\bm{a}_t)^{1/2}\) dominates the predicted reward \(\langle \bm{a}_t, \mathcal{P}(\widehat{\bm{\theta}}_{t-1})\rangle\). This pushes actions toward the least-explored eigenspaces of \(\bm{\Lambda}_{t-1}\). We formalize the resulting spectral growth below. 

\begin{prop}
\label{prop:first-stage}

For LinUCB (Algorithm~\ref{alg:linucb}), there exists \(t_1=\Theta(\beta^2 d)\) such that

\begin{itemize}
    \item \textbf{(Linear growth of the minimum eigenvalue):}  
    For all $t\leq t_1$, the minimum eigenvalue \(\lambda_{t,d}\) of  \(\bm{\Lambda}_{t}\) grows at least linearly, i.e.,
    \begin{align}
        \lambda_{t,d} \asymp \frac{t}{d}.
    \end{align}


    \item \textbf{(Spectral gap between top two eigenvalues):}  
    With probability greater than $1-1/T$, the eigengap between the largest and the second largest eigenvalues at time \(t_1\) is lower bounded as
    \begin{align}
    \label{eqn:phase-1-eigengap}
        \lambda_{t_1,1} - \lambda_{t_1,2} \gtrsim t_1.
    \end{align}
\end{itemize}
\end{prop}

Proposition~\ref{prop:first-stage} certifies an initial \emph{exploration-dominated} phase in which every action is explored approximately at the same rate, and eigenvalues grow roughly uniformly, that \(\lambda_{t,d}\asymp t/d\). 
By time \(t_1\), a pronounced gap separates the top two eigenvalues, signaling the onset of an \emph{exploitation-dominated} stage in which one direction is revisited frequently. The scale \(t_1=\Theta(\beta^2 d)\) is natural: under \(\lambda_{t,d}\asymp t/d\), the typical UCB width behaves as \(\beta/\sqrt{\lambda_{t,d}}\asymp \beta/\sqrt{t/d}\) and becomes order one when \(t\asymp \beta^2 d\). After that, the process transitions into the next phase.

    

\paragraph{Phase II: $O(\sqrt{t})$ growth of non-leading eigenvalues.}
In the subsequent phase, once the minimum eigenvalue \(\lambda_{t,d}\) exceeds a fixed threshold, the estimator \(\widehat{\bm{\theta}}_t\) attains a non-trivial correlation with the true signal \(\bm{\theta}^{\star}\). From that point onward, action selection tilts toward the reward direction.
Consequently, the top eigenvalue \(\lambda_{t,1}\) absorbs most of the trace growth, while the non-leading eigenvalues grow with a slower diffusive rate. The precise growth rates and directional convergence are formalized next.

\begin{prop}
    \label{prop:second-stage}
    With probability greater than $1-1/T$, the following conditions hold for all $t_1\leq t\leq T$.
    \begin{itemize}
        \item \textbf{($O(\sqrt{t})$ growth of the non-leading eigenvalues)} The non-leading eigenvalues grow with comparable speed. Specifically, the mean of non-leading eigenvalues $\overline{\lambda}_t$, 
        and the minimum eigenvalue $\lambda_{t,d}$ satisfies
        \begin{align}
            \lambda_{t,d}\asymp\overline{\lambda}_t\asymp \beta\sqrt{\frac{t}{d}}.
        \end{align}
        \item \textbf{(Concentration of leading eigenvector)} The leading eigenvector $\bm{v}_{t,1}$ satisfies
        \begin{align}
            \Vert\bm{v}_{t,1}-\bm{\theta}^{\star}\Vert_2\lesssim\frac{\beta}{\sqrt{\lambda_{t,d}}}.
        \end{align}
    \end{itemize}
\end{prop}
The first part of this proposition guarantees that, beyond $t_1$, the non-leading spectrum is essentially flat and grows at the diffusive rate \(\Theta(\sqrt{t})\):
\(\lambda_{t,d}\asymp \overline{\lambda}_t \asymp \beta \sqrt{t/d}\).
Since each unit-norm action adds one to the trace, $\lambda_{t,1}
= t - \Theta \bigl(\beta \sqrt{td}\bigr)$, so \(\lambda_{t,1}\asymp t\) and the top–second eigengap is linear in \(t\).
The second part controls the distance between the leading eigenvector with $\thetastar$ via standard eigenvector perturbation bounds.
Notably, \(\beta/\sqrt{\lambda_{t,d}}\) coincides with the worst-direction UCB width $\beta\cdot (\bm{a}^\top\bm{\Lambda}_t^{-1}\bm{a})^{1/2}$,
so the leader’s misalignment is controlled by the same quantity that governs exploration.

\paragraph{Phase III: Refined concentration of the top eigenvector.}
In this phase, we provide a sharper concentration guarantee for the leading eigenvector \(\bm{v}_{t,1}\) of the cumulative covariance \(\bm{\Lambda}_t\).
Recall that, we have established a uniform high probability bound that \(\|\bm{v}_{t,1}-\bm{\theta}^\star\|_2 = O\big(\beta/\sqrt{\lambda_{t,d}}\big)\), matching the worst-direction exploration bonus, our goal for this stage is to obtain a finer result on this quantity.  

\begin{prop}
\label{prop:third-stage}
There exists $t_2 = O\!\left(\beta^{8}/(\sigma^{6} d^{2})\right)$ such that, with probability at least $1-\frac{1}{\log T}$, the following holds simultaneously for $t_3 \le t\leq T$: the leading eigenvector $\bm{v}_{t,1}$ satisfies

\begin{align}
\label{equ:top-refined-concentration}
\|\bm{v}_{t,1}-\bm{\theta}^\star\|_2 \;\lesssim\; \frac{\sigma\sqrt{d+\log\log T}+1}{\sqrt{\lambda_{t,d}}},
\qquad
\|\bm{v}_{t,1}-\widehat{\bm{\theta}}_t\|_2 \;\lesssim\; \frac{\sigma\sqrt{d+\log\log T}+1}{\sqrt{\lambda_{t,d}}}.
\end{align}
\end{prop}
Proposition~\ref{prop:third-stage} shows that, for \(t\ge t_2\), the misalignment of the leading direction $\|\bm{v}_{t,1}-\bm{\theta}^\star\|_2$ is of the same order as the maximum estimator error that shown in Corollary \ref{cor:uniform-para-error}, both decaying like \((\sigma\sqrt{d+\log\log T}+1)/\sqrt{\lambda_{t,d}}\). Since both the estimator $\widehat{\bm{\theta}}_t$ and the leading direction $\bm{v}_{t,1}$ concentrate tightly around the true signal $\bm{\theta}^\star$, they can be effectively treated as  ``quasi-deterministic’’ compared to the exploration bonus which scales as $O(\beta/\sqrt{\lambda_{t,d}})$. The evolution of non-leading eigenvalues, in the orthogonal space to $\thetastar$, depends mainly on the bonus term. 
This allows the non-leading eigenvalues to grow at comparable rates and concentrate accordingly—a mechanism that we formalize in the next phase.


\begin{figure}[t]
    \centering
    \includegraphics[width=0.6\linewidth]{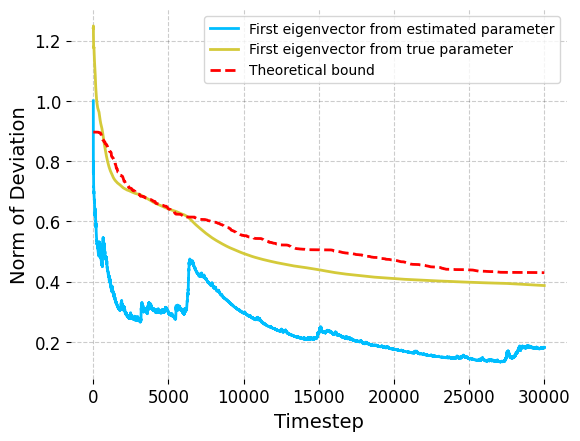}
    \caption{Concentration of the top eigenvector from the parameter estimate obtained through ridge regression and the true parameter. This is compared with the refined theoretical bound in Proposition \ref{prop:third-stage}.}
    \label{fig:top-eigenvector}
\end{figure}

\paragraph{Phase IV: Concentration of the non-leading eigenvalues.}
In this phase, we turn to the evolution of the \emph{non-leading} eigenvalues of the cumulative covariance.
Earlier phases established sharp control over the leading  eigenvector, which capture the dominant direction of variability.
Building on that, we now characterize the spectral structure on the subspace orthogonal to this top direction.

\begin{prop}
\label{prop:fourth-stage}
There exists $t_3 = O\left(\beta^{8}/\sigma^{6}\right)$ such that, with probability at least $1-\frac{1}{\log T}$, the following holds simultaneously for all $t_3 \le t\le T$: 
\begin{itemize}
    \item \textbf{(Near equality of non-leading eigenvalues)} The non-leading eigenvalues concentrate as follows:
    \[
\lambda_{t,2}
= \left[\,1+ O\!\left(\frac{d(\sigma \sqrt{d+\log\log T}+1)}{\beta}\right)\right]\lambda_{t,d}.
\]
    \item \textbf{(Deterministic benchmark and deviations)} The non-leading eigenvalues $\lambda_{t,i}$ ($i \ge 2$) satisfies: 
    \begin{align}
        \lambda_{t,i} = (1+\Delta_{t,i})\sqrt{\frac{2\beta^2t}{d+1}},
    \end{align}
    where $\Delta_{t,i}$ can be upper bounded as 
    \begin{align}
    \label{equ:Delta-t}
        |\Delta_{t,i}|\lesssim d\left( \frac{\beta^8}{t\sigma^6}\right)^{\frac{d+1}{d-1}} + \frac{(\sigma\sqrt{d+\log\log T}+1)^{1/2}}{\sqrt{\beta}},
    \end{align}
    whenever $\beta\gtrsim d^2(\sigma\sqrt{d+\log\log T}+1)$.
\end{itemize}
\end{prop}

\noindent
In words, Proposition~\ref{prop:fourth-stage} ensures that once $t\ge t_3$, the non-leading eigenvalues \(\lambda_{t,2},\ldots,\lambda_{t,d}\) are nearly equal and begin to concentrate around a deterministic value \(\bigl(2\beta^2 t/(d+1)\bigr)^{1/2}\). More specifically, they satisfy 
\[
  \frac{\lambda_{t,2}}{\lambda_{t,d}} = 1 + o(1),
\]
so the associated eigenspace is approximately \emph{isotropic}: restricted to the subspace orthogonal to the top eigenvector, the covariance is close to a scalar multiple of the identity. This near-isotropy indicates that, beyond the dominant signal direction, LinUCB explores the remaining directions at a nearly uniform rate.
As a result, in the end when \(t=T\),  non-leading eigenvalues satisfy 
\[
  \lambda_{T,i}
  = (1+\Delta_{T,i})\left(\frac{2\beta^2 T}{d+1}\right)^{1/2},
  \qquad 2\le i\le d,
\]
with \(\Delta_{T,i}\) controlled by~(\ref{equ:Delta-t}). This makes explicit the \(\sqrt{t}\)-scaling of the non-leading eigenvalues in directions orthogonal to the signal: the factor \(\beta^2\) reflects the choice of LinUCB exploration bonus (which inflates uncertainty to promote exploration), while the normalization \(2/(d+1)\) captures how the non-leading eigenvalues scale with dimension—when $d$ increases, exploration is spread across more dimensions, and consequently each non-leading eigenvalue decreases.


\begin{figure}[t]
    \centering
    \includegraphics[width=0.7\linewidth]{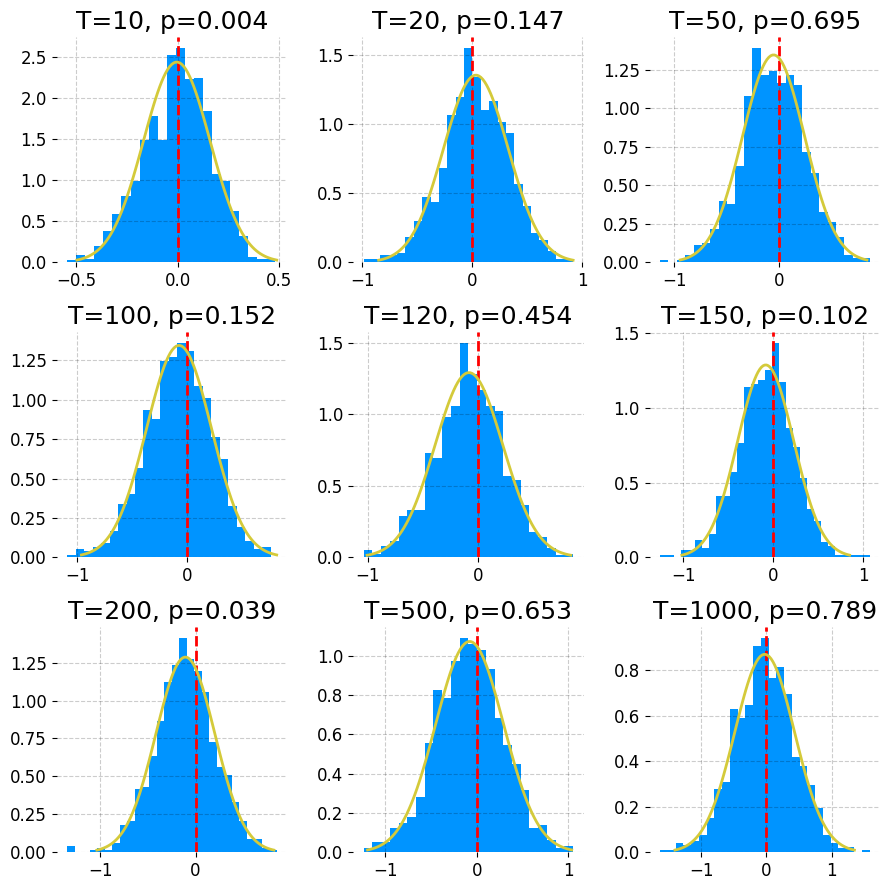}
    \caption{Asymptotic normality of the LinUCB algorithm in the setup of the MovieLens experiment within \citet{kausik2024leveragingofflinedatalinear}. For some random vector $\bm{u}$ on the unit ball, we plot $\sqrt{T} \cdot \bm{u}^\top ( \widehat{\btheta}_T - \btheta^\star)$ over 1000 independent trials, with KDE estimate overlaid as well as Shapiro-Wilk p-values provided as a test for non-normality. A finite number of possible movies to recommend are sampled at each round, and the learner has to recommend the best one available. As each possible action has good coverage on average, the minimum eigenvalue scales quickly on the order of $T$, and asymptotic normality is quickly achieved.}
    \label{fig:normality}
\end{figure}

\section{Conclusion and future work}
\label{sec:conclusion}

In this work, we characterize the asymptotic behavior of the LinUCB algorithm, and in doing so derive inference procedures that remain valid under adaptivity. Our main result shows the asymptotic normality of the terminal estimator of $\bm{\theta}^{\star}$:  after an explicit rescaling by $\left(2\beta^2 T/(d+1)\right)^{1/4}$, the estimation error projected onto the tangent space satisfies a central limit theorem with variance $\sigma^2 \bI_{d-1}$. This is accomplished through a thorough non-asymptotic characterization of the asymptotic behavior of the LinUCB feature covariance matrix, showing that it decomposes into a rank-one direction aligning with the true parameter and an isotropic bulk growing at a $\sqrt{T}$ rate. 

However, our results pertain to the case where the action set is the unit ball. This is partly by design -- rich action sets like these allow the learner to achieve good coverage over the feature space, but crucially, the learner does not need to (and in fact cannot) do so in order to achieve sublinear regret. With other choices of action sets, one should be able to achieve spiritually similar results, but the exact result will very much depend on whether these action sets ensure good 
coverage. For instance, within the MovieLens experiment within \cite{kausik2024leveragingofflinedatalinear}, a finite number of possible movies to recommend are sampled at each round, and the learner has to recommend the best one available. Here, as demonstrated in Figure \ref{fig:normality}, each possible action has good coverage on average, the minimum eigenvalue scales quickly on the order of $T$, and asymptotic normality is quickly achieved. We aim to explore this phenomenon in future work.

In addition, Berry-Esseen bounds characterizing distributional rates of convergence to the asymptotic distribution, as well as extensions to nonlinear function approximation methods and reinforcement learning \citep{wu2024statistical,wu2025uncertainty}, would also be welcome future directions to explore. 


\section*{Ackowledgement}

This work is supported in part by the NSF grants CCF-2106778, CCF-2418156 and CAREER award DMS-2143215.

\appendix 

\section{Technical preparations}
We begin with the technical preliminaries needed for the proofs of the main theorems. Appendix~\ref{sec:spectral} details notation on the spectral decomposition of the design covariance, and Appendix~\ref{sec:technical-lemma} lists auxiliary technical lemmas.
\subsection{Spectral decompositions}
\label{sec:spectral}
Throughout the proof, to analyze the spectral evolution of the design covariance matrix $\bm{\Lambda}_t$, we decompose $\bm{\Lambda}_t$ into its eigenvalues and eigenvectors. Moreover, we also decompose decompose the action vector $\bm{a}_t$ and the estimated signal $\bm{\widehat{\theta}}_t$ onto the orthogonal basis formed by the eigenvectors of $\bm{\Lambda}_t$.

\paragraph{Spectral decomposition $\bm{\Lambda}_t$.}  
As covariance matrix $\bm{\Lambda}_t$ is symmetric and positive semi-definite, it admits a spectral (or eigenvalue) decomposition of the form:
\begin{align}\label{equ:design-decom}
\bm{\Lambda}_t = \sum_{i=1}^d \lambda_{t,i} \bm{v}_{t,i} \bm{v}_{t,i}^\top,
\end{align}
where $\lambda_{t,1} \geq \lambda_{t,2} \geq \cdots \geq \lambda_{t,d}$ are the eigenvalues of $\bm{\Lambda}_t$ arranged in non-increasing order, and $\bm{v}_{t,1}, \ldots, \bm{v}_{t,d} \in \mathbb{R}^d$ are the corresponding eigenvectors. These eigenvectors satisfy the orthogonality condition $\bm{v}_{t,i}^\top \bm{v}_{t,j} = \delta_{ij}$, where $\delta_{ij}$ is the Kronecker delta. As a result, the collection $\{\bm{v}_{t,1}, \ldots, \bm{v}_{t,d}\}$ forms an orthonormal basis of $\mathbb{R}^d$, aligned with the principal directions of the covariance structure at time $t$.

\paragraph{Decomposition of $\bm{a}_t$ and $\bm{\widehat{\theta}}_t$.}  To facilitate component-wise analysis, we further express the vectors $\bm{a}_t$, $\widehat{\bm{\theta}}_t$ in terms of the eigenbasis $\{\bm{v}_{t,1}, \ldots, \bm{v}_{t,d}\}$ derived from the spectral decomposition of $\bm{\Lambda}_t$. That is,
\begin{align}
\label{equ:main-decom}
    \bm{a}_t = \sum_{i=1}^{d} \kappa_{t,i} \bm{v}_{t,i}, \quad
    \widehat{\bm{\theta}}_t = \sum_{i=1}^{d} \nu_{t,i} \bm{v}_{t,i},
\end{align}
where the coefficients $\kappa_{t,i}, \nu_{t,i}\in \mathbb{R}$ represent the projections of $\bm{a}_t$, $\widehat{\bm{\theta}}_t$ onto the $i$-th eigenvector $\bm{v}_{t,i}$, respectively. These coefficients are explicitly given by the inner products $\kappa_{t,i} = \bm{v}_{t,i}^\top \bm{a}_t$, $\nu_{t,i} = \bm{v}_{t,i}^\top \widehat{\bm{\theta}}_t$. With this decomposition, throughout the analysis, the vectors $\bm{a}_t$ and $\widehat{\bm{\theta}}_t$ can be expressed in terms of their coordinates with respect to the orthonormal basis $\{\bm{v}_{t,1}, \ldots, \bm{v}_{t,d}\}$.

We further decompose the action vector $\bm{a}_t$ into its component along the estimated signal and its orthogonal complement:
\begin{align}
\label{equ:a-decompose}
\bm{a}_t = \alpha_t \widehat{\bm{\theta}}_t + \bm{\xi}_t,
\end{align}
where $\widehat{\bm\theta}_t$ is the ridge estimator projected onto the unit sphere, $\bm\xi_t \in \mathbb{R}^d$ satisfies $\bm\xi_t^\top \widehat{\bm\theta}_t = 0$, and $\alpha_t \in \mathbb{R}$ is a scalar. Under Assumption~\ref{aspt:unconstrained}, the LinUCB maximizer over the unit ball lies on the unit sphere—indeed, the radial projection $\mathcal P(\bm a)$ of any interior point weakly increases $\mathrm{UCB}_t(\bm a)$ defined in (\ref{equ:ucb})—so $\|\bm a_t\|_2=1$, which implies $\alpha_t^2 + \Vert\bm{\xi}_t\Vert^2 = 1$. Intuitively, $\alpha_t$ measures alignment with the estimated signal, while $\bm\xi_t$ collects directions orthogonal to $\widehat{\bm\theta}_t$. For fixed $(\bm a_t,\widehat{\bm\theta}_t)$, this orthogonal decomposition is unique.



\subsection{Auxiliary lemmas}
\label{sec:technical-lemma}
We now establish a series of technical lemmas that will be employed in the proof of our results. The proofs of the lemmas stated in this section are deferred to Appendix \ref{sec:appendix-auxiliary}.

\paragraph{Characterizing action vector $\bm{a}_t$.}
Equipped with the spectral decompositions from Section~\ref{sec:spectral}, we characterize the LinUCB action $\bm a_t$ in the orthonormal basis $\{\bm v_{t,1},\ldots,\bm v_{t,d}\}$, which is key to tracking the evolution of the design covariance $\bm\Lambda_t$. Our goal is to obtain a closed-form description of the coefficients $\kappa_{t,i}$. We begin with an equivalent representation of $\bm a_t$.

\begin{lem}[An equivalent representation of $\bm{a}_t$]
\label{lem:a-t-characterization}
The LinUCB action admits the following representation:
\begin{align}
\label{equ:w-t-characterize}
    \bm{w}_t &= \mathrm{arg}\max_{\Vert \bm{w}\Vert_2 =1}\left\Vert\bm{\widehat{\theta}}_t + \beta\cdot \bm{\Lambda}_t^{-1/2} \bm{w}\right\Vert_2,\\
\label{equ:a-t-characterize}
    \bm{a}_t &= \mathcal{P}\left(\bm{\widehat{\theta}}_t + \beta \bm{\Lambda}_t^{-1/2} \bm{w}_t\right),
\end{align}
where $\mathcal P:\mathbb R^d\to\mathcal S^{d-1}$ denotes the projection onto the unit sphere $\mathcal S^{d-1}$.
\end{lem}
Lemma~\ref{lem:a-t-characterization} shows that $\bm a_t$ is the projection of the sum of the current parameter estimate and an exploration shift. The shift $\beta\,\bm\Lambda_t^{-1/2}\bm w_t$ points toward under-explored (high-variance) directions, enabling a recursive description of the selected actions in terms of the projected ridge estimator $\widehat{\bm\theta}_t$ in (\ref{equ:hat-theta}).  Consequently, the coefficients $\kappa_{t,i}$ admit a closed form. Writing $\bm w_t=(w_{t,1},\ldots,w_{t,d})$, we have the equivalent optimization
\begin{align}
\label{equ:w-t-expression}
    \bm{w}_t = \mathrm{arg}\max_{\Vert\bm{w}\Vert_2 = 1}\sum_{i=1}^{d}\left(\nu_{t,i}+\frac{\beta w_{i}}{\sqrt{\lambda_{t,i}}}\right)^2,\quad \text{where }\sum_{i=1}^{d}\nu_{t,i}^2 = 1,
\end{align} 
and the coefficients $\kappa_{t,i}$ are given by
\begin{align}
\label{equ:kappa-t-expression}
    \kappa_{t,i} = \frac{\nu_{t,i}+\frac{\beta w_{t,i}}{\sqrt{\lambda_{t,i}}}}{\sqrt{\sum_{j=1}^{d}\left(\nu_{t,j}+\frac{\beta w_{t,j}}{\sqrt{\lambda_{t,j}}}\right)^2}}.
\end{align}

This, in turn, enables us to quantify how the current action $\bm a_t$ decomposes into the component aligned with the estimated signal $\widehat{\bm\theta}_t$ and the orthogonal component $\bm\xi_t$, which is also related to the $\lambda_{t,d}$, the minimum eigenvalue of $\bm{\Lambda}_t$.

\begin{lem}[Spectral decomposition of $\bm{a}_t$]
\label{lem:decom-a}
Suppose that 
\begin{align}\label{equ:concen-top}
    \Vert\bm{v}_{t,1}-\bm{\widehat{\theta}}_{t}\Vert_2\leq h_t,
\end{align}
 the spectral decompositions in (\ref{equ:main-decom}) and (\ref{equ:a-decompose}) satisfies
\begin{itemize}
  \item The decomposition of $\bm{\widehat{\theta}}_t$  on orthogonal basis formed by eigenvectors $\{\bm{v}_{t,1},\ldots,\bm{v}_{t,d}\}$ satisfies:
  \begin{align}
      \nu_{t,1}\geq 1-h_t^2,\quad \nu_{t,i}\leq h_t.
  \end{align}
  \item The decomposition of $\bm{a}_t$ on $\bm{\widehat{\theta}}_t$ satisfies:
  \begin{align}
      \alpha_t = 1-O\left(\frac{\beta^2}{\lambda_{t,d}}\right),\quad \Vert\bm{\xi}_t\Vert = O\left(\frac{\beta}{\sqrt{\lambda_{t,d}}}\right).
  \end{align}
  \item The decomposition of $\bm{a}_t$ on orthogonal basis formed by eigenvectors $\{\bm{v}_{t,1},\ldots,\bm{v}_{t,d}\}$ satisfies:
  \begin{align}
      \kappa_{t,1} = 1-O\left(h_t^2+\frac{\beta^2}{\lambda_{t,d}}\right),\quad \kappa_{t,i} = O\left(h_t + \frac{\beta}{\sqrt{\lambda_{t,d}}}\right), \text{ for any }i\geq 2.
  \end{align}
\end{itemize}

\end{lem}

In other words, we provide a quantitative characterization of the intuition that if the estimator $\widehat{\bm{\theta}}_t$ aligns well with $\bm{v}_{t,1}$, the leading eigenvector of $\bm{\Lambda}_{t}$, then the action $\bm{a}_t$  aligns well with both $\bm{\widehat{\theta}}_t$ and $\bm{v}_{t,1}$.

\paragraph{Fine grained characterization of $\bm{a}_t$.} The lemmas above provide a coarse characterization of $\bm a_t$. While these formulas are clean and simple, they are not sufficient to establish finer properties of $\bm a_t$—in particular, to explain why $\bm a_t$ drives the non-leading eigenvalues of $\bm\Lambda_t$ to concentrate. To address this, we recast (\ref{equ:w-t-expression}) as a constrained optimization problem. Fix a radius $c_0>0$, a scalar $\beta\in\mathbb R$, a signal vector $\bm\nu=(\nu_1,\ldots,\nu_d)\in\mathbb R^d$, and a spectrum $\bm\lambda=(\lambda_1,\ldots,\lambda_d)$ with $\lambda_i>0$ for all $i$. Define the objective
\[
    g(\bm{w}) \;:=\; \sum_{i=1}^d\left(\nu_i+\frac{\beta w_i}{\sqrt{\lambda_i}}\right)^{2},
\]
and the maximizer over the $\ell_2$–sphere of radius $c_0$,
\begin{align}\label{equ:opt-problem}
    \bm{w}^\star(c_0,\bm{\nu},\bm{\lambda})
    \;:=\; \arg\max_{\|\bm{w}\|_2=c_0}\; g(\bm{w}).
\end{align}
For bookkeeping, we also define the (scaled) coordinate contributions at the maximizer,
\begin{align}
\label{equ:kappa-star}
    \kappa_i^{\star}(c_0,\bm{\nu},\bm{\lambda})
    \;:=\; \frac{\bigl(\nu_i+\frac{\beta w_i^{\star}}{\sqrt{\lambda_i}}\bigr)^{2}}
    {\sqrt{\,\sum_{j=1}^d\bigl(\nu_j+\frac{\beta w_j^{\star}}{\sqrt{\lambda_j}}\bigr)^{2}}}
    \;=\; \frac{\bigl(\nu_i+\frac{\beta w_i^{\star}}{\sqrt{\lambda_i}}\bigr)^{2}}{\sqrt{g(\bm{w}^\star)}}.
\end{align}

Our first result establishes a lower bound on the projection of $\bm{a}_t$ onto the eigenspace associated with the ``small eigenvalues". To avoid the extreme case where the optimization is driven solely by the signal vector $\bm{\nu}$, we impose the following structural condition.

\begin{aspt}\label{aspt:condition-opt}
There exists a constant $c>0$ such that
\[
    \max_{1\le i\le d}\frac{\beta^{2}}{\lambda_i}
    \;\;\ge\;\; \frac{c}{c_0}\|\bm{\nu}\|_2^2 .
\]
\end{aspt}

Equivalently, $\beta^2/\lambda_{\min}\!\ge (c/c_0)\|\bm{\nu}\|_2^2$, so at least one rescaled coordinate (governed by $\beta^2/\lambda_i$) can compete with the signal energy $\|\bm{\nu}\|_2^2$. Hence the optimizer is influenced by the eigenstructure $\{\lambda_i\}$ rather than aligning with $\bm{\nu}$ alone. We then define for a fixed constant $c_1> 1$ the index set of relatively small eigenvalues
\[
    \mathcal{L} \;:=\; \bigl\{i:\, \lambda_i \le c_1\,\min_{1\le j\le d}\lambda_j \,\bigr\}.
\]
Then we have the following concentration property.

\begin{lem}\label{lem:opt-concentration}
For $\bm{\kappa}^{\star}$ defined in (\ref{equ:kappa-star}), under Assumption~\ref{aspt:condition-opt}, there exists a constant $C=C(c,c_1)>0$ such that
\[
    \sum_{i\in\mathcal{L}} \bigl(\kappa_i^\star\bigr)^2\;\ge\; C \cdot c_0^2.
\]
\end{lem}

This ``spectral concentration’’ means the optimal solution cannot spread its scaled mass arbitrarily; it must allocate a non-negligible portion to indices with small $\lambda_i$, where the factor $\beta/\sqrt{\lambda_i}$ enhances the coordinate-wise effect. The optimization balances two forces: the signal $\bm{\nu}$ and the spectral scaling $\beta/\sqrt{\lambda_i}$. Assumption~\ref{aspt:condition-opt} ensures the latter is sufficiently strong, and Lemma~\ref{lem:opt-concentration} shows the optimizer reflects this by concentrating on a subset of small-$\lambda$ coordinates. We will use this later to control early-stage eigenvalue growth, in particular to show that certain eigenvalue ratios remain uniformly bounded under Assumption~\ref{aspt:condition-opt}.



We then show another result related to the constrained optimization problem (\ref{equ:opt-problem}). Consider the modified problem with the canonical signal $\bm{\tilde{\nu}}=(1,0,\ldots,0)$. Then the solution $\bm{w}^{\star}$ is defined as  
\begin{align}\label{equ:opt-modified}
    \widetilde{\bm{w}}^{\star}(c_0,\bm{\lambda})
    \;:=\; \bm{w}^{\star}(c_0,c_0\bm{\tilde{\nu}},\bm{\lambda})
    \;=\; \arg\max_{\|\bm{w}\|_2=c_0}\;
    \Bigl(c_0+\frac{\beta w_1}{\sqrt{\lambda_1}}\Bigr)^{\!2}
    \;+\; \sum_{i=2}^d \frac{\beta^2 w_i^2}{\lambda_i}.
\end{align}
We compare the first coordinate of the optimizer in the general  case to this canonical instance.

\begin{lem}
\label{lem:opt-comparison}
Let $\bm{\nu}\in\mathbb{R}^d$ satisfy $\|\bm{\nu}\|_2=c_0$ and set
$\bm{w}^{\star}(c_0,\bm{\nu},\bm{\lambda})$ as in (\ref{equ:opt-problem}), and
$\widetilde{\bm{w}}^{\star}(c_0,\bm{\lambda})$ as in (\ref{equ:opt-modified}).
Then
\begin{align*}
    \bigl|w^{\star}_1(c_0,\bm{\nu},\bm{\lambda})\bigr|
    \;\le\;
    \bigl|\widetilde{w}^{\star}_1(c_0,\bm{\lambda})\bigr|.
\end{align*}
\end{lem}

This reduction is particularly convenient when translating coordinate bounds into statements about normalized contributions.
For instance, any bound on $|\widetilde{w}_1^\star|$ immediately limits how much the term
$\bigl(\nu_1+\beta w_1^\star/\sqrt{\lambda_1}\bigr)^2$ can dominate the objective, and hence lower bound the contributions of $\bm{a}_t$ on non-leading eigenvalue directions. In later sections, we will exploit this “canonical-to-general” transfer to lower bound the growth of non-leading eigenvalues.

\paragraph{Rank-one update of $\bm{\Lambda}_t$.} To relate the action vector $\bm{a}_t$ to the evolution of the design covariance $\bm{\Lambda}_t$, we present the following lemma, which characterizes how the eigenvalues and eigenvectors of a positive-definite matrix evolve under a rank-one perturbation.
\begin{lem}[Theorem 8.4.3, \cite{golub2013matrix}]\label{lem:rank-one-update}
   Let $\bm{A} \in \mathbb{R}^{d \times d}$ be a positive definite matrix with eigenvalues 
$\lambda_1 \geq \lambda_2 \geq \cdots \geq \lambda_d > 0$, and corresponding eigenvectors $\bm{v}_1, \ldots, \bm{v}_d$. Let $\bm{u} = \sum_{i=1}^{d} \alpha_i \bm{v}_i$
for some scalars $\alpha_1, \ldots, \alpha_d \in \mathbb{R}$. Define the rank-one updated matrix $\widetilde{\bm{A}} = \bm{A} + \bm{u} \bm{u}^{\top}$. 
\begin{enumerate}
    \item The eigenvalues $\widetilde{\lambda}_1, \ldots, \widetilde{\lambda}_d$ of $\widetilde{\bm{A}}$ are the solutions to the secular equation
$$
f(\lambda) = 1 + \sum_{i=1}^{d} \frac{\alpha_i^2}{\lambda_i - \lambda} = 0.
$$
    \item The eigenvector $\widetilde{\bm{v}}_{1},\ldots,\widetilde{\bm{v}}_{d}$ of $\widetilde{A}$ are unit vectors that satisfies
    $$\widetilde{\bm{v}}_i\propto \sum_{j=1}^{d}\frac{\alpha_i}{\lambda_j-\widetilde{\lambda}_i}\bm{v}_j.$$
\end{enumerate}
\end{lem}

A direct consequence of the rank-one update lemma is listed as follows, where we quantify the growth of the largest eigenvalue of $\bm{\Lambda}_t$.
\begin{lem}[Growth of the largest eigenvalue]\label{lem:growth-largest-eigenvalue}
    When $\Vert\bm{v}_{t,1}-\bm{\widehat{\theta}}_t\Vert\lesssim \beta/\sqrt{\lambda_{t,d}}$, the largest eigenvalue of $\bm{\Lambda}_t$ evolves according to the update rule
    $$\lambda_{t+1,1} =\lambda_{t,1} + \kappa_{t,1}^2 + O(t^{-1}) =  \lambda_{t,1}+1 -O(\beta^2/\lambda_{t,d}).$$
\end{lem}

\paragraph{Multivariate Martingale Lindeberg CLT.} A multivariate Lindeberg Central Limit Theorem (CLT) for triangular arrays says that sums of many small, row-wise independent random vectors converge in distribution to a multivariate normal—provided their covariances stabilize and no single term has too much mass in its tails. We present this result as a lemma as follows.
\begin{lem}
\label{lem:CLT}
Let \(\{\bm{X}_{n,k}\}\) (\(n\in\mathbb{N}\), \(1\le k\le m_n\)) be an array of \(\mathbb{R}^d\)-valued random vectors. For each \(n\in\mathbb{N}\), let \((\mathcal{F}_{n,k})_{1\le k\le m_n}\) be a filtration to which \(\{\bm{X}_{n,k}\}_{1\le k\le m_n}\) is adapted, and suppose the array satisfies the following properties.
    \begin{align*}
        \mathbb{E}[\bm{X}_{n,k}|\mathcal{F}_{n,k-1}] = 0,\quad \bm{V}_n = \sum_{k=1}^{m_n}\mathrm{Var}(\bm{X}_{n,k}|\mathcal{F}_{n,k-1})\to\bm{\Sigma},
    \end{align*}
    where $\bm{\Sigma}$ is a fixed, positive definite $d\times d$ matrix. Furthermore, the array  $\{\bm{X}_{n,k}\}$ satisfies Lyapunov condition, i.e. there exists $\delta>0$ such that 
    \begin{align*}
        \sum_{k=1}^{m_n}\mathbb{E}\left[\Vert\bm{X}_{n,k}\Vert^{2+\delta}|\mathcal{F}_{n,k-1}\right]\to 0.
    \end{align*}
    Define the row sum $\bm{S}_n = \sum_{k=1}^{m_n}\bm{X}_{n,k}$. Then we have 
    \begin{align*}
        \bm{S}_n\stackrel{d}{\longrightarrow}\mathcal{N}(0,\bm{\Sigma}).
    \end{align*}
\end{lem}
The multivariate CLT with the Lyapunov condition stated above can be established for example by applying the one-dimensional Lindeberg CLT (Theorem~27.3 in \cite{billingsley2013convergence}) for one-dimensional projection $\boldsymbol{\theta}^\top \mathbf{X}_{n,k}$ ($\bm{\theta}\in\mathbb{R}^d$) and then using the Cramér--Wold theorem to conclude convergence in distribution of the vector.

\section{Proof of Theorem \ref{thm:ucb-unconstrained}}
\label{sec:proof-thm-1}
We present the proof of Theorem \ref{thm:ucb-unconstrained}, based on the conclusion of Theorem \ref{thm:covariance}. Without loss of generality, throughout the proof, we assume that $\bm{\theta}^{\star} =\bm{e}_1$. 
For general $\bm{\theta}^{\star}$, the result follows from the same analysis.
\paragraph{Step 1: show that $\Vert\widetilde{\bm{\Sigma}}_T^{-1}\widetilde{\bm{\Lambda}}_T - \bm{I}_{d-1}\Vert_2\to 0$.} Recall in (\ref{equ:design-decom}), we decompose $\bm{\Lambda}_T$ as 
\begin{align*}
    \bm{\Lambda}_T = \lambda_{T,1}\bm{v}_{T,1}\bm{v}_{T,1}^\top + \sum_{i=2}^{d}\lambda_{T,i}\bm{v}_{T,i}\bm{v}_{T,i}^\top. 
\end{align*}
From Proposition \ref{prop:fourth-stage}, when $\beta\gg d^2(\sigma\sqrt{d+\log\log T}+1)$, the non-leading eigenvalues satisfy 
\begin{align*}
    \lambda_{T,i} = (1+\Delta_{T,i})\sqrt{\frac{2\beta^2T}{d+1}},\qquad \;\forall\;i\geq 2.
\end{align*}
with the size of $\Delta_{T,i}$ obeying~(\ref{equ:Delta-t}), and consequently, the leading eigenvalue $\lambda_{T,1}$ can be characterized as
\begin{align*}
    \lambda_{T,1} = \sum_{i=1}^{d}\lambda_{T,i} - \sum_{i=2}^{d}\lambda_{T,i} = T+d -\left(d-1 + \sum_{i=2}^{d}\Delta_{T,i}\right)\sqrt{\frac{2\beta^2 T}{d+1}}.
\end{align*}
As a result, we can express $\bm{\Lambda}_T$ as follows
\begin{align}
      \bm{\Lambda}_T & = \lambda_{T,1}\bm{v}_{T,1}\bm{v}_{T,1}^\top + \sum_{i=2}^{d}\lambda_{T,i}\bm{v}_{T,i}\bm{v}_{T,i}^\top\nonumber\\
    & = \lambda_{T,1}\bm{v}_{T,1}\bm{v}_{T,1}^\top + \sqrt{\frac{2\beta^2 T}{d+1}}\sum_{i=2}^{d}\bm{v}_{T,i}\bm{v}_{T,i}^\top + \sqrt{\frac{2\beta^2 T}{d+1}}\sum_{i=2}^{d}\Delta_{T,i}\bm{v}_{T,i}\bm{v}_{T,i}^\top\nonumber\\
    & = \sqrt{\frac{2\beta^2 T}{d+1}}\bm{I}_d + \left(\lambda_{T,1}-\sqrt{\frac{2\beta^2 T}{d+1}}\right)\bm{v}_{T,1}\bm{v}_{T,1}^\top + \sqrt{\frac{2\beta^2 T}{d+1}}\sum_{i=2}^{d}\Delta_{T,i}\bm{v}_{T,i}\bm{v}_{T,i}^\top,
\end{align}
Here, as we set
$$\widetilde{\bm\Sigma}_T =\sqrt{\frac{2\beta^2T}{d+1}}\bm{I}_{d-1},$$
we can write 
\begin{align}
\label{equ:diff-tilde-lambda-sigma}
     \widetilde{\bm\Lambda}_T -\widetilde{\bm\Sigma}_T =   \left(\lambda_{T,1}-\sqrt{\frac{2\beta^2 T}{d+1}}\right) \cdot \bm{U}^{\top}\bm{v}_{T,1}\bm{v}_{T,1}^\top\bm{U} + \sqrt{\frac{2\beta^2 T}{d+1}}\bm{U}^{\top}\left(\sum_{i=2}^{d}\Delta_{T,i}\bm{v}_{T,i}\bm{v}_{T,i}^\top\right)\bm{U}.
\end{align}
We proceed to calculate the following quantities
\begin{align}
\label{equ:leading-U}
\left\Vert\bm{U}^{\top}\bm{v}_{T,1}\bm{v}_{T,1}^\top\bm{U}\right\Vert_2 & = \Vert\bm{U}^{\top}\bm{v}_{T,1} \Vert_2^{2}\leq \Vert\bm{v}_{T,1}-\bm{\theta}^{\star}\Vert_2^{2}\lesssim\frac{(\sigma\sqrt{d+\log\log T}+1)^2}{\lambda_{T,d}},\\
\label{equ:delta-v}    \left\Vert\sum_{i=2}^{d}\Delta_{T,i}\bm{v}_{T,i}\bm{v}_{T,i}^\top\right\Vert_2& \lesssim \max_{2\leq i\leq d}|\Delta_{T,i}|\cdot \left\Vert\sum_{i=1}^{d}\bm{v}_{T,i}\bm{v}_{T,i}^\top\right\Vert_2 = \max_{2\leq i\leq d}|\Delta_{T,i}|,
\end{align}
where the first inequality of $(\ref{equ:leading-U})$ holds true as
\begin{align*}
    \Vert\bm{U}^{\top}\bm{v}_{T,1} \Vert_2^{2} \leq\left(\Vert\bm{U}^{\top}\bm{\theta}^{\star} \Vert_2+\Vert\bm{U}^{\top}(\bm{v}_{T,1}-\bm{\theta}^{\star})\Vert_2\right)^2 = \Vert\bm{U}^{\top}(\bm{v}_{T,1}-\bm{\theta}^{\star})\Vert_2^2\leq \Vert\bm{v}_{T,1}-\bm{\theta}^{\star}\Vert_2^{2}.
\end{align*}
Consequently, with (\ref{equ:leading-U}) and (\ref{equ:delta-v}), we can upper bound the norm of (\ref{equ:diff-tilde-lambda-sigma}) as 
\begin{align*}
    \Vert \widetilde{\bm\Lambda}_T -\widetilde{\bm\Sigma}_T \Vert_2\lesssim T\cdot \frac{(\sigma\sqrt{d+\log\log T}+1)^2}{\lambda_{T,d}} +  \sqrt{\frac{2\beta^2 T}{d+1}}\cdot \max_{2\leq i\leq T}|\Delta_{T,i}|,
\end{align*}
which leads to the following upper bound when  $\beta\gg d^2(\sigma\sqrt{d+\log\log T}+1)$
\begin{align}
\label{equ:bound-asymptotic-tilde-lambda}
    \left\Vert \widetilde{\bm{\Sigma}}_T^{-1}\widetilde{\bm{\Lambda}}_T - \bm{I}_{d-1}\right\Vert_2 
    & \leq \Vert\widetilde{\bm{\Sigma}}_T^{-1}\Vert_2\cdot\Vert\widetilde{\bm{\Lambda}}_T- \widetilde{\bm{\Sigma}}_T \Vert_2\nonumber\\
    & {\lesssim T \sqrt{\frac{d+1}{2\beta^2 T}}\cdot \frac{(\sigma\sqrt{d+\log\log T}+1)^2}{\lambda_{T,d}} + \max_{2\leq i\leq d}|\Delta_{T,i}|}\nonumber\\
    &\lesssim \frac{d(\sigma\sqrt{d+\log\log T}+1)^2}{\beta^2} + d\left( \frac{\beta^8}{T\sigma^6}\right)^{\frac{d+1}{d-1}} + \frac{(\sigma\sqrt{d+\log\log T}+1)^{1/2}}{\sqrt{\beta}}\nonumber\\
    & \lesssim d\left( \frac{\beta^8}{T\sigma^6}\right)^{\frac{d+1}{d-1}} + \frac{(\sigma\sqrt{d+\log\log T}+1)^{1/2}}{\sqrt{\beta}},
\end{align}
where the third inequality follows from (\ref{equ:Delta-t}) and the fact that $\lambda_{T,d} \asymp \sqrt{\frac{2\beta^2 T}{d+1}}$, which establishes that
\begin{align*}
    \left\Vert \widetilde{\bm{\Sigma}}_T^{-1}\widetilde{\bm{\Lambda}}_T - \bm{I}_{d-1}\right\Vert_2  = o(1),
\end{align*}
whenever $\beta\gg d^2(\sigma\sqrt{d+\log\log T}+1)$.

\paragraph{Step 2: Asymptotic normality of  $\widetilde{\bm{\Sigma}}_{T}^{1/2}\bm{U}^{\top}(\overline{\bm{\theta}}_{T}-\bm{\theta}^{\star})$.}
Before diving into the analysis for $\widehat{\bm\theta}_T$, we first deal with the asymptotic of $\bm{\overline{\theta}}_T$ (defined in (\ref{equ:overline-theta})), which does not require projection to the unit sphere. To this end, we decompose $\overline{\bm\theta}_{T}$ as 
\begin{align}
\label{equ:decom-overline-theta}
    \overline{\bm\theta}_{T} = \bm{\Lambda}_{T}^{-1}\left[(\bm{\Lambda}_{T}-\bm{I}_d)\bm{\theta}^{\star} +  \bm{\eta}_{T}\right] = \bm{\theta}^{\star} - \bm{\Lambda}_{T}^{-1}\bm{\theta}^{\star} +  \bm{\Lambda}_{T}^{-1}\bm{\eta}_{T},
\end{align}
where $\bm{\eta}_{T} = \sum_{t=1}^{T}\bm{a}_t\epsilon_t$. In terms of this decomposition, we write 
\begin{align}
\label{eqn:thetabar-asymp-tmp}
    \widetilde{\bm{\Sigma}}_{T}^{1/2}\bm{U}^{\top}(\overline{\bm{\theta}}_{T}-\bm{\theta}^{\star}) =  \widetilde{\bm{\Sigma}}_{T}^{1/2}\bm{U}^{\top}\bm{\Lambda}_{T}^{-1}\bm{\eta}_T -\widetilde{\bm{\Sigma}}_{T}^{1/2}\bm{U}^{\top} \bm{\theta}^{\star}.
\end{align}

We shall deal with these two terms separately. Let us start with the first term since it is more involved. Construct a diagonal matrix
\begin{align}
    \bm{\Sigma}_T =  \sqrt{\frac{2\beta^2 T}{d+1}}\bm{I}_d.  
\end{align}
Intuitively, $\bm{\Sigma}_T$ approximates $\bm{\Lambda}_T$, except for the first row and column. We decompose the first term as follows
\begin{align}
\widetilde{\bm{\Sigma}}_{T}^{1/2}\bm{U}^{\top}\bm{\Lambda}_{T}^{-1}\bm{\eta}_T = \widetilde{\bm{\Sigma}}_{T}^{1/2}\bm{U}^{\top}\bm{\Sigma}_{T}^{-1}\bm{\eta}_T + \widetilde{\bm{\Sigma}}_{T}^{1/2}\bm{U}^{\top}(\bm{\Lambda}_{T}^{-1}-\bm{\Sigma}_{T}^{-1})\bm{\eta}_T.
\end{align}

We will first derive the asymptotic of the first term above. 
Here we need to use the triangular array argument, which was stated in Lemma \ref{lem:CLT}. 
Formally, let $\{\bm{a}_{t,s}\}$ and $\{\epsilon_{t,s}\}$ $(t\in\mathbb{N}, s\leq t)$ be arrays of actions and noise when implementing LinUCB. Then we rewrite $\bm{\Lambda}_T$ and $\bm{\eta}_T$ in the following way
\begin{align}
    \bm{\Lambda}_T & = \bm{I}_d + \sum_{s=1}^{t}\bm{a}_{T,s}\bm{a}_{T,s}^{\top},\qquad \bm{\eta}_T = \sum_{s=1}^{t}\bm{a}_{T,s}\epsilon_{T,s}.
\end{align}
Let $\mathcal{F}_{T,s} = \sigma(\bm{a}_{T,1},\epsilon_{T,1},\ldots,\bm{a}_{T,s},\epsilon_{T,s})$ and define
\begin{align}
    \bm{X}_{T,s} = \widetilde{\bm{\Sigma}}_{T}^{1/2}\bm{U}^{\top}\bm{\Sigma}_{T}^{-1}\bm{a}_{T,s}\epsilon_{T,s}.
\end{align}
Consequently, $\sum_{s=1}^{T}\bm{X}_{T,s} = \widetilde{\bm{\Sigma}}_{T}^{1/2}\bm{U}^{\top}\bm{\Sigma}_{T}^{-1}\bm{\eta}_T$. In addition, $\bm{X}_{T,s}\in\mathcal{F}_{T,s}$, i.e. $\bm{X}_{T,s}$ is adapted to the filtration $\mathcal{F}_{T,s}$, and the mean and variance of $\bm{X}_{T,s}$ conditioned on $\mathcal{F}_{T,s-1}$ are given as 
\begin{align}
    \mathbb{E}[\bm{X}_{T,s}|\mathcal{F}_{T,s-1}] = 0,\quad \mathrm{Var}[\bm{X}_{T,s}|\mathcal{F}_{T,s-1}] = \sigma^2\widetilde{\bm{\Sigma}}_{T}^{1/2}\bm{U}^{\top}\bm{\Sigma}_{T}^{-1}\bm{a}_{T,s}\bm{a}_{T,s}^{\top}\bm{\Sigma}_{T}^{-1}\bm{U}\widetilde{\bm{\Sigma}}_{T}^{1/2}.
\end{align}
As a matter of fact, it holds that
\begin{align}
    \sum_{s=1}^{T} \mathrm{Var}[\bm{X}_{T,s}|\mathcal{F}_{T,s-1}] & = \sigma^2\widetilde{\bm{\Sigma}}_{T}^{1/2}\bm{U}^{\top}\bm{\Sigma}_{T}^{-1}(\bm{\Lambda}_T-\bm{I}_d)\bm{\Sigma}_{T}^{-1}\bm{U}\widetilde{\bm{\Sigma}}_{T}^{1/2}\nonumber\\
    & = \sigma^2(\bm{0}_{d-1}, \widetilde{\bm{\Sigma}}_{T}^{-1/2})(\bm{\Lambda}_T-\bm{I}_d)(\bm{0}_{d-1}, \widetilde{\bm{\Sigma}}_{T}^{-1/2})^{\top}\nonumber\\
    & = \sigma^2\sqrt{\frac{d+1}{2\beta^2 T}}\bm{U}^{\top} (\bm{\Lambda}_T-\bm{I}_d)\bm{U}.
\end{align}
Here, since
\begin{align}
\label{equ:claim-lambda-t}
    \sqrt{\frac{d+1}{2\beta^2 T}}\bm{U}^{\top}\bm{\Lambda}_T\bm{U} = \widetilde{\bm\Sigma}_T^{-1}\widetilde{\bm\Lambda}_T\longrightarrow \bm{I}_d,
\end{align}
which was shown in Step 1, we can conclude 
\begin{align}
    \sum_{s=1}^{T} \mathrm{Var}[\bm{X}_{T,s}|\mathcal{F}_{T,s-1}]\longrightarrow\sigma^2\bm{I}_d.
\end{align}
To apply the Martingale Lindeberg CLT (Lemma~\ref{lem:CLT}), let us further verify the Lyapunov condition. For any $\delta>0$, we have
\begin{align*}
    \sum_{s=1}^{T}\mathbb{E}[\Vert\bm{X}_{t,s}\Vert^{2+\delta}|\mathcal{F}_{t,s-1}] = \max_{1\leq s\leq t }\mathbb{E}\left[\Vert\bm{X}_{t,s}\Vert^{\delta}|\mathcal{F}_{t,s-1}\right]\cdot  \sum_{s=1}^{t}\mathbb{E}\left[\Vert\bm{X}_{t,s}\Vert^{2}|\mathcal{F}_{t,s-1}\right]\leq \max_{1\leq s\leq t }\mathbb{E}\left[\Vert\bm{X}_{t,s}\Vert^{\delta}|\mathcal{F}_{t,s-1}\right].
\end{align*}
We upper bound $\mathbb{E}\left[\Vert\bm{X}_{t,s}\Vert^{\delta}|\mathcal{F}_{t,s-1}\right]$ by
\begin{align*}
    \mathbb{E}\left[\Vert\bm{X}_{T,s}\Vert^{\delta}|\mathcal{F}_{T,s-1}\right]\lesssim \left\Vert\widetilde{\bm{\Sigma}}_{T}^{1/2}\bm{U}^{\top}\bm{\Sigma}_{T}^{-1}\bm{a}_{T,s}\epsilon_{T,s}\right\Vert^{\delta}\lesssim \left( \sqrt{\frac{d+1}{2\beta^2 T}}\right)^{\delta/2}\longrightarrow 0,
\end{align*}
as $T\to\infty$, since $\lambda_{t,d} \to\infty$. Then, by virtue of the Martingale Lindeberg CLT (Lemma~\ref{lem:CLT}), we obtain 
\begin{align}
    \widetilde{\bm{\Sigma}}_{T}^{1/2}\bm{U}^{\top}\bm{\Sigma}_{T}^{-1}\bm{\eta}_T\longrightarrow\mathcal{N}(0,\sigma^2\bm{I}_d).
\end{align}

Let us then consider the term $ \widetilde{\bm{\Sigma}}_{T}^{1/2}\bm{U}^{\top}(\bm{\Lambda}_{T}^{-1}-\bm{\Sigma}_{T}^{-1})\bm{\eta}_T$. In view of the Cauchy-Schwarz inequality, the term of interest can be written as 
\begin{align}
    \Vert\bm{U}^{\top}(\bm{\Lambda}_{T}^{-1}-\bm{\Sigma}_{T}^{-1})\bm{\eta}_T\Vert_2 \leq  
    \left\Vert\bm{U}^{\top}(\bm{\Lambda}_{T}^{-1}-\bm{\Sigma}_{T}^{-1})\bm{\Lambda}_T^{1/2}\right\Vert_2\cdot \left\Vert\bm{\Lambda}_T^{-1/2}\bm{\eta}_T\right\Vert_2.
\end{align}
To bound the right hand side of the above inequality, we first note that
\begin{align*}
    \mathbb{E}\left[\Big\Vert\bm{\Lambda}_T^{-1/2}\bm{\eta}_T\Big\Vert_2^2\right] & = \sigma^2\,\mathrm{tr}\Big(\bm{\Lambda}_T^{-1/2}\Big(\sum_{s=1}^T \bm{a}_s\bm{a}_s^\top\Big)\bm{\Lambda}_T^{-1/2}\Big)\leq\sigma^2 d,
\end{align*}
which leads to $\Vert\bm{\Lambda}_T^{-1/2}\bm{\eta}_T\Vert_2 = O_p(1)$. 
It is therefore only left for us to control the quantity $\Vert\bm{U}^{\top}(\bm{\Lambda}_{T}^{-1}-\bm{\Sigma}_{T}^{-1})\bm{\Lambda}_T^{1/2}\Vert_2$.
Recalling the definition $\bm{\Sigma}_T^{-1} =  \sqrt{\frac{d+1}{2\beta^2 T}}\bm{I}_d$, we write 
\begin{align}
    \left\Vert\bm{U}^{\top}(\bm{\Lambda}_{T}^{-1}-\bm{\Sigma}_{T}^{-1})\bm{\Lambda}_T^{1/2}\right\Vert_2
    & = \left\Vert\bm{U}^{\top}\left(\bm{\Lambda}_{T}^{-1/2}-\sqrt{\frac{d+1}{2\beta^2 T}}\bm{\Lambda}_T^{1/2}\right)\right\Vert_2.
\end{align}
To control the right-hand side of the above equality, invoking the decomposition~(\ref{equ:design-decom}) yields that
\begin{align}
    \bm{\Lambda}_{T}^{-1/2}-\sqrt{\frac{d+1}{2\beta^2 T}}\bm{\Lambda}_T^{1/2}=\sum_{i=1}^{d}\left(\frac{1}{\sqrt{\lambda_{T,i}}} -\sqrt{\frac{d+1}{2\beta^2 T}}\sqrt{\lambda_{T,i}}\right)\bm{v}_{T,i}\bm{v}_{T,i}^{\top}.
\end{align}
For $i\geq 2$, notice that $\lambda_{T,i} = (1+o(1))\sqrt{(d+1)/(2\beta^2T)}$. As a result, we have
\begin{align*}
    \frac{1}{\sqrt{\lambda_{T,i}}} -\sqrt{\frac{d+1}{2\beta^2 T}}\sqrt{\lambda_{T,i}} = \frac{1}{\sqrt{\lambda_{T,i}}}\cdot\left(1-\sqrt{\frac{d+1}{2\beta^2 T}}\lambda_{T,i}\right) = o(1)\cdot\frac{1}{\sqrt{\lambda_{T,i}}}.
\end{align*}
For $i=1$, recall  
\begin{align*}
    \Vert\bm{U}^{\top}\bm{v}_{T,1}\Vert_2\leq \Vert\bm{v}_{T,1}-\bm{\theta}^{\star}\Vert_2\lesssim\frac{\sigma\sqrt{d+\log\log T}+1}{\sqrt{\lambda_{T,d}}},
\end{align*}
to arrive at 
\begin{align}
    \left\Vert\bm{U}^{\top}\left(\frac{1}{\sqrt{\lambda_{T,1}}} -\sqrt{\frac{d+1}{2\beta^2 T}}\cdot\sqrt{\lambda_{T,1}}\right)\bm{v}_{T,1}\bm{v}_{T,1}^{\top}\right\Vert_2& \lesssim \frac{\sigma\sqrt{d+\log\log T}+1}{\sqrt{\lambda_{T,d}}}\cdot \left(\frac{1}{\sqrt{T}} + \sqrt{T}\cdot \frac{\sqrt{d}}{\beta\sqrt{T}}\right)\nonumber\\
    & \lesssim \frac{\sigma\sqrt{d+\log\log T}+1}{\sqrt{\lambda_{T,d}}}\cdot \frac{\sqrt{d}}{\beta},
\end{align}
where the last inequality holds as $\beta = O\left(\mathrm{poly}\log T\right)\ll T$. Putting pieces together, we obtain the upper bound
\begin{align}
    \left\Vert\bm{U}^{\top}(\bm{\Lambda}_{T}^{-1}-\bm{\Sigma}_{T}^{-1})\bm{\Lambda}_T^{1/2}\right\Vert_2\lesssim\frac{1}{\sqrt{\lambda_{T,d}}}\left(d\cdot o(1) + \frac{\sqrt{d}(\sigma\sqrt{d+\log\log T}+1)}{\beta}\right) = \frac{1}{\sqrt{\lambda_{T,d}}}\cdot o(1),
\end{align}
as we set $\beta\gg d^2(\sigma\sqrt{d+\log\log T}+1)$. As a result, we conclude that
\begin{align}
    \left\Vert\widetilde{\bm{\Sigma}}_{T}^{1/2}\bm{U}^{\top}(\bm{\Lambda}_{T}^{-1}-\bm{\Sigma}_{T}^{-1})\bm{\eta}_T\right\Vert_2 = \left(\frac{2\beta^2 T}{d+1}\right)^{1/4}\cdot\frac{1}{\sqrt{\lambda_{T,d}}}\cdot o(1)\cdot O_p(1) = o_p(1).
\end{align}
In other words, the term $\Vert\widetilde{\bm{\Sigma}}_{T}^{1/2}\bm{U}^{\top}(\bm{\Lambda}_{T}^{-1}-\bm{\Sigma}_{T}^{-1})\bm{\eta}_T\Vert_2$ converges to $0$ in probability. As a result, one has
\begin{align}
\label{eqn:brahms}\widetilde{\bm{\Sigma}}_{T}^{1/2}\bm{U}^{\top}\bm{\Lambda}_{T}^{-1}\bm{\eta}_T\longrightarrow\mathcal{N}(0,\sigma^2\bm{I}_d).
\end{align}

For the second term, we note that 
\begin{align*}
    \left\Vert\widetilde{\bm{\Sigma}}_{T}^{1/2}\bm{U}^{\top}\bm{\Lambda}_{T}^{-1}\bm{\theta}^{\star}\right\Vert_2\leq \left\Vert\widetilde{\bm{\Sigma}}_{T}^{1/2}\bm{U}^{\top}\bm{\Lambda}_{T}^{-1}\right\Vert_2\leq \left(\frac{2\beta^2 T}{d+1}\right)^{1/4}\cdot\frac{1}{\lambda_{T,d}}\longrightarrow 0,
\end{align*}
as $T\to\infty$. Combining this with~(\ref{eqn:thetabar-asymp-tmp}) and~(\ref{eqn:brahms}), we arrive at 
\begin{align}
    \widetilde{\bm{\Sigma}}_{T}^{1/2}\bm{U}^{\top}(\overline{\bm{\theta}}_{T}-\bm{\theta}^{\star})\longrightarrow\mathcal{N}(0,\sigma^2\bm{I}_d).
\end{align}

\paragraph{Step 3: Asymptotic normality of  $\widetilde{\bm{\Sigma}}_{T}^{1/2}\bm{U}^{\top}(\widehat{\bm{\theta}}_{T}-\bm{\theta}^{\star})$. }
In view of the decomposition in~(\ref{equ:decom-overline-theta}), we bound the norm of $\overline{\bm{\theta}}_{T}$ by triangle's inequality as  
\begin{align}
   \left\vert \Vert\overline{\bm{\theta}}_{T} \Vert_2 - 1\right\vert\leq  \left\Vert\bm{\Lambda}_{T}^{-1}\bm{\theta}^{\star}\right\Vert_2+\left\Vert\bm{\Lambda}_{T}^{-1}\bm{\eta}_{T}\right\Vert_2\leq\frac{1}{\lambda_{T,d}} + \frac{1}{\sqrt{\lambda_{T,d}}}\left\Vert\bm{\Lambda}_{T}^{-1/2}\bm{\eta}_{T}\right\Vert_2 = O_p\left(\frac{1}{\sqrt{\lambda_{T,d}}}\right).
\end{align}
As a result, for the projection $\widehat{\bm{\theta}}_{T} = \overline{\bm{\theta}}_{T}/\|\overline{\bm{\theta}}_{T}\|_2$, it is easily seen that
\begin{align}
    \Vert\widehat{\bm{\theta}}_{T}-\bm{\theta}^{\star}\Vert\leq \Vert\overline{\bm{\theta}}_{T}-\bm{\theta}^{\star}\Vert =  O_p\left(\frac{1}{\sqrt{\lambda_{T,d}}}\right).
\end{align}
Furthermore, as $\widehat{\bm\theta}_T$ being the projection of $\overline{\bm\theta}_T$, we have
\begin{align*}
    (\bm{\theta}^{\star})^{\top}(\overline{\bm{\theta}}_{T}-\widehat{\bm{\theta}}_{T}) =   (\bm{\theta}^{\star})^{\top}\widehat{\bm{\theta}}_{T}\cdot \Vert\overline{\bm{\theta}}_{T}-\widehat{\bm{\theta}}_{T} \Vert_2.
\end{align*}
Observing the relation  
\begin{align*}
    \left\Vert(\bm{\theta}^{\star})^{\top}(\overline{\bm{\theta}}_{T}-\widehat{\bm{\theta}}_{T})\right\Vert_2^2 + \left\Vert\bm{U}^{\top}\left(\overline{\bm{\theta}}_{T}-\widehat{\bm{\theta}}_{T}\right)\right\Vert_2^2 = \Vert\overline{\bm{\theta}}_{T}-\widehat{\bm{\theta}}_{T}\Vert_2^2,
\end{align*}
leads to 
\begin{align}
    \left\Vert\bm{U}^{\top}\left(\overline{\bm{\theta}}_{T}-\widehat{\bm{\theta}}_{T}\right)\right\Vert_2 & = \sqrt{1-\left[(\bm{\theta}^{\star})^{\top}\widehat{\bm{\theta}}_{T}\right]^2}\cdot  \Vert\overline{\bm{\theta}}_{T}-\widehat{\bm{\theta}}_{T}\Vert_2 \lesssim \Vert\widehat{\bm{\theta}}_{T}-\bm{\theta}^{\star}\Vert\cdot \Vert\overline{\bm{\theta}}_{T}-\widehat{\bm{\theta}}_{T} \Vert_2\nonumber\\
    & = O_p\left(\frac{1}{\lambda_{T,d}}\right),
\end{align}
where the first inequality holds as 
\[
\left\|\widehat{\bm\theta}_T-\bm{\theta}^{\star}\right\|_2^2
=2\bigl(1-\langle \widehat{\bm\theta}_T,\bm{\theta}^{\star}\rangle\bigr)
\ge \bigl(1-\langle \widehat{\bm\theta}_T,\bm{\theta}^{\star}\rangle\bigr)\bigl(1+\langle \widehat{\bm\theta}_T,\bm{\theta}^{\star}\rangle\bigr)
=1-\langle \widehat{\bm\theta}_T,\bm{\theta}^{\star}\rangle^2.
\]
Therefore, we obtain the following result
\begin{align}
    \left(\frac{2\beta^{2}T}{d+1}\right)^{1/4}\bm{U}^{\top}\left(\overline{\bm{\theta}}_{T}-\widehat{\bm{\theta}}_{T}\right) = O_p\left(\frac{1}{\sqrt{\lambda_{T,d}}}\right).
\end{align}
Since $\lambda_{T,d}\to\infty$, combining this with the asymptotic result in Step 2, we conclude that 
\begin{align*}
     \left(\frac{2\beta^{2}T}{d+1}\right)^{1/4}{\bm U}^{\top}\left(\widehat{\bm{\theta}}_{T} - \bm{\theta}^{\star}\right)\to\mathcal{N}(0,\bm{I}_{d-1}).
\end{align*}

\section{Proof of Theorem \ref{thm:noise-bound}}
\label{app:noise-bound}


We present the full proof of Theorem \ref{thm:noise-bound} in this section. Although a similar result was established in Lemma 5.1 of \citet{khamaru2024inference} based on a previous result (see this \href{http://blog.wouterkoolen.info/QnD_LIL/post.html}{blog}), the proof in our setting is much more involved. This increased complexity arises from two key challenges: first, the noise term $\bm{\eta}_t=\sum_{s=1}^{t}\bm{a}_s\epsilon_s$ is multi-dimensional, and second, the sample covariance matrix $\bm{\Lambda}_t$ evolves in a non-stationary manner over time. To proceed, let us denote 
\begin{align}
    \bm{\widetilde{\theta}}_t = \mathbb{E}\left[\left(\sum_{s=1}^{t-1}\bm{a}_s\bm{a}_s^{T}+\bm{I}_d\right)^{-1}\left(\sum_{s=1}^{t-1}\bm{a}_s y_s\right)\right] = \bm{\Lambda}_{t-1}^{-1}(\bm{\Lambda}_{t-1}-\bm{I}_d)\bm{\theta}^{\star},
\end{align}
 be the expectation of ridge regression regression estimator $\overline{\bm\theta}_t$, and set $\bm{\eta}_t = \sum_{s=1}^{t}\bm{a}_s\epsilon_s$ be the corresponding noise part. Then we can rewrite $\overline{\bm\theta}_t$ as 
 \begin{align}
     \overline{\bm\theta}_t = \bm{\Lambda}_{t-1}^{-1}\left(\sum_{s=1}^{t-1}\bm{a}_s y_s\right) = \widetilde{\bm\theta}_t + \bm{\Lambda}_{t-1}^{-1}\bm{\eta}_{t-1}.
 \end{align}
 Note that $\bm{\widetilde{\theta}}_t$ is the estimator of ridge regression when no noise is present. We can easily bound the difference between $\bm{\widetilde{\theta}}_t$ and $\bm{\theta}^{\star}$ as follows:
\begin{align}
    \bm{\widetilde{\theta}}_{t} - \bm{\theta}^{\star} = \bm{\Lambda}_{t-1}^{-1}(\bm{\Lambda}_{t-1}-\bm{I}_d)\bm{\theta}^{\star}  - \bm{\theta}^{\star}= -\bm{\Lambda}_{t-1}^{-1}\bm{\theta}^{\star},
\end{align}
then one has
\begin{align}
     \left\Vert{\bm \Lambda}_{t}^{1/2}(\bm{\widetilde{\theta}}_{t} - \bm{\theta}^{\star})\right\Vert_2\leq \frac{1}{\sqrt{\lambda_{t,d}}}.
\end{align}
It remains for one to establish 
\begin{align*}
    \max_{1\leq t\leq T}\left\Vert\bm{\Lambda}_t^{-1/2}\bm{\eta}_t\right\Vert_2\lesssim\sigma\sqrt{d+\log\log T}.
\end{align*}


Towards this, we begin by observing the following identity:
\begin{align}\label{equ:noise-equivalent}
\frac{1}{2\sigma^2}\left\Vert\bm{\Lambda}_t^{-1/2}\bm{\eta}_t\right\Vert_2^2 = \frac{1}{2\sigma^2}\bm{\eta}_t^{T}\bm{\Lambda}_t^{-1}\bm{\eta}_t = \max_{\bm{\lambda}\in\mathbb{R}^d} \left\{ \bm{\lambda}^{T}\bm{\eta}_t - \frac{\sigma^2}{2}\bm{\lambda}^{T}\bm{\Lambda}_t\bm{\lambda} \right\}.
\end{align}
Consequently, in order to control the left-hand side of (\ref{equ:noise-equivalent})—which plays a key role in our proof—it suffices to control the right-hand side, that is, the supremum over $\bm{\lambda} \in \mathbb{R}^d$ of the random process $\bm{\lambda}^{T}\bm{\eta}_t - \frac{\sigma^2}{2}\bm{\lambda}^{T}\bm{\Lambda}_t\bm{\lambda}$.

To facilitate this, we introduce the following exponential process indexed by $\bm{\lambda}$ and adapted to the filtration up to time $t$:
\begin{align}\label{equ:super-martingale}
M_t(\bm{\lambda}) = \exp\left(\bm{\lambda}^{T}\bm{\eta}_t - \frac{\sigma^2}{2}\bm{\lambda}^{T}\bm{\Lambda}_t\bm{\lambda}\right).
\end{align}

Here, we claim that
\begin{align}
\label{claim:supermartingale}
    M_t(\bm{\lambda})\text{ is a supermartingale for any }\bm{\lambda}\in\mathbb{R}^d,
\end{align}
which is proved at the end of this section. 
To extend this to a uniform bound over $\bm{\lambda}$, we follow a strategy similar to that in~\cite{khamaru2024inference} and consider a weighted aggregation of the processes $M_t(\bm{\lambda})$. Specifically, we define:
\begin{align}
    Z_t = \int \gamma(\bm{\lambda}) M_t(\bm{\lambda}) \, d\bm{\lambda},
\end{align}
where $\gamma(\bm{\lambda})$ is a prior density (or mass function) over $\mathbb{R}^d$ satisfying $\int \gamma(\bm{\lambda}) \, d\bm{\lambda} = 1$. Since $M_t(\bm{\lambda})$ is a non-negative supermartingale and $\gamma(\bm{\lambda})$ integrates to one, standard results ensure that $Z_t$ is also a supermartingale.

For analytical convenience, we may take $\gamma(\bm{\lambda})$ to be a discrete prior supported on a countable set $\{ \bm{\lambda}_i \}_{i=1}^\infty$ with weights $\{ \gamma_i \}_{i=1}^\infty$, where $\sum_{i=1}^\infty \gamma_i = 1$. In this case, $Z_t$ admits the form:
\begin{align}\label{equ:Z-t}
    Z_t = \sum_{i=1}^{\infty} \gamma_i M_t(\bm{\lambda}_i),
\end{align}
which is again a supermartingale as a convex combination of supermartingales. This construction allows us to control the supremum over $\bm{\lambda} \in \mathbb{R}^d$ via a union bound or concentration argument over the discrete support, thereby paving the way for a high-probability bound on~(\ref{equ:noise-equivalent}).
 As a matter of fact, by uniform concentration of supermartingale, one may show that
\begin{align}\label{equ:martingale-concentration}
    P\left(\exists t: Z_t\geq\frac{1}{\delta}\right)\leq\delta.
\end{align}
Consider the event $E_{i,t} = \{\gamma_i M_t(\bm{\lambda}_i) \geq \delta^{-1}\}$. Observe that $E_{i,t} \subset \{Z_t \geq \delta^{-1}\}$, since $Z_t$ is a weighted sum over the $M_t(\bm{\lambda}_i)$. Taking a union over all indices $i$ and all time steps $t$, we define the event
\begin{align}
    E := \bigcup_{i=1}^{\infty} \bigcup_{t=1}^{T} E_{i,t} \subset \left\{ \exists\, t \in [T]: Z_t \geq \frac{1}{\delta} \right\},
\end{align}
which, by the concentration inequality in~(\ref{equ:martingale-concentration}), implies that $\mathbb{P}(E) \leq \delta$.
Moreover, each event $E_{i,t}$ can be equivalently rewritten in terms of an inequality involving $\bm{\eta}_t$ and $\bm{\Lambda}_t$:
\begin{align}\label{equ:E-equivalent}
    E_{i,t} = \left\{ \gamma_i \exp\left( \bm{\lambda}_i^{\top} \bm{\eta}_t - \frac{\sigma^2}{2} \bm{\lambda}_i^{\top} \bm{\Lambda}_t \bm{\lambda}_i \right) \geq \frac{1}{\delta} \right\} 
    = \left\{ \bm{\lambda}_i^{\top} \bm{\eta}_t - \frac{\sigma^2}{2} \bm{\lambda}_i^{\top} \bm{\Lambda}_t \bm{\lambda}_i \geq \log\left( \frac{1}{\gamma_i \delta} \right) \right\}. 
\end{align}

Conditioned on the complement event $E^c$, which occurs with probability at least $1 - \delta$, none of the events $E_{i,t}$ hold for any $i$ or $t$. This allows us to uniformly control the values of the process in~(\ref{equ:super-martingale}) across the net points $\{\bm{\lambda}_i\}$.

To translate this control into a bound on~(\ref{equ:noise-equivalent}), we construct a weighted net, denoted by $\mathcal{N} = \{\bm{\lambda}_i\}$ with associated weights $\{\gamma_i\}$, satisfying the following approximation guarantee: with probability $1-\delta$, for any $(\bm{\eta}_t, \bm{\Lambda}_t)$ such that $\|\bm{\Lambda}_t^{-1/2} \bm{\eta}_t\|_2 \geq \sigma^2$, there exists an index $i_t$ such that
\begin{align}\label{equ:net-property}
    \max_{\bm{\lambda} \in \mathbb{R}^d} \left\{ \bm{\lambda}^{T} \bm{\eta}_t - \frac{B^2}{2} \bm{\lambda}^{T} \bm{\Lambda}_t \bm{\lambda} \right\} 
    \leq \frac{16}{7} \left( \bm{\lambda}_{i_t}^{T} \bm{\eta}_t - \frac{\sigma^2}{2} \bm{\lambda}_{i_t}^{T} \bm{\Lambda}_t \bm{\lambda}_{i_t} \right), 
    \quad \text{with} \quad \gamma_{i_t} \gtrsim (13^d \cdot\mathrm{poly}(d)\cdot\mathrm{poly}(\log T))^{-1}.
\end{align}

As a result, we conclude that with probability at least $1 - \delta$, for every $t \in [T]$, the following bound holds:
\begin{align}
    \|\bm{\Lambda}_t^{-1/2} \bm{\eta}_t\|_2 
    \lesssim \sqrt{2}\sigma \cdot \max\left(\sigma, \sqrt{ \log\left( \frac{1}{\gamma_{i_t} \delta} \right) } \right) 
    \lesssim \sigma \sqrt{d + \log\left( \frac{\log T}{\delta} \right)}.
\end{align}
Setting $\delta = (\log T)^{-1}$ completes the proof of the desired result.



 To complete the proof of Theorem \ref{thm:noise-bound}, it remains to construct a net that satisfies the covering condition in (\ref{equ:net-property}). As previously discussed, the primary difficulty lies in the fact that the pair $(\bm{\Lambda}_t, \bm{\eta}_t)$ evolves with time and may exhibit instability, making it challenging to construct a uniform net over all $t$. To make progress, we first consider a simplified setting where $t$ is fixed. 
\paragraph{Building a net for a single pair $(\bm{\Lambda}_t, \bm{\eta}_t)$.} We begin by constructing a covering net for a single instance of the pair $(\bm{\Lambda}_t, \bm{\eta}_t)$, which serves as a foundational step toward addressing the more general case where these quantities vary with time. In this simplified setting, we do not assign weights to the elements of the net $\mathcal{N}_t$; instead, our goal is to control the size of $\mathcal{N}_t$ and ensure that it remains as small as possible while still providing sufficient coverage, i.e. with high probability there exists index $i$ such that 
\begin{align}
\label{equ:coverage-condition}
    \max_{\bm{\lambda} \in \mathbb{R}^d} \left\{ \bm{\lambda}^{T} \bm{\eta}_t - \frac{\sigma^2}{2} \bm{\lambda}^{T} \bm{\Lambda}_t \bm{\lambda} \right\} 
    \leq \frac{9}{5} \left( \bm{\lambda}_{i}^{T} \bm{\eta}_t - \frac{\sigma^2}{2} \bm{\lambda}_{i}^{T} \bm{\Lambda}_t \bm{\lambda}_i \right),\quad\text{if }\Vert\bm{\Lambda}_t^{-1/2}\bm{\eta}_t\Vert_2\geq \sigma^2.
\end{align}

To help us construct the grid, we claim that for each $(\bm{\Lambda}_t,\bm{\eta}_t)$, the points that satisfy the above condition form a ball in $\mathbb{R}^d$. We characterize that in the following lemma, which is characterized in the following lemma.

\begin{lem}\label{lem:approximate-condition}
    Let $f(\bm{\lambda}) = \bm{\lambda}^{T}\bm{\eta} - (\sigma^2/2)\bm{\lambda}^{T}\bm{\Lambda}\bm{\lambda}$. For any $\kappa\in(0,1),$ if we define the set $\mathcal{C}$ to be $\mathcal{C} = \{\bm{\lambda}_0: f(\bm{\lambda}_0)\geq (1-\kappa^2)\max_{\bm{\lambda}\in\mathbb{R}^d} f(\bm{\lambda})\}$. Then $\mathcal{C}$ can also be characterized as
    $$\mathcal{C} = \left\{\bm{\lambda}_0: \left\Vert\bm{\Lambda}^{1/2}\bm{\lambda}_0 - \frac{1}{\sigma^2}\bm{\Lambda}^{-1/2}\bm{\eta}\right\Vert_2\leq \frac{\kappa}{\sigma^2}\left\Vert\bm{\Lambda}^{-1/2}\bm{\eta}\right\Vert_2\right\}.$$
\end{lem}

By Lemma~\ref{lem:approximate-condition}, an equivalent condition for approximating the maximizer in our variational bound is the existence of a point $\bm{\lambda}_i \in \mathcal{N}_t$ such that
\begin{align}
\left\Vert \bm{\Lambda}_t^{1/2} \bm{\lambda}_i - \frac{1}{B^2} \bm{\Lambda}_t^{-1/2} \bm{\eta}_t \right\Vert_2
\leq \frac{2}{3\sigma^2} \left\Vert \bm{\Lambda}_t^{-1/2} \bm{\eta}_t \right\Vert_2,
\quad \text{whenever } \Vert \bm{\Lambda}_t^{-1/2} \bm{\eta}_t \Vert_2 \geq \sigma^2.
\end{align}
This means that at least one point $\bm{\lambda}_i \in \mathcal{N}_t$ lies within a Euclidean ball centered at $\frac{1}{B^2} \bm{\Lambda}_t^{-1/2} \bm{\eta}_t$ with radius $\frac{1}{2B^2} \Vert \bm{\Lambda}_t^{-1/2} \bm{\eta}_t \Vert_2$. To analyze this further, we state the following lemma that gives us an upper bound of $\Vert \bm{\Lambda}_t^{-1/2} \bm{\eta}_t \Vert_2$.

\begin{lem}[Adapted Version of Theorem 6.3.2, \citet{vershynin2018high}]\label{lem:linear-subgaussian}
Let $\bA$ be an $n\times d$ matrix and let $\bm{X}=(X_1,\ldots,X_n)$ be a random vector with independent mean 0, $\sigma^2$-variance and $\sigma$ sub-gaussian coordinates. Then with probability $1-\delta$, 
$$\Vert\bm{AX}\Vert_2 \lesssim   \sigma^2\left(\Vert\bA\Vert_F + \sqrt{\log\left(\frac{1}{\delta}\right)} \cdot \Vert\bA\Vert_{\mathrm{op}}\right).$$
\end{lem}

To derive the upper bound of $\Vert \bm{\Lambda}_t^{-1/2} \bm{\eta}_t \Vert_2$ from Lemma \ref{lem:linear-subgaussian}, define the matrix $\bm{A}_t = \bm{\Lambda}_t^{-1/2} (\ba_1, \ldots, \ba_t) \in \mathbb{R}^{d \times t}$ and let $\bm{X}_t = (\epsilon_1, \ldots, \epsilon_t)^{\mathsf{T}} \in \mathbb{R}^{t}$ denote the vector of noise terms. 
Noting that $\Vert\bA_t\Vert_{\mathrm{F}} = \sqrt{d}$ and $\Vert\bA_t\Vert_{\mathrm{op}} = 1$, we apply Lemma~\ref{lem:linear-subgaussian} to obtain  
\begin{align}\label{equ:upper-noise}
\left\Vert \bm{\Lambda}_t^{-1/2} \bm{\eta}_t \right\Vert_2
= \left\Vert \bm{A}_t \bm{X}_t \right\Vert_2
\lesssim \sigma^2  \sqrt{d+\log\left(\frac{1}{\delta}\right)},
\end{align}
with probability at least $1 - \delta$.

Since this is the only information we have about $\bm{\Lambda}_t^{-1/2} \bm{\eta}_t$, we shall  use this upper bound to guide the construction of the covering net. In particular, the net $\mathcal{N}_t$ must satisfy the following condition: for any vector $\bm{x} \in \mathbb{R}^d$ such that $\sigma^2\leq \Vert\bm{x}\Vert_2 \lesssim \sigma^2 \sqrt{d + \log(\delta^{-1})}$, there exists a point $\bm{\lambda} \in \mathcal{N}_t$ such that
\begin{align}\label{equ:approximate-property}
\left\Vert \bm{\Lambda}_t^{1/2} \bm{\lambda} - \frac{1}{B^2}\bm{x} \right\Vert_2
\leq \frac{2}{3\sigma^2} \left\Vert \bm{x} \right\Vert_2.
\end{align}
To construct the net $\mathcal{N}_t$ for a fixed time step $t$, we proceed in three stages:

\begin{itemize}
    \item \textbf{Step 1: Construct a base $\epsilon$-net on the unit sphere.}  
    We begin by constructing an $\epsilon$-net on the unit Euclidean sphere $\mathcal{S}^{d-1}$ in $\mathbb{R}^d$, where we set $\epsilon = 1/6$. Let $\mathcal{N}_0$ denote this covering. By a standard volume argument (see, e.g., \citet{vershynin2018high}), such a net exists with cardinality bounded by
    $$|\mathcal{N}_0| \leq \left(1 + \frac{2}{\epsilon}\right)^d \leq 13^d.$$
    This ensures that every point on the unit sphere lies within Euclidean distance $\epsilon$ of some point in $\mathcal{N}_0$.

    \item \textbf{Step 2: Construct nets on expanding spheres.}  
    To cover the full range of norms that may arise in the transformed space (i.e., after rescaling by $\bm{\Lambda}_t^{-1/2}$ and $B^{-2}$), we scale $\mathcal{N}_0$ to form nets on concentric spheres of increasing radii. Specifically, for $j=0,1,\ldots,k$ we construct $(3/2)^j\epsilon$ covering on sphere with radii $(3/2)^j$, where $k$ is chosen to ensure that the largest radius exceeds the typical scale of the transformed vector. In particular, we set 
    $$k = O\left( \log(d + \log(\delta^{-1})) \right),$$ 
    so that the maximum radius $(3/2)^k$ covers the high-probability upper bound of $\sigma^{-2}\| \bm{\Lambda}_t^{-1/2} \bm{\eta}_t\|_2$.

    \item \textbf{Step 3: Construct the final net in the transformed space.}  
    For each $j$, we map the scaled points from the unit sphere net through the transformation $\bm{\lambda} = (3/2)^j \cdot \bm{\Lambda}_t^{-1/2} \bm{v}$, where $\bm{v} \in \mathcal{N}_0$. This results in a net that discretizes the ellipsoidal region defined by the inverse covariance geometry of $\bm{\Lambda}_t$. The final net for time $t$ is thus given by
    \begin{align} \label{equ:net-construction-single}
        \mathcal{N}_t = \left\{ \bm{\lambda} = 2^j \cdot \bm{\Lambda}_t^{-1/2} \bm{v} : \bm{v} \in \mathcal{N}_0,\; j=0,1,\ldots,k \right\},
    \end{align}
    to account for the full effective range of relevant $\bm{x}$ vectors satisfying $\sigma^2 \leq \|\bm{x}\|_2 \lesssim \sigma^2 \sqrt{d + \log(\delta^{-1})}$.

\end{itemize}

With this construction, the total size of each individual $\mathcal{N}_t$ satisfies
\begin{equation}\label{equ:single-net-size}
N_{\delta} = k|\mathcal{N}_0| \lesssim 13^d \log(d + \log(\delta^{-1})).
\end{equation}
We now show that $\mathcal{N}$ satisfies the approximate property stated in~(\ref{equ:approximate-property}), which concludes our construction of $\mathcal{N}$ in this case.


Let $j = \left\lfloor \log_2\left(\sigma^{-2} \|\bm{x}\|_2\right) \right\rfloor$ be the largest index such that $2^j \leq \sigma^{-2} \|\bm{x}\|_2$. By construction, we include a $2^j \epsilon$-net on the sphere $2^j \mathcal{S}^{d-1}$ in $\mathcal{N}$, where $\epsilon = 1/6$. We first observe that
\begin{equation} \label{equ:approximate-1-single}
d\left(\sigma^{-2} \bm{x}, 2^j \mathcal{S}^{d-1}\right) \leq \frac{1}{3\sigma^2} \|\bm{x}\|_2,
\end{equation}
where $d(\cdot, \cdot)$ denotes the Euclidean distance in $\mathbb{R}^d$. Let
\[
\bm{z} = \arg\min_{\bm{y} \in 2^j \mathcal{S}^{d-1}} \left\| \sigma^{-2} \bm{x} - \bm{y} \right\|_2.
\]
Then $d(\sigma^{-2} \bm{x}, 2^j \mathcal{S}^{d-1}) = \|\sigma^{-2} \bm{x} - \bm{z}\|_2$. Since $\mathcal{N}$ contains a $2^j \epsilon$-net on $2^j \mathcal{S}^{d-1}$, there exists $\bm{\lambda} \in \mathcal{N}$ such that
\begin{equation} \label{equ:approximate-2-single}
\|\bm{z} - \bm{\lambda}\|_2 \leq \frac{2^j}{6} \leq \frac{1}{6\sigma^2} \|\bm{x}\|_2.
\end{equation}

Combining~(\ref{equ:approximate-1-single}) and~(\ref{equ:approximate-2-single}) using the triangle inequality, we obtain
\begin{align}
\left\| \sigma^{-2} \bm{x} - \bm{\lambda} \right\|_2 
\leq \left\| \sigma^{-2} \bm{x} - \bm{z} \right\|_2 + \left\| \bm{z} - \bm{\lambda} \right\|_2 \leq \frac{1}{3\sigma^2} \|\bm{x}\|_2 + \frac{1}{6\sigma^2} \|\bm{x}\|_2 = \frac{1}{2\sigma^2} \|\bm{x}\|_2.
\end{align}
This proves that $\mathcal{N}$ satisfies the desired approximation.

\paragraph{Building a net with uniform approximation on all $(\bm{\Lambda}_t, \bm{\eta}_t)$.} We now turn to the more general case, where the goal is to construct a net that uniformly approximates all pairs $(\bm{\Lambda}_t, \bm{\eta}_t)$ for every $t \in [T]$. A straightforward extension of the previous construction would be to build a separate net $\mathcal{N}_t$ for each $t$. However, this naive approach yields a total of at least $O(T)$ points in the combined net, implying that some individual points would be assigned weights of order $O(T^{-1})$. This clearly contradicts the desired guarantee in~(\ref{equ:net-property}). Consequently, we require an alternative construction of a global net $\mathcal{N}$ whose size scales only as $O(\log T)$, while still ensuring uniform approximation across all time indices $t \in [T]$.


Our construction is inspired by rare-switching techniques commonly used in online reinforcement learning. However, in contrast to \citep{he2023nearly,tan2025actor}, where rare-switching is applied to policy updates to control the growth of the function class, we adopt a different perspective. Specifically, we apply the rare-switching principle to select a small number—only $O(\log T)$—of representative time indices $\{t_i\}$, and construct nets $\mathcal{N}_{t_i}$ based solely on the pairs $(\bm{\Lambda}_{t_i}, \bm{\eta}_{t_i})$.

The key insight is that even though we are not covering every $(\bm{\Lambda}_t, \bm{\eta}_t)$ individually, these $O(\log T)$ representative nets are sufficient to uniformly approximate all $(\bm{\Lambda}_t, \bm{\eta}_t)$ across $t \in [T]$. In other words, for any $t$, there exists some $t_i$ such that the corresponding net $\mathcal{N}_{t_i}$ provides a good approximation for $(\bm{\Lambda}_t, \bm{\eta}_t)$. This significantly reduces the size of the overall net $\mathcal{N}$ while preserving the desired approximation guarantees.

We now describe the construction procedure for the weighted net $\mathcal{N}$, which leverages this rare-switching idea to achieve efficient coverage over all $(\bm{\Lambda}_t, \bm{\eta}_t)$ with only logarithmically many representative components. The construction proceeds as follows:

\begin{itemize}
    \item \textbf{Step 1: Selection of representative time indices.}  
    Initialize $s_1 = 1$, and recursively define the sequence $\{s_2, s_3, \ldots, s_j\}$ according to the rule:
    \begin{align}\label{equ:rare-switching}
        s_{i+1} = \min\left\{t \leq T : \det(\bm{\Lambda}_t) > 2 \cdot \det(\bm{\Lambda}_{s_i})\right\}.
    \end{align}
    This procedure continues until no further such $t$ exists. We claim that the following properties hold true for the switching times defined in (\ref{equ:rare-switching}). 
    \begin{align}
    \label{equ:property-switching}
        j = O(d\log T);\quad s_{i+1}-s_i\leq 2^{i-1}d,\;\forall\; 1\leq i\leq j.
    \end{align}
   The intuition is that each $\bm{\Lambda}_{s_i}$ represents a “scale” of covariance growth, and doubling the determinant indicates a substantial geometric change in the feature space. With this setup, we ensure that $\bm{\Lambda}_t$ ($s_i\leq t<s_{i+1})$ does not change much from $\bm{\Lambda}_{s_i}$.

    \item \textbf{Step 2: Construction of local nets.}  
    For each selected time index $s_i$, construct a net $\mathcal{N}_{s_i}$ following the procedure described in~(\ref{equ:net-construction-single}), with $\delta_i = \frac{1}{2^{i-1}d^2\log^2 T} $. Set $N_i = |\mathcal{N}_{s_i}|$.

    \item \textbf{Step 3: Aggregation into a weighted net.}  
    Let $\mathcal{N}$ be the weighted union of all constructed nets $\{\mathcal{N}_{s_i}\}_{i=1}^j$. Assign each net $\mathcal{N}_{s_i}$ a total weight of $\gamma_i = \frac{1}{i(i+1)}$, distributed uniformly across its $N$ points. The remaining probability mass is assigned to the zero vector $\bm{0}$ with weight $\gamma_0 =\frac{1}{j+1}$. This ensures the total weight sums to one. The resulting weighted net can be written explicitly as:
    \begin{align}\label{equ:multiple-net}
        \mathcal{N} = \left\{
        \left(\overline{\gamma}_i, \bm{\lambda}_{i,m}\right) :
        1 \leq i \leq j,\;
        1 \leq m \leq N_i,\;
        \overline{\gamma}_i = \frac{1}{i(i+1)N_i},\;
        \bm{\lambda}_{i,m} \in \mathcal{N}_{s_i}
        \right\} \cup \left\{(\gamma_0, \bm{0})\right\}.
    \end{align}
\end{itemize}



We will then show below that the net $\mathcal{N}$ satisfies property  (\ref{equ:net-property}) for any $(\eta_t,\bm{\Lambda}_t)$.
     By Lemma \ref{lem:approximate-condition}, we only need to show that for any  $t$ that satisfies $s_i\leq t<s_{i+1}$, there exists $\bm{\lambda}\in\mathcal{N}_{s_i}$ such that
\begin{align}\label{equ:net-coverage-2}
\left\Vert \bm{\Lambda}_t^{1/2} \bm{\lambda} - \frac{1}{\sigma^2} \bm{\Lambda}_t^{-1/2} \bm{\eta}_t \right\Vert_2
\leq \frac{3}{4\sigma^2} \left\Vert \bm{\Lambda}_t^{-1/2} \bm{\eta}_t \right\Vert_2,
\quad \text{whenever } \Vert \bm{\Lambda}_t^{-1/2} \bm{\eta}_t \Vert_2 \geq \sigma^2.
\end{align}
As we note that 
\begin{align*}
    1\leq \lambda_{\min}( \bm{\Lambda}_t^{1/2}\bm{\Lambda}_{s_i}^{-1/2})\leq \lambda_{\max}( \bm{\Lambda}_t^{1/2}\bm{\Lambda}_{s_i}^{-1/2}) \leq \sqrt{2},\quad \frac{1}{\sqrt{2}}\leq \lambda_{\min}(\bm{\Lambda}_{s_i}^{1/2}\bm{\Lambda}_{t}^{-1/2})\leq \lambda_{\max}(\bm{\Lambda}_{s_i}^{1/2}\bm{\Lambda}_{t}^{-1/2})\leq 1,
\end{align*}
   conditioned on $\sqrt{2}\leq \sigma^{-2}\Vert\bm{\Lambda}_{t}^{-1/2}\bm{\eta}_t\Vert_2\leq R$, it holds that
\begin{align}
    \left\Vert \bm{\Lambda}_t^{1/2} \bm{\lambda} - \frac{1}{B^2} \bm{\Lambda}_t^{-1/2} \bm{\eta}_t \right\Vert_2 
    &= \left\Vert (\bm{\Lambda}_{t}^{1/2}\bm{\Lambda}_{s_i}^{-1/2})\bm{\Lambda}_{s_i}^{1/2} \bm{\lambda} - \frac{1}{\sigma^2} (\bm{\Lambda}_t^{1/2}\bm{\Lambda}_{s_i}^{-1/2})(\bm{\Lambda}_{s_i}^{1/2}\bm{\Lambda}_{t}^{-1/2}) \bm{\Lambda}_{t}^{1/2}\bm{\eta}_t \right\Vert_2\nonumber\\
    & \leq \sqrt{2} \left\Vert \bm{\Lambda}_{s_i}^{1/2} \bm{\lambda} - \frac{1}{\sigma^2} (\bm{\Lambda}_{s_i}^{1/2}\bm{\Lambda}_{t}^{-1/2}) \bm{\Lambda}_{t}^{-1/2}\bm{\eta}_t \right\Vert_2\nonumber\\
    & \leq \sqrt{2} \max_{1\leq \Vert\bm{x}\Vert_2\leq R} \left\Vert\bm{\Lambda}_{s_i}^{1/2} \bm{\lambda} - \bm{x}\right\Vert_2.
\end{align}
Next, we are going to give a uniform upper bound for $\max_{s_i\leq t<s_{i+1}}\sigma^{-2}\Vert\bm{\Lambda}_{t}^{-1/2}\bm{\eta}_t \Vert_2$. We utilize the following lemma to show this result

From Lemma \ref{lem:linear-subgaussian}, we know that with probability $1-\frac{1}{2^{i-1}d^2\log^2T}$, we have 
\begin{align}
    \Vert\bm{\Lambda}_{t}^{-1/2}\bm{\eta}_t \Vert_2\lesssim \sigma^2\sqrt{d+i+\log\log T}\lesssim \sigma^2\sqrt{d\log T}.
\end{align}
From Claim (\ref{equ:property-switching}) that $s_{i+1}-s_i\leq 2^{i-1}d$, we obtain that with probability $1-\frac{1}{d^2\log^2 T}$, one has 
\begin{align}
    \max_{s_i\leq t<s_{i+1}} \Vert\bm{\Lambda}_{t}^{-1/2}\bm{\eta}_t \Vert_2\lesssim \sigma^2\sqrt{d\log T}.
\end{align}
Therefore, from (\ref{equ:single-net-size}), we know that the size of $\mathcal{N}_{s_i}$ can be bounded as 
\begin{align}
    N_i\lesssim 13^d\log \left(\frac{1}{\sigma^2}\max_{s_i\leq t<s_{i+1}} \Vert\bm{\Lambda}_{t}^{-1/2}\bm{\eta}_t \Vert_2\right)\lesssim 13^d(\log d+\log\log T).
\end{align}
With this setup, we ensure that with probability $1-\frac{1}{d\log^2T}$, for any $s_i\leq t<s_{i+1}$, there exists $\bm{\lambda}\in\mathcal{N}_{s_i}$ such that 
\begin{align}
    \left\Vert \bm{\Lambda}_t^{1/2} \bm{\lambda} - \frac{1}{\sigma^2} \bm{\Lambda}_t^{-1/2} \bm{\eta}_t \right\Vert_2
&\leq \sqrt{2} \max_{1\leq \Vert\bm{x}\Vert_2\leq R} \left\Vert\bm{\Lambda}_{s_i}^{1/2} \bm{\lambda} - \bm{x}\right\Vert_2\nonumber\\
&\leq\frac{3}{4\sigma^2}\left\Vert \bm{\Lambda}_t^{-1/2} \bm{\eta}_t \right\Vert_2,
\quad \text{whenever } \Vert \bm{\Lambda}_t^{-1/2} \bm{\eta}_t \Vert_2 \geq \sigma^2.
\end{align}
In this way, we ensure that with probability $1-\frac{1}{\log T}$, for any $1\leq t\leq T$, there exists some $\mathcal{N}_{s_i}$ and  $\bm{\lambda}\in\mathcal{N}_{s_i}$ such that 
    \begin{align}
\left\Vert \bm{\Lambda}_t^{1/2} \bm{\lambda} - \frac{1}{\sigma^2} \bm{\Lambda}_t^{-1/2} \bm{\eta}_t \right\Vert_2
\leq \frac{3}{4\sigma^2} \left\Vert \bm{\Lambda}_t^{-1/2} \bm{\eta}_t \right\Vert_2,
\quad \text{whenever } \Vert \bm{\Lambda}_t^{-1/2} \bm{\eta}_t \Vert_2 \geq \sigma^2.
\end{align}

We guarantee that the weight $\overline{\gamma}_i$ that assigned on any points satisfies
\begin{align}
    \overline{\gamma}_i = \frac{1}{i(i+1)N_i}\gtrsim \frac{1}{13^d d^2\log^2T(\log d +\log\log T)} = \frac{1}{13^d\mathrm{poly}(d)\mathrm{poly}(\log T)},
\end{align}
which concludes (\ref{equ:net-property}).
\paragraph{Proof of Claim~(\ref{claim:supermartingale}).}

    We note that for any $\bm{\lambda}\in\mathbb{R}^d$, the following inequality holds.
    \begin{align}
\mathbb{E}[M_{T}(\bm{\lambda}) \mid \mathcal{F}_t] 
&= \mathbb{E}\left[
    M_t(\bm{\lambda}) \exp\left(
        \bm{\lambda}^\top \ba_{T} \epsilon_{T} 
        - \frac{\sigma^2}{2} \bm{\lambda}^\top \ba_{T} \ba_{T}^\top \bm{\lambda}
    \right)
    \,\middle|\, \mathcal{F}_t
\right] \nonumber \\
&= M_t(\bm{\lambda}) 
   \, \mathbb{E}\left[
       \exp\left(\bm{\lambda}^\top \ba_{T} \epsilon_{T}\right)
       \,\middle|\, \mathcal{F}_t
   \right] 
   \exp\left(
       -\frac{\sigma^2}{2} \bm{\lambda}^\top \ba_{T} \ba_{T}^\top \bm{\lambda}
   \right) \nonumber \\
&\leq M_t(\bm{\lambda}),
\end{align}
where the last inequality holds because $\epsilon_t$ is $\sigma$ sub-Gaussian random variable. As a result, we show that $M_t(\bm{\lambda})$ is a supermartingale.

\paragraph{Proof of Claim~(\ref{equ:property-switching}).}

    Note that $\mathrm{tr}(\bm{\Lambda}_1) = d$ and $\mathrm{tr}(\bm{\Lambda}_t) = d+t-1$. As a result, the determinant of $\mathrm{det}(\bm{\Lambda}_T)$ can be upper bounded as 
    \begin{align}
       \mathrm{det}(\bm{\Lambda}_T)\leq \left(\frac{d+T-1}{d}\right)^{d},
    \end{align}
    and as a result, the number of switch times $j$ can be bounded as 
    \begin{align}
      j \leq O\left[\log\left(\frac{\mathrm{det}(\bm{\Lambda}_T)}{\mathrm{det}(\bm{\Lambda}_1)}\right)\right] \leq O(d\log(d+T-1))\leq O(d\log T).
    \end{align}
    For the second inequality, we note that $\mathrm{tr}(\bm{\Lambda}_{s_{i+1}})\leq 2\mathrm{tr}(\bm{\Lambda}_{s_{i}})$. Therefore, $s_{i+1}+d-1\leq 2(s_i+d-1)$. By induction, we note that $s_{i}+d-1\leq 2^{i-1}(s_1+d-1) = 2^{i-1}d$. As a result, we establish the inequality that 
    \begin{align}
        s_{i+1} - s_i\leq s_i + d-1\leq 2^{i-1}d.
    \end{align}

\label{app:action-vector-decomp}

\section{Proof of Theorem \ref{thm:covariance}}
\label{sec:proof-thm-3}
We present the complete proof of Theorem~\ref{thm:covariance} in this section. We begin with a high-level overview of the argument by phases in Appendix~\ref{sec:proof-ideas-thm3}, and then provide the detailed proofs of each phase from Appendix~\ref{sec:phase-1}–\ref{sec:proof-phase-4}.

For simplicity, throughout this section, we write with slight abuse of notation $\bm{\phi}(\bm{x}_t,\bm{a}_t) = \bm{a}_t$. 
While in general, the context can shift the mean and covariance of the action set, we can always transform it back to the unit ball, and so we utilize this for ease of notation. We also assume without loss of generality that the true signal $\bm{\theta}^\star = \bm{e}_1$ for our theoretical analysis, since the final result does not change up to a rotation. We also denote $\lambda_{t,1},\ldots,\lambda_{t,d}$ as the eigenvalues of $\bm{\Lambda}_t$ ranked in decreasing order and $\bm{v}_{t,1},\ldots,\bm{v}_{t,d}$ be the corresponding eigenvectors.
\label{sec:proof-thm-2}
\subsection{Key proof ideas}
\label{sec:proof-ideas-thm3}
We summarize the high level proof ideas of each phase throughout the whole process.

\paragraph{Key analysis steps in Phase \#1.}  The analysis of Phase \#1 is divided into the following steps.


\begin{enumerate}
    \item \emph{Nontrivial mass on the bottom subspace.}
In the initial phase, the ratio satisfies $\beta/\sqrt{\lambda_{t,d}}\gtrsim 1$. Consequently, $\mathrm{UCB}_t(\bm a)$ is dominated by the exploration bonus rather than the estimated reward, favoring less–explored directions. Hence $\bm a_t$ necessarily places a constant fraction of its energy on the small–eigenvalue subspace (Lemma \ref{lem:opt-concentration}), i.e., on indices $i$ with $\lambda_{t,i}\le C\,\lambda_{t,d}$ for a fixed constant $C>1$.


    \item \emph{Rank-one updates lift the bottom.}
Each update is $\bm a_t\bm a_t^\top$, and as a result, the eigenvalue increments of $\bm{\Lambda}_t$ is approximated by the squared projections of $\bm a_t$ onto the eigenbasis (Lemma~\ref{lem:rank-one-update}). Since $\bm a_t$ places a constant mass on the small-eigenvalue subspace, the small-eigenvalue block gains a uniformly positive aggregate amount each round. 

   \item \emph{Lower bound for the minimum eigenvalue.}
Distributing the persistent aggregate gain over at most $d$ coordinates forces the minimum to rise at least linearly: the bottom block accrues a constant total increase each round, so its per–coordinate average grows by at least a constant multiple of $1/d$ per round, yielding $\lambda_{t,d}\gtrsim t/d$, for all $t$ in the first stage when $\beta/\sqrt{\lambda_{t,d}}\gtrsim 1$, and in fact $\lambda_{t,d}\asymp t/d$ since it cannot exceed the average eigenvalue.

    \item \emph{Exit and eigengap.}
Meanwhile, at some point when $\beta/\sqrt{\lambda_{t,d}}\lesssim 1$, the top eigenvalue grows strictly faster (by Theorem~\ref{thm:noise-bound} and Lemma~\ref{lem:decom-a}), while the rest are constrained by the mass-splitting above. Hence there exists $t_1$ such that $\beta/\sqrt{\lambda_{t_1,d}}\asymp 1$,  an eigengap emerges, completing the first stage.
\end{enumerate}

\paragraph{Key analysis steps in Phase \#2.}  Set $\lambda_{t,d}=c_t\,\beta\sqrt{t}$ for $t\ge t_1$. It suffices to show $c_t\asymp d^{-1/2}$, which pins down $\lambda_{t,d}\asymp \beta\sqrt{t/d}$.
    \begin{enumerate}
        \item \emph{Alignment of the top eigenvector.} A crude concentration bound gives
$\|\bm{a}_t-\bm{\theta}^\star\|_2^2\lesssim \beta/\bigl( c_t\sqrt t\bigr)$.
 Since $\operatorname{tr}(\bm{\Lambda}_t)\approx t$, a Rayleigh–Ritz/Davis–Kahan argument yields
$\|\bm{v}_{t,1}-\bm{\theta}^\star\|_2^2\lesssim \beta/\bigl(\underline c_t\sqrt t\bigr)$, with $\underline c_t=\min_{t_1\le s\le t}c_s$. Here $\underline c_t$ is the historical normalized floor—the smallest value the normalized minimum eigenvalue has attained up to time $t$. The bound depends on $1/\underline c_t$ (rather than $1/c_t$) because it aggregates past rounds: the process cannot “forget’’ earlier times when the floor was lower and exploration bonuses were larger.

        \item \emph{Non-leading spectrum grows at the $\beta/\sqrt t$ scale.}
Let $\overline\lambda_t$ be the mean of the non-leading eigenvalues. The one-step change $\overline\lambda_{t+1}-\overline\lambda_t$ can be upper bounded as $\overline\lambda_{t+1}-\overline\lambda_t
= O\bigl(\beta/(d\underline c_t\sqrt t)\bigr)$. Consequently, it is controlled by the alignment of $\bm v_{t,1}$ (better alignment means less spillover into non-leading directions), hence its dependence on the historical bottleneck $\underline c_t$. One-step changes can also be lower bounded as  $\overline\lambda_{t+1}-\overline\lambda_t= \Omega\bigl(\beta/(dc_t\sqrt t)\bigr)$, where $c_t$ comes from the exploration of the UCB objective, which allocates nontrivial weight to under-explored directions and thus scales with the current floor $c_t$. As the result, the upper/lower envelopes match up to the ratio $c_t/\underline c_t$.

        \item \emph{Relative concentration on small eigenvalues.} For times with $c_t\le \widetilde c$ (a fixed constant $\widetilde c>0$) and $c_t/\underline c_t=O(1)$, we are in a regime where the per–coordinate exploration bonus
$\beta w_{t,i}/\sqrt{\lambda_{t,i}}$
is comparable to the perturbation induced by the top direction (controlled by $\|\bm v_{t,1}-\bm\theta^\star\|$) for the non-leading coordinates.
Consequently, the optimizer of $\mathrm{UCB}_t$ places a constant fraction of its non-leading mass on indices with small $\lambda_{t,i}$ (by Lemma~\ref{lem:opt-concentration}).
Repeating the rank-one update argument from Proposition~\ref{prop:first-stage} (cf. Lemma~\ref{lem:rank-one-update}), this relative mass split transfers to growth: a fixed constant fraction of the total non-leading increment in each round is captured by the small-eigenvalue block.
         \item \emph{Forcing up the minimum and closing the loop.} Because a fixed fraction of the non-leading increment lands on the bottom block each round, and that
$\overline{\lambda}_t$ enjoys a per-step lower bound shown in Step 2, we can show that $\lambda_{t,d}$ grows at least at the rate of $\Omega(\beta/(dc_t\sqrt{t}))$ in this regime, as a result, $c_t\gtrsim 1/(dc_t)$, leading to a lower bound $\lambda_{t,d}\gtrsim \beta\sqrt t/\sqrt {d}$,  when $t\ge t_1$. Conversely, $\lambda_{t,d}$ cannot exceed the non-leading average, yielding the matching upper bound and hence
$\lambda_{t,d}\asymp \beta\sqrt{t/d}$. Plugging this back into the alignment bound gives $\|\bm{v}_{t,1}-\bm{\theta}^\star\|_2^2\lesssim \beta^2/\lambda_{t,d}$, which completes the argument.
    \end{enumerate}

\paragraph{Key analysis steps in Phase \#3.} We present the main steps of the proof of Phase \#3 as follows.
\begin{enumerate}
\item We track how the leading eigenvector $\bm{v}_{t,1}$ evolves as the sample size increases from $t$ to $t{+}1$. Using that $\bm{a}_t$ optimizes the UCB objective, the new sample induces a rank-one perturbation linking $(\bm{v}_{t,1},\widehat{\bm{\theta}}_t)$ to $\bm{v}_{t+1,1}$:
\[
\bm{v}_{t+1,1}
= \bm{v}_{t,1}
+ \frac{\widehat{\bm{\theta}}_t-\bm{v}_{t,1}}{t}
+ \bm{\zeta}_t.
\]
The fluctuation term $\bm{\zeta}_t$ has, in the worst case, only higher-order adverse effect on alignment with $\widehat{\bm{\theta}}_t$, whereas any component that improves alignment may be non-negligible. Consequently, the update makes $\bm{v}_{t+1,1}$ closer to $\widehat{\bm{\theta}}_t$ than $\bm{v}_{t,1}$. As a result, $\|\bm{v}_{t,1}-\bm{\theta}^{\star}\|_2$ decreases until it is of the same order as $\|\widehat{\bm{\theta}}_{t}-\bm{\theta}^{\star}\|_2$.

\item Leveraging upon the above update, we study how error propogates with time step $t$ and obtain the following 
\begin{align}
\bigl\|\bm{v}_{t+1,1}-\bm{\theta}^\star\bigr\|_2^2
\;\le\;
\left(\Bigl(1-\frac{1}{t}\Bigr)\bigl\|\bm{v}_{t,1}-\bm{\theta}^\star\bigr\|_2
+\widetilde{O}(t^{-5/4})\right)^2
+\widetilde{O}(t^{-2}),
\end{align}
so the error contracts by roughly $(1-1/t)$ with faster-vanishing additive terms. Combining this with the Phase~\#2 initialization, there exists $t_2 = O\!\left(\beta^{8}/(\sigma^6 d^2)\right)$ such that for all $t \ge t_2$, $\bm{v}_{t,1}$ attains the desired concentration around $\bm{\theta}^\star$, matching the concentration order of $\widehat{\bm{\theta}}_t$.

\end{enumerate}

\paragraph{Key analysis steps in Phase \#4.}     We outline the key steps in Phase \#4 as follows.
    \begin{enumerate}
        \item \emph{Precise decomposition of $\bm{a}_t$.} Using the refined concentration of the top eigenvector (the high-probability bound on $\|\bm v_{t,1}-\widehat{\bm\theta}_t\|_2$), we sharpen Lemma~\ref{lem:decom-a} and make explicit how $\bm a_t$ splits between the leading direction and its orthogonal complement:

\begin{itemize}
\item \emph{Non-leading directions.}
The non-leading coordinates of $\bm w_t$ (defined in (\ref{equ:w-t-expression})) satisfy
\[
\sum_{i=2}^{d} w_{t,i}^2\!\left(1-\frac{\lambda_{t,d}}{\lambda_{t,i}}\right)=o(1).
\]
This equation implies that $\bm w_t$, and hence $\bm a_t$, concentrates most of its non-leading mass on directions whose eigenvalues are close to $\lambda_{t,d}$.

\item \emph{Mass off the top.}
The portion of $\bm a_t$ lying outside the top direction admits the precise characterization
\[
\sum_{i=2}^{d}\kappa_{t,i}^2
=\frac{\beta^2}{\lambda_{t,d}}
\Bigl(1-\frac{\lambda_{t,d}^2}{\beta^2 t}+o(1)\Bigr),
\]
as long as $\lambda_{t,d}\leq\beta\sqrt{t}$. Therefore, we can precisely characterize the growth speed of non-leading eigenvalues under this regime.
\end{itemize}

    
  \item \emph{Preliminary growth-rate control for $\lambda_{t,d}$.} We prove that there exists $t_2' = O(t_2)$ such that for all $t \ge t_2'$, $\lambda_{t,d} \le c\,\beta \sqrt{t}$ for some constant $c<1$. This bound is essential: as shown in the previous step, the growth of the non-leading eigenvalues can only be characterized sharply when $\lambda_{t,d}$ is smaller than $\beta\sqrt{t}$. We need this requirement to establish further fine-grained arguments.
    
    \item \emph{Limiting the projection of $\bm{a}_t$ onto eigenspaces with large eigenvalues.} Using the precise decomposition from Step~1, define the set of large eigenvalues as those exceeding $(1+\delta)\lambda_{t,d}$ with $\delta=o(1)$. Then $\bm{a}_t$ allocates at most an $O(1/d)$ fraction of its non-leading mass to these large eigenvalues.
    
    \item \emph{Controlling the growth of large non-leading eigenvalues.} We show that eigenvalues above the $(1+d\delta)\lambda_{t,d}$ threshold do not grow faster than $\overline{\lambda}_t$. Consequently, the non-leading eigenvalues concentrate, differing in magnitude by at most a $(1+o(1))$ factor.
    \item \emph{Precise characterization of non-leading eigenvalues.}  By previous results, the non-leading eigenvalues coalesce: for \(t\ge t_4\),
\(\lambda_{t,2}\approx\cdots\approx\lambda_{t,d}\). It is therefore natural to track their common level via the average
\(\overline{\lambda}_t := \frac{1}{d-1}\sum_{i=2}^d \lambda_{t,i}\).
For \(t\ge t_4\), \(\overline{\lambda}_t\) evolves according to
\[
\overline{\lambda}_{T}
= \overline{\lambda}_t
+ \frac{1}{d-1}\,\frac{\beta^2}{\overline{\lambda}_t}
\left(1-\frac{\overline{\lambda}_t^2}{\beta^2 t}+o(1)\right).
\]
This recursion makes the growth of \(\overline{\lambda}_t\) explicit: to leading order it increases at rate \(\beta^2/((d-1)\overline{\lambda}_t)\), with a vanishing correction of order \(\overline{\lambda}_t/(\beta^2 t)\).
In turn, it delivers a precise large-\(t\) asymptotic for \(\overline{\lambda}_t\) together with its first-order correction. This constructs our desired result.
    \end{enumerate}

\subsection{Analysis of Phase \#1 (proof of Proposition \ref{prop:first-stage})}
\label{sec:phase-1}

    


We aim to show that whenever $\beta/\sqrt{\lambda_{t,d}}\geq C'$, there exists constant $C = C(C')$ such that 
$$\lambda_{t,d}\geq C\cdot\frac{t}{d}.$$

\paragraph{Step 1: lower bound the projection of $\bm{a}_t$ on eigenspaces with ``low eigenvalues."}  
At the beginning, the following conditions are satisfied:  
\[
\frac{\beta}{\sqrt{\lambda_{t,d}}} \;\geq\; C', 
\qquad 
\sum_{i=1}^{d}\nu_{t,i}^2 = 1,
\]  
for some constant $C'>0$. The first inequality guarantees that the effective signal-to-noise ratio remains bounded away from zero, while the second condition normalizes the direction vector $\nu_t$.  

Under these assumptions, the optimization problem (\ref{equ:w-t-expression}) falls within the scope of Lemma~\ref{lem:opt-concentration}. By direct application of this lemma, there exist absolute constants $C_1=C_1(C')$ and $C_2\in(0,1]$ such that, if we define  
\[
k_t \;=\; \max\bigl\{i : \lambda_{t,i} > C_1 \lambda_{t,d} \,\bigr\},
\]  
then the coefficients $\kappa_{t,i}$ in the expansion of $\bm{a}_t$ must satisfy  
\begin{align}
\sum_{i = k_t + 1}^d \kappa_{t,i}^2 \;\;\geq\;\; C_2.
\end{align}  
This inequality formalizes the intuition that a nontrivial fraction of the action vector necessarily lies in the less dominant eigenspaces—those associated with eigenvalues not substantially larger than $\lambda_{t,d}$. In other words, $\bm{a}_t$ cannot concentrate exclusively on the top eigen-directions; a uniformly positive share of its energy always projects onto the ``low-eigenvalue'' space.  

\paragraph{Step 2: control the growth of ``large eigenvalues".}  
Next, define the set of ``large eigenvalues'' to be those exceeding $2C_1\lambda_{t,d}$. Let  
\[
k_t = \max\left\{j: \lambda_{t,j}> 2C_1\lambda_{t,d}\right\}.
\]  

From Lemma~\ref{lem:rank-one-update}, the updated eigenvalues $\lambda_{t+1,1},\ldots,\lambda_{t+1,d}$ are precisely the roots of  
\[
f(\lambda) = 1+\sum_{i=1}^{d}\frac{\kappa_{t,i}^2}{\lambda_{t,i}-\lambda}.
\]  

To isolate the contribution of the ``large'' coordinates, we define the auxiliary function  
\[
f_1(\lambda) = 1+\sum_{i=1}^{k_t}\frac{\kappa_{t,i}^2}{\lambda_{t,i}-\lambda}.
\]  
For any $i\leq k_t$, one has  
\begin{align}
f(\lambda_{t+1,i}) 
= f_1(\lambda_{t+1,i}) + \sum_{i=k_t+1}^{d}\frac{\kappa_{t,i}^2}{\lambda_{t,i}-\lambda_{t+1,i}} \geq f_1(\lambda_{t+1,i}) - \frac{1}{C_1\lambda_{t,d}}\sum_{i=k_t+1}^{d}\kappa_{t,i}^2,
\end{align}  
which implies  
\begin{align} \label{equ:upper-bound-large}
f_1(\lambda_{t+1,i}) \;\leq\; \frac{1}{C_1}\sum_{i=k_t+1}^{d}\kappa_{t,i}^2.
\end{align}  

Let $\widetilde{\lambda}_{t,1},\ldots,\widetilde{\lambda}_{t,k_t}$ be the solutions of  
\[
f_1(\lambda) = \frac{1}{C_1}\sum_{i=k_t+1}^{d}\kappa_{t,i}^2.
\]  
Equivalently, these roots satisfy  
\begin{align}
f_1(\lambda) - \frac{1}{C_1}\sum_{i=k_t+1}^{d}\kappa_{t,i}^2 
&= 1-\frac{1}{C_1}\sum_{i=k_t+1}^{d}\kappa_{t,i}^2 + \sum_{i=1}^{k_t}\frac{\kappa_{t,i}^2}{\lambda_{t,i}-\lambda}\nonumber\\
&= \frac{\left(1-\frac{1}{C_1}\sum_{i=k_t+1}^{d}\kappa_{t,i}^2\right)\cdot \prod_{i=1}^{k_t}(\lambda_{t,i}-\lambda) +\sum_{i=1}^{k_t} \kappa_{t,i}^2 \prod_{j\neq i}^{k_t}(\lambda_{t,j}-\lambda)}{\prod_{i=1}^{k_t}(\lambda_{t,i}-\lambda)}.
\end{align}  

By examining the coefficients of the characteristic polynomial, we identify the leading coefficients:  
\[
m_1 =(-1)^{k_t}\left(1-\frac{1}{C_1}\sum_{i=k_t+1}^{d}\kappa_{t,i}^2\right),
\]  
and  
\[
m_2 = (-1)^{k_t-1}\left[\left(1-\frac{1}{C_1}\sum_{i=k_t+1}^{d}\kappa_{t,i}^2\right)\cdot \sum_{i=1}^{k_t}\lambda_i+\sum_{i=1}^{k_t}\kappa_{t,i}^2\right].
\]  

Hence, the sum of the auxiliary roots satisfies  
\begin{align}
\sum_{i=1}^{k_t}\widetilde{\lambda}_{t+1,i} 
= -\frac{m_2}{m_1} 
= \sum_{i=1}^{k_t}\lambda_i + \frac{\sum_{i=1}^{k_t}\kappa_{t,i}^2}{1-\frac{1}{C_1}\sum_{i=k_t+1}^{d}\kappa_{t,i}^2}
\;\leq\; \sum_{i=1}^{k_t}\lambda_i + \frac{1-C_2}{1-\frac{C_2}{C_1}}.
\end{align}  
Combining this with (\ref{equ:upper-bound-large}), we deduce that for all $i\leq j_t$,  
\[
\lambda_{t+1,i}\leq\widetilde{\lambda}_{t+1,i}.
\]  
Therefore, the total mass of the top $j_t$ eigenvalues after the update is bounded as  
\begin{align}
\sum_{i=1}^{j_t}\lambda_{t+1,i} 
&\leq \sum_{i=1}^{j_t}\widetilde{\lambda}_{t+1,i} 
= \sum_{i=1}^{j_t}\lambda_{t,i} + \sum_{i=1}^{j_t}(\widetilde{\lambda}_{t+1,i} - \lambda_{t,i}) \leq \sum_{i=1}^{j_t}\lambda_{t,i} + \frac{1-C_2}{1-\frac{C_2}{C_1}} \nonumber\\
&\leq \sum_{i=1}^{j_t}\lambda_{t,i} + \frac{1-C_2}{1-\frac{C_2}{2}} \leq \sum_{i=1}^{j_t}\lambda_{t,i} + 1-\frac{C_2}{2}.
\end{align}  

Consequently, the remaining eigenvalues necessarily gain at least a fixed amount of mass:  
\begin{align}
\sum_{i=j_t+1}^{d}\lambda_{t+1,i}\;\;\geq\;\;\sum_{i=j_t+1}^{d}\lambda_{t,i} + \frac{C_2}{2}.
\end{align}  

\paragraph{Step 3: lower bound on the smallest eigenvalue.}  
We now use the results from the previous parts to establish a quantitative lower bound for $\lambda_{t,d}$. Recall that in Step~2 we showed that for the set  
\[
\mathcal{L}_t = \{\, i:\lambda_{t,i}\leq 2C_1\lambda_{t,d}\,\},
\]  
the projection of $\bm{a}_t$ satisfies  
\begin{align*}
    \sum_{i\in\mathcal{L}_t}\kappa_{t,i}^2 \;\;\geq\;\; \frac{C_2}{2}.
\end{align*}

To track the cumulative effect across time, we extend this notation. For each $t\leq t'$, define  
\[
\overline{\mathcal{L}}_{t,t'} = \{\, i:\lambda_{t,i}\leq 2C_1\lambda_{t',d}\,\}.
\]  
Clearly, $\mathcal{L}_t \subseteq \overline{\mathcal{L}}_{t,t'}$. Next, define $A_t$ to be the total gain accumulated by the eigenvalues in $\mathcal{L}_s$ up to time $t$:  
\begin{align}
\label{equ:A-t-lower}
    A_t = \sum_{s=0}^{t-1}\sum_{i\in\mathcal{L}_s}\bigl(\lambda_{s+1,i}-\lambda_{s,i}\bigr).
\end{align}  
From Step~2, each summand contributes at least $C_2/2$, which yields the lower bound  
\[
A_t \;\;\geq\;\; \sum_{s=0}^{t-1}\frac{C_2}{2} = \frac{C_2 t}{2}.
\]

On the other hand, we can upper bound $A_t$ using the enlarged sets $\overline{\mathcal{L}}_{s,t}$:  
\begin{align}
    A_t &\leq \sum_{s=0}^{t-1}\sum_{i\in\overline{\mathcal{L}}_{s,t}}\bigl(\lambda_{s+1,i}-\lambda_{s,i}\bigr) \nonumber\\
    &= \sum_{s=0}^{t-1}\sum_{i=1}^{d} \bigl(\lambda_{s+1,i}-\lambda_{s,i}\bigr)\,\mathbf{1}\{\lambda_{s,i}\leq 2C_1\lambda_{t,d}\}\nonumber\\
    &= \sum_{i=1}^{d}\sum_{s=0}^{t-1} \bigl(\lambda_{s+1,i}-\lambda_{s,i}\bigr)\,\mathbf{1}\{\lambda_{s,i}\leq 2C_1\lambda_{t,d}\}.
\end{align}  
Since each $\lambda_{s,i}$ is monotone increasing, the inner summation is a telescoping sum. Moreover, whenever $\lambda_{s,i}\leq 2C_1\lambda_{t,d}$ we have  
\[
\lambda_{s+1,i}\leq \lambda_{s,i}+1 \leq 2C_1\lambda_{t,d}+1.
\]  
Thus the total increase of each eigenvalue is bounded by $2C_1\lambda_{t,d}$, recalling that $\bm{\Lambda}_0=\bm{I}_d$ (so $\lambda_{0,i}=1$ for all $i$). Therefore,  
\begin{align}
\label{equ:A-t-upper}
    A_t \;\leq\; \sum_{i=1}^{d} 2C_1\lambda_{t,d} \;\leq\; 2C_1 d\,\lambda_{t,d}.
\end{align}

Combining the lower bound (\ref{equ:A-t-lower}) with the upper bound (\ref{equ:A-t-upper}), we obtain  
\[
\lambda_{t,d}\;\;\geq\;\;\frac{C_2 t}{4C_1 d}.
\]  
Therefore, we have shown that whenever $\beta/\sqrt{\lambda_{t,d}}\geq 1/2$, there exists a universal constant $C>0$ such that  
\[
\lambda_{t,d}\;\;\geq\;\; C\cdot\frac{t}{d}.
\]  

\paragraph{Step 4: growth of the leading eigenvalue.}  
Define the stopping time  
\[
t_1' = \min\left\{t:\;\frac{\beta}{\sqrt{\lambda_{t,d}}}\leq C'\right\},
\]  
 for some constant $C'$. From the conclusion of Step~3, it follows immediately that $t_1' = \Theta(\beta^2 d)$. For any $t\geq t_1'$, with probability $1-1/T$, the estimation error can be bounded as  
\begin{align}
\label{equ:con-action-initial}
    \Vert\bm{a}_t-\bm{\theta}^{\star}\Vert_2& \leq \left\Vert\bm{a}_t-\bm{\widehat{\theta}}_{t}\right\Vert_2 + \left\Vert\bm{\widehat{\theta}}_{t}-\bm{\theta}^{\star}\right\Vert_2\lesssim\frac{\beta}{\sqrt{\lambda_{t,d}}} + \frac{\beta}{\sqrt{\lambda_{t,d}}}\leq \frac{3}{5},
\end{align}
when $C'$ is set small enough, where the second inequality holds from the confidence set in (\ref{eq:linucb-confidence-sets}), where we set $\delta = 1/T$.
As a consequence, for all $t>t_1'$ we obtain the correlation guarantee  
\[
\langle\bm{a}_t,\bm{\theta}^{\star}\rangle \;\;\geq\;\; \sqrt{1-\left(\frac{3}{5}\right)^2} \;=\; \frac{4}{5}.
\]  
Hence, setting $t_1 = 6t_1'$, we deduce  
\begin{align}
    {\bm{\theta}^{\star}}^{\!T}\bm{\Lambda}_{t_1}\bm{\theta}^{\star} 
    = \sum_{s=t_1'+1}^{t_1}\langle \bm{\theta}^{\star},\bm{a}_t\rangle^2 \geq \sum_{s=t_1'+1}^{t_1}\left(\frac{4}{5}\right)^2 \geq \frac{8}{15}\,t_1,
\end{align}
which implies  
\[
\lambda_{t_1,1} \;\;\geq\;\; \frac{8}{15}t_1.
\]  
Therefore, the leading eigenvalue separates from the rest by  
\[
\lambda_{t_1,1}-\lambda_{t_1,2}\;\;\geq\;\;\frac{t_1}{15}.
\]  



\subsection{Analysis of Phase \#2 (proof of Proposition \ref{prop:second-stage})}

We reparameterize $\lambda_{t,d}$ as follows: define $c_t$ such that
\begin{align}
    \lambda_{t,d} = c_t\beta\sqrt{t}.
\end{align}
We further introduce the notation
\[
\underline{c}_t = \min_{t_1\leq s\leq t} c_s,
\]
and define $\overline{\lambda}_t$ as the empirical average of the non-leading eigenvalues:
\begin{align}
    \overline{\lambda}_t = \frac{1}{d-1}\sum_{i=2}^{d}\lambda_{t,i}.
\end{align}
The key to establishing the desired result is to show that there exists a constant $c^{\star}>0$ such that
\[
0<c^{\star}<\underline{c}_t \qquad \text{for all } t>t_1.
\]

\paragraph{Step 1: characterize the distance between $\bm{v}_{t,1}$ and $\bm{\theta}^{\star}$.}  
We begin with the standard confidence bound on the estimation error:  
\begin{align}     
\Vert\bm{a}_t - \bm{\widehat{\theta}}_t\Vert_2 \lesssim \frac{\beta}{\sqrt{\lambda_{t,d}}}.  
\end{align} 
This shows that the action $\bm{a}_t$ chosen at time $t$ is close to the current estimate $\bm{\widehat{\theta}}_t$, with the error shrinking as the smallest eigenvalue $\lambda_{t,d}$ grows. By the triangle inequality, this further implies closeness to the true parameter: with probability $1-1/T$, 
\begin{align}     
\Vert\bm{a}_t - \bm{\theta}^{\star}\Vert_2
    &\lesssim \frac{\beta}{\sqrt{\lambda_{t,d}}} + \frac{\beta}{\sqrt{\lambda_{t,d}}} \lesssim \frac{\beta}{\sqrt{\lambda_{t,d}}},
\end{align} 
where the first inequality follow directly from (\ref{equ:con-action-initial}). Intuitively, this condition ensures that the statistical noise is dominated by the exploration parameter $\beta$, keeping $\bm{a}_t$ well aligned with $\bm{\theta}^{\star}$.  Next, consider the magnitude of the ``transformed signal'' $\bm{\Lambda}_t^{1/2}\bm{\theta}^{\star}$. Its squared norm is the cumulative signal energy collected along $\bm{\theta}^{\star}$:  
\begin{align}     
\left\Vert\bm{\Lambda}_t^{1/2}\bm{\theta}^{\star}\right\Vert_2      
    &= \sqrt{{\bm{\theta}^{\star}}^{T}\bm{\Lambda}_t\bm{\theta}^{\star}}     
    = \sqrt{\sum_{s=1}^{t}\langle \bm{a}_s,\bm{\theta}^{\star}\rangle^2} \nonumber\\     
    &= \sqrt{\sum_{s=1}^{t}\left(1 - O\left(\frac{\beta^2}{\lambda_{s,d}}\right)\right)}      
    = \sqrt{\sum_{s=1}^{t}\left(1 - O\left(\frac{\beta}{c_s\sqrt{s}}\right)\right)} \nonumber\\     
    &= \sqrt{t - O\!\left(\frac{\beta\sqrt{t}}{\underline{c}_t}\right)}.  
\end{align} 
Thus, the information collected in the direction of $\bm{\theta}^{\star}$ grows like $\sqrt{t}$, with only a mild correction due to imperfect exploration.  By comparison, the transformed leading eigenvector has energy bounded by the trace:  
\begin{align}      
\left\Vert\bm{\Lambda}_t^{1/2}\bm{v}_{t,1}\right\Vert_2       
    = \sqrt{\bm{v}_{t,1}^{T}\bm{\Lambda}_t\bm{v}_{t,1}}      
    \leq t+d.  
\end{align} 
This is a crude upper bound, but sufficient for our purposes.  We now turn to the eigenvalue structure. The largest eigenvalue satisfies  
\begin{align}     
\lambda_{t,1}\geq \left\Vert\bm{\Lambda}_t^{1/2}\bm{\theta}^{\star}\right\Vert_2^2
     = t-O\!\left(\frac{\beta\sqrt{t}}{\underline{c}_t}\right).  
\end{align} 


To relate $\bm{v}_{t,1}$ to $\bm{\theta}^{\star}$, expand $\bm{\theta}^{\star}$ in the eigenbasis:  
\begin{align} 
\left\Vert\bm{\Lambda}_t^{1/2}\bm{\theta}^{\star}\right\Vert_2^2 
&=\sum_{i=1}^d \lambda_{t,i}\,\langle \bm{\theta}^{\star},\bm{v}_{t,i}\rangle^2 \nonumber\\ 
&\le \lambda_{t,1}\langle \bm{\theta}^{\star},\bm{v}_{t,1}\rangle^2 
   + \lambda_{t,2}\left(1-\langle \bm{\theta}^{\star},\bm{v}_{t,1}\rangle^2\right)\nonumber\\ 
&\le \lambda_{t,2} + \left(\lambda_{t,1}-\lambda_{t,2}\right)\cdot \langle \bm{\theta}^{\star},\bm{v}_{t,1}\rangle^2.  
\end{align} 

Since $\operatorname{tr}(\bm{\Lambda}_t)=t+d$, the contribution from the smaller eigenvalues is limited:  
\begin{align} 
\sum_{i\ge 2}\lambda_{t,i} = t+d - \lambda_{t,1} = O\!\left(\frac{\beta\sqrt{t}}{\underline{c}_t}\right).  
\end{align} 
Thus, $\lambda_{t,2}\leq O\!\left(\frac{\beta\sqrt{t}}{\underline{c}_t}\right)$. It is also straightforward that $\lambda_{t,1}-\lambda_{t,2}\leq t$. Plugging this into the expansion gives  
\begin{align} 
\left\Vert\bm{\Lambda}_t^{1/2}\bm{\theta}^{\star}\right\Vert_2^2 
&\le  O\!\left(\frac{\beta\sqrt{t}}{\underline{c}_t}\right) + t\cdot \langle \bm{\theta}^{\star},\bm{v}_{t,1}\rangle^2.  
\end{align} 

On the other hand, we already have the lower bound  
\[
\left\Vert\bm{\Lambda}_t^{1/2}\bm{\theta}^{\star}\right\Vert_2^2 \ge t - O\!\left(\frac{\beta\sqrt{t}}{\underline{c}_t}\right).
\]  
Together, these inequalities imply  
\begin{align}     
t - O\!\left(\frac{\beta\sqrt{t}}{\underline{c}_t}\right) 
    \;\leq\; t\cdot \langle \bm{\theta}^{\star},\bm{v}_{t,1}\rangle^2 + O\!\left(\frac{\beta\sqrt{t}}{\underline{c}_t}\right).  
\end{align} 
Rearranging yields  
\begin{align}     
\langle\bm{\theta}^{\star},\bm{v}_{t,1}\rangle^2 \geq 1 - O\!\left(\frac{\beta}{\underline{c}_t\sqrt{t}}\right).  
\end{align}  

\medskip

Finally, for unit vectors it is well known that  
\begin{align*}     
\langle\bm{\theta}^{\star},\bm{v}_{t,1}\rangle^2 
    = 1- \Omega\!\left(\Vert \bm{v}_{t,1}-\bm{\theta}^{\star}\Vert_2^2\right).  
\end{align*} 
Hence the above inequality translates into the key control  
\begin{align} \label{equ:upper-bound-topeig}     
\Vert \bm{v}_{t,1}-\bm{\theta}^{\star}\Vert_2^2       
    = O\!\left(\frac{\beta}{\underline{c}_t\sqrt{t}}\right).  
\end{align}


\paragraph{Step 2: upper and lower bound the growth speed of non-leading eigenvalues.}  
We next analyze how the non-leading eigenvalues evolve over time. Establishing both upper and lower bounds on their growth is important: the upper bound ensures that they do not accumulate too much mass relative to the leading eigenvalue, while the lower bound guarantees that they still grow at a sufficient rate to prevent degeneration.  

\medskip

\noindent{\bf Upper bound.}  
By Lemma~\ref{lem:decom-a}, the decomposition coefficients satisfy
\begin{align}
    \nu_{t,1} = 1- O\!\left(\frac{\beta}{\underline{c}_t\sqrt{t}}\right),
    \qquad
    \kappa_{t,1} = 1- O\!\left(\frac{\beta}{\underline{c}_t\sqrt{t}}\right).
\end{align}
This shows that the majority of the mass concentrates in the leading eigen-direction, with only a small error proportional to $\frac{\beta}{\underline{c}_t\sqrt{t}}$.  Applying Lemma~\ref{lem:growth-largest-eigenvalue}, the leading eigenvalue evolves as
\begin{align}
    \lambda_{t+1,1} = \lambda_{t,1} + \kappa_{t,1}^2 + O\!\left(t^{-1}\right).
\end{align}
Because $\kappa_{t,1}^2 \approx 1$, this increment essentially captures the rate at which the leading eigenvalue dominates. Consequently, the mean of the non-leading eigenvalues evolves according to
\begin{align}
    \overline{\lambda}_{T} 
    = \overline{\lambda}_t + \frac{1}{d-1}\cdot O\!\left(\frac{\beta}{\underline{c}_t\sqrt{t}}\right).
\end{align}
Since $\underline{c}_t$ is non-increasing, telescoping this recursion yields the global upper bound
\begin{align}
\label{equ:upper-bound-mean-nonleading}
    \overline{\lambda}_t = O\!\left(\frac{\beta\sqrt{t}}{d\underline{c}_t}\right).
\end{align}
In other words, the mean non-leading eigenvalue cannot grow faster than $\sqrt{t}$, up to a factor depending on $\beta$ and the stability term $\underline{c}_t$.  

\medskip

\noindent{\bf Lower bound.}  
To complement the above, we construct a lower bound by considering the auxiliary optimization problem
\begin{align}
     \widetilde{\bm{w}}_t = \arg\max_{\|\bm{w}\|_2=1}\;
    \Bigl(1+\frac{\beta w_1}{\sqrt{\lambda_{t,1}}}\Bigr)^{\!2}
    + \sum_{i=2}^d \frac{\beta^2 w_i^2}{\lambda_{t,i}}.
\end{align}
This problem identifies the direction that maximizes the quadratic growth contribution, balancing the leading component with those from the non-leading directions.  Expanding the expression shows that this is equivalent to maximizing
\begin{align}
\max_{\Vert\bm{w}\Vert_2 = 1} \;& \frac{2\beta w_{t,1}}{\sqrt{\lambda_{t,1}}} 
- \beta^2 w_{t,1}^2\left(\frac{1}{\lambda_{t,d}} - \frac{1}{\lambda_{t,1}}\right).
\end{align}
Thus, the optimal weight on the leading coordinate, denoted $\tilde{w}_{t,1}^{\star}$, must satisfy
\begin{align}
    \tilde{w}_{t,1}^{\star} = \min\!\left(1,\; 
    \frac{1}{\beta\sqrt{\lambda_{t,1}}}\left(\frac{1}{\lambda_{t,d}} - \frac{1}{\lambda_{t,1}}\right)^{-1}\right).
\end{align}

In Phase~2, since $\lambda_{t,1}= \Omega(t)$ and the eigengap obeys $\lambda_{t,1}-\lambda_{t,d} = \Omega(\lambda_{t,1})$, this reduces to
\begin{align}
      \tilde{w}_{t,1}^{\star} 
      = \min\!\left(1,O\!\left(\frac{\lambda_{t,d}}{\beta\sqrt{t}}\right)\right)
      \;\leq\; O\!\left(\frac{\lambda_{t,d}}{\beta\sqrt{t}}\right).
\end{align}
Hence the mass allocated to the leading coordinate is negligible whenever $\lambda_{t,d}$ is small, meaning most weight shifts to non-leading directions.  Therefore, the contribution from the non-leading coordinates satisfies
\begin{align}
    \sum_{i=2}^{d}\frac{\beta^2 (\widetilde{w}_{t,i}^\star)^2}{\lambda_{t,i}} 
    =  \frac{\beta^2}{\lambda_{t,d}}\left(1-O\!\left(\frac{\lambda_{t,d}^2}{\beta^2 t}\right)\right).
\end{align}
By Lemma~\ref{lem:opt-comparison}, this translates to the inequality
\begin{align}
    \sum_{i=2}^{d}\left(\nu_{t,i}+\frac{\beta w_{t,i}^{\star}}{\sqrt{\lambda_{t,i}}}\right)^2
    \geq \frac{\beta^2\left(1-(w_{t,1}^{\star})^2\right)}{\lambda_{t,d}}
    \geq \frac{\beta^2\left(1-(\tilde{w}_{t,1}^{\star})^2\right)}{\lambda_{t,d}}.
\end{align}
Meanwhile, the total contribution is bounded by
\begin{align}
    \sum_{i=1}^{d}\left(\nu_{t,i}+\frac{\beta w_{t,i}^{\star}}{\sqrt{\lambda_{t,i}}}\right)^2
    \leq 1+\frac{\beta}{\sqrt{\lambda_{t,d}}} = O(1),
\end{align}
so it follows that
\begin{align}
    \sum_{i=2}^{d}\kappa_{t,i}^2 
    \gtrsim \frac{\beta^2}{\lambda_{t,d}}\left(1-O\!\left(\frac{\lambda_{t,d}^2}{\beta^2 t}\right)\right).
\end{align}
 
Finally, recalling that the average of the non-leading eigenvalues evolves as
\begin{align}
    \overline{\lambda}_{t+1} = \overline{\lambda}_{t} 
    + \frac{\sum_{i=2}^{d}\kappa_{t,i}^2}{d-1} + O\!\left(t^{-1}\right),
\end{align}
we conclude that whenever $\lambda_{t,d}\leq \widetilde{c}\beta\sqrt{t}$ (for some universal constant $\widetilde{c}$), the recursion satisfies
\begin{align}
     \overline{\lambda}_{t+1} = \overline{\lambda}_{t} 
     +\frac{1}{d-1}\cdot \Omega\!\left(\frac{\beta^2}{\lambda_{t,d}}\right) 
     = \overline{\lambda}_{t} +\frac{1}{d-1}\cdot \Omega\!\left(\frac{\beta}{c_t\sqrt{t}}\right).
\end{align}


\paragraph{Step 3: lower bound the growth of ``small eigenvalues".}

We now aim to establish a lower bound on the growth rate of the set of ``small eigenvalues’’ in the regime where $t$ belongs to the set  
\[
\mathcal{S} = \left\{t: c_t \leq \widetilde{c}, \ \frac{c_t}{\underline{c}_t} \leq 2 \right\},
\]  
that is, the regime in which the smallest eigenvalue risks dropping below the desired growth rate and has already crossed a critical threshold. Our proof strategy parallels that of Phase~\#1, but with a crucial modification. In this phase, we bound the growth rate of the small eigenvalues from below by a constant multiple of the growth rate of the mean of the non-leading eigenvalues, rather than by the mean of all eigenvalues. This adjustment is necessary because the leading eigenvalue $\lambda_{t,1}$ grows too rapidly to serve as a meaningful reference point in this regime.  

To proceed, consider the optimization problem
\begin{align}
& \max_{\Vert\bm{w}_t\Vert_2 = 1} \sum_{i=1}^{d}\left(\nu_{t,i} + \frac{\beta w_{t,i}}{\sqrt{\lambda_{t,i}}}\right)^2 \nonumber \\
&= \max_{w_1} \ \max_{\Vert\bm{w}_{-1}\Vert_2 = \sqrt{1-w_1^2}} \Bigg[ \left(\nu_{t,1}+\frac{\beta w_{t,1}}{\sqrt{\lambda_{t,1}}}\right)^2 + \sum_{i=2}^{d}\left(\nu_{t,i}+\frac{\beta w_{t,i}}{\sqrt{\lambda_{t,i}}}\right)^2 \Bigg]\nonumber\\
&= \max_{w_1} \left(\nu_{t,1}+\frac{\beta w_{t,1}}{\sqrt{\lambda_{t,1}}}\right)^2 + \max_{\Vert\bm{w}_{-1}\Vert_2 = \sqrt{1-w_1^2}} \sum_{i=2}^{d}\left(\nu_{t,i}+\frac{\beta w_{t,i}}{\sqrt{\lambda_{t,i}}}\right)^2. 
\end{align}
Thus, to characterize the contribution from $(w_{t,2}, \ldots, w_{t,d})$, it suffices to analyze the sub-optimization problem
\begin{align}
\max_{\Vert\bm{w}_{-1}\Vert_2 = \sqrt{1-w_1^2}} \sum_{i=2}^{d}\left(\nu_{t,i}+\frac{\beta w_{t,i}}{\sqrt{\lambda_{t,i}}}\right)^2.
\end{align}

Since $c_t/\underline{c}_t \leq 2$, it follows that
\begin{align}
\nu_{t,1} 
= 1 - O\left(\frac{\beta}{\underline{c}_t \sqrt{t}}\right) = 1 - O\left(\frac{\beta}{c_t \sqrt{t}}\right) = 1 - O\left(\frac{\beta^2}{\lambda_{t,d}}\right),
\end{align}
where the last step uses $\lambda_{t,d} = c_t \beta \sqrt{t}$. Consequently,
\begin{align}
\sqrt{\sum_{i=2}^{d} \nu_{t,i}^2} = O\left(\frac{\beta}{\sqrt{\lambda_{t,d}}}\right).
\end{align}

Returning to the optimization, we obtain
\begin{align}
\label{equ:upper-bound-opt}
& \left(\nu_{t,1}+\frac{\beta w_{t,1}}{\sqrt{\lambda_{t,1}}}\right)^2 
+ \max_{\Vert\bm{w}_{-1}\Vert_2 = \sqrt{1-w_1^2}} \sum_{i=2}^{d}\left(\nu_{t,i}+\frac{\beta w_{t,i}}{\sqrt{\lambda_{t,i}}}\right)^2 \nonumber\\
&\leq \left(\nu_{t,1}+\frac{\beta w_{t,1}}{\sqrt{\lambda_{t,1}}}\right)^2 
+ \left(\sqrt{\sum_{i=2}^{d}\nu_{t,i}^2} + \frac{\beta}{\sqrt{\lambda_{t,d}}}\sqrt{\sum_{i=2}^{d}w_{t,i}^2}\right)^2 \nonumber\\
&\leq 1 + O\left(\frac{\beta}{\sqrt{\lambda_{t,1}}}\right) + O\left(\frac{\beta^2}{\lambda_{t,d}}\right) \cdot \sqrt{\sum_{i=2}^{d}w_{t,i}^2}.
\end{align}

At the same time, the optimization also admits a simple lower bound:
\begin{align}
\label{equ:opt-lower-bound}
\max_{\Vert\bm{w}_t\Vert_2 = 1} \sum_{i=1}^{d}\left(\nu_{t,i}+\frac{\beta w_{t,i}}{\sqrt{\lambda_{t,i}}}\right)^2 \geq 1 + \frac{\beta^2}{\lambda_{t,d}}.
\end{align}

Since $\lambda_{t,1} \geq t/2$ and $\lambda_{t,d} = c_t \beta \sqrt{t}$, we deduce that when $c_t$ is sufficiently small (which can be achieved by choosing a sufficiently small constant $\widetilde{c}$),  it holds that
\begin{align}
\label{equ:bound-beta-lambda-1}
    O\left(\frac{\beta}{\sqrt{\lambda_{t,1}}}\right)\leq \frac{\beta^2}{2\lambda_{t,d}},
\end{align}
Plugging (\ref{equ:bound-beta-lambda-1}) into (\ref{equ:upper-bound-opt}) gives
\begin{align}
\label{equ:upper-bound-opt-2}
& \left(\nu_{t,1}+\frac{\beta w_{t,1}}{\sqrt{\lambda_{t,1}}}\right)^2 
+ \max_{\Vert\bm{w}_{-1}\Vert_2 = \sqrt{1-w_1^2}} \sum_{i=2}^{d}\left(\nu_{t,i}+\frac{\beta w_{t,i}}{\sqrt{\lambda_{t,i}}}\right)^2 \nonumber\\
&\leq 1 + \frac{\beta^2}{2\lambda_{t,d}} + O\left(\frac{\beta^2}{\lambda_{t,d}}\right) \cdot \sqrt{\sum_{i=2}^{d}w_{t,i}^2}.
\end{align}

Comparing the (\ref{equ:opt-lower-bound}) and (\ref{equ:upper-bound-opt-2}) ensures that there exists a universal constant $c_0 > 0$ such that
\begin{align}
\sqrt{\sum_{i=2}^{d} w_{t,i}^2} \geq c_0.
\end{align}

Therefore, by a direct application of Lemma \ref{lem:opt-concentration}, for some constant $C_1$, we define  
\[
k_t = \{i:\lambda_{t,i}\leq C_1\lambda_{t,d}\}.
\]  
In words, this set indexes the eigenvalues that are not ``too large'' compared to the smallest one. Then there exists a constant $C_2$ such that  
\begin{align}
    \sum_{i\in\mathcal{L}_t}\kappa_{t,i}^2\geq C_2\cdot \sum_{i=2}^{d}\kappa_{t,i}^2.
\end{align}  
This ensures that a nontrivial fraction of the total variance (as measured by the $\kappa_{t,i}^2$) is concentrated in the set $\mathcal{L}_t$.  We now proceed by following a strategy similar to the one employed in Phase 1. Specifically, we enlarge the set slightly and redefine  
\[
\mathcal{L}_t = \{i:\lambda_{t,i}\leq 2C_1\lambda_{t,d}\}.
\]  
Recall from Lemma \ref{lem:rank-one-update} that the updated eigenvalues $\lambda_{t+1,1},\ldots,\lambda_{t+1,d}$ are the roots of the secular equation  
\begin{align*}
    f(\lambda) = 1+\sum_{i=1}^{d}\frac{\kappa_{t,i}^2}{\lambda_{t,i}-\lambda}.
\end{align*}

To study the behavior of the largest eigenvalues (and control their possible growth), we introduce the auxiliary function  
\[
f_1(\lambda) = 1+\sum_{i=1}^{k_t}\frac{\kappa_{t,i}^2}{\lambda_{t,i}-\lambda},
\]  
which only accounts for the ``small'' eigenvalues indexed by $k_t$. The idea is that contributions from the larger eigenvalues (outside this set) can be bounded separately.  For any eigenvalue $\lambda_{t+1,i}$ with $i\leq j_t$, we have  
\begin{align}
    f(\lambda_{t+1,i})  
    = f_1(\lambda_{t+1,i}) + \sum_{i=k_t+1}^{d}\frac{\kappa_{t,i}^2}{\lambda_{t,i}-\lambda}
    \geq f_1(\lambda_{t+1,i})-\frac{1}{C_1\lambda_{t,d}}\sum_{i=k_t+1}^{d}\kappa_{t,i}^2.
\end{align}  
Here, the last inequality comes from the fact that each denominator is at least $C_1\lambda_{t,d}$. Thus, it follows that  
\begin{align} \label{equ:upper-bound-large-2}
    f_1(\lambda_{t+1,i})\leq \frac{1}{C_1}\sum_{i=k_t+1}^{d}\kappa_{t,i}^2.
\end{align}  
Intuitively, this means that the influence of the ``large'' eigenvalues (outside $k_t$) limits how big the roots of $f_1$ can be.  Next, define $\widetilde{\lambda}_{t,1},\ldots,\widetilde{\lambda}_{t,k_t}$ as the solutions of  
\begin{align}
    f_1(\lambda) =  \frac{1}{C_1}\sum_{i=k_t+1}^{d}\kappa_{t,i}^2.
\end{align}  
Equivalently, these values are the roots of the polynomial identity  
\begin{align}
    f_1(\lambda) - \frac{1}{C_1}\sum_{i=k_t+1}^{d}\kappa_{t,i}^2 
    & = 1-\frac{1}{C_1}\sum_{i=k_t+1}^{d}\kappa_{t,i}^2 + \sum_{i=1}^{k_t}\frac{\kappa_{t,i}^2}{\lambda_{t,i}-\lambda}\nonumber\\
    & = \frac{\left(1-\frac{1}{C_1}\sum_{i=k_t+1}^{d}\kappa_{t,i}^2\right)\cdot \prod_{i=1}^{k_t}(\lambda_{t,i}-\lambda) +\sum_{i=1}^{k_t} \kappa_{t,i}^2 \prod_{j\neq i}^{k_t}(\lambda_{t,j}-\lambda)}{\prod_{i=1}^{k_t}(\lambda_{t,i}-\lambda)}.
\end{align}  

By examining the coefficients of this polynomial (via Vieta’s formulas), we obtain:  
\[
m_1 =(-1)^{k_t}\left(1-\frac{1}{C_1}\sum_{i=k_t+1}^{d}\kappa_{t,i}^2\right),
\]  
\[
m_2 = (-1)^{k_t-1}\left[\left(1-\frac{1}{C_1}\sum_{i=k_t+1}^{d}\kappa_{t,i}^2\right)\cdot \sum_{i=1}^{k_t}\lambda_i+\sum_{i=1}^{k_t}\kappa_{t,i}^2\right].
\]  

This gives the bound  
\begin{align}
    \sum_{i=1}^{k_t}\widetilde{\lambda}_{t+1,k_t} = -\frac{m_2}{m_1} 
    = \sum_{i=1}^{k_t}\lambda_i + \frac{\sum_{i=1}^{k_t}\kappa_{t,i}^2}{1-\frac{1}{C_1}\sum_{i=k_t+1}^{d}\kappa_{t,i}^2}
    \leq  \sum_{i=1}^{k_t}\lambda_i + \frac{1-C_2\cdot \sum_{i=2}^{d}\kappa_{t,i}^2}{1-\frac{C_2}{C_1}\cdot\sum_{i=2}^{d}\kappa_{t,i}^2}.
\end{align}  
This shows that the sum of these ``controlled'' eigenvalues cannot grow too quickly, since the denominator penalizes large contributions from $\sum_{i=2}^{d}\kappa_{t,i}^2$.  From (\ref{equ:upper-bound-large-2}), we know that for all $i\leq j_t$,  
\[
\lambda_{t+1,i}\leq\widetilde{\lambda}_{t+1,i}.
\]  
Hence, the sum of the first $j_t$ eigenvalues of $\bm{\Lambda}_{T}$ can be bounded by  
\begin{align}
    \sum_{i=1}^{j_t}\lambda_{t+1,i}  
    &\leq  \sum_{i=1}^{j_t}\widetilde{\lambda}_{t+1,i} 
    =\sum_{i=1}^{j_t}\lambda_{t,i} + \sum_{i=1}^{j_t}(\widetilde{\lambda}_{t+1,i} - \lambda_{t,i})\nonumber\\
    & \leq \sum_{i=1}^{j_t}\lambda_{t,i} + \frac{1-C_2}{1-\frac{C_2}{C_1}}
    \leq \sum_{i=1}^{j_t}\lambda_{t,i} + \frac{1-C_2\sum_{i=2}^{d}\kappa_{t,i}^2}{1-\frac{C_2}{2}\cdot \sum_{i=2}^{d}\kappa_{t,i}^2}\nonumber\\
    &\leq  \sum_{i=1}^{j_t}\lambda_{t,i} + 1-\frac{C_2}{2}\cdot \sum_{i=2}^{d}\kappa_{t,i}^2.
\end{align}  

The key takeaway is that the leading eigenvalues are essentially ``capped'' in their growth: they can increase by at most a bounded additive amount, while the remaining eigenvalues must absorb a proportional share of the increase.  This balance becomes explicit when we look at the complement:  
\begin{align}
    \sum_{i=j_t+1}^{d}\lambda_{t+1,i}\geq \sum_{i=j_t+1}^{d}\lambda_{t,i} + \frac{C_2}{2}\cdot \sum_{i=2}^{d}\kappa_{t,i}^2.
\end{align}  
That is, the ``smaller'' eigenvalues (those beyond $j_t$) are guaranteed to grow by a nontrivial amount.  Moreover, since  
\begin{align}
    \lambda_{t+1,1}\geq \lambda_{t,1}+\kappa_{t,1}^2,
\end{align}  
we also have  
\begin{align}
    \sum_{i=2}^{d}\lambda_{t+1,i}\leq \sum_{i=2}^{d}\lambda_{t,i}+\sum_{i=2}^{d}\kappa_{t,i}^2.
\end{align}

Putting everything together, we finally obtain  
\begin{align}
    \sum_{i\in\mathcal{L}_t}(\lambda_{t+1,i}-\lambda_{t,i})\geq\frac{C_2}{2}\cdot  \sum_{i=2}^{d}\kappa_{t,i}^2\geq \frac{C_2}{2}\cdot \sum_{i=2}^{d}(\lambda_{t+1,i}-\lambda_{t,i}).
\end{align}  


\paragraph{Step 4: lower bound the minimum eigenvalue.}
We are now in position to establish the main result of this phase. Specifically, we show that when $\underline{c}_t$ is sufficiently small and throughout any period in which $c_t / \underline{c}_t \leq 2$, the minimum eigenvalue $\lambda_{t,d}$ grows at least on the order of $O(\beta \sqrt{t}/\sqrt{d})$. In the previous step we observed that the threshold $\widetilde{c}$ can be chosen arbitrarily small. We therefore set $\widetilde{c} = O(d^{-1/2})$ and define
\[
\mathcal{S} = \left\{ t : c_t \leq \widetilde{c},\; \frac{c_t}{\underline{c}_t} \leq 2 \right\}.
\]
By construction, $\mathcal{S}$ is a finite union of \emph{contiguous integer intervals} (segments) of time. Our objective is to obtain a uniform lower bound on $c_t$ (and hence on $\lambda_{t,d}$); it suffices to analyze $c_t$ on each segment in $\mathcal{S}$ separately. Throughout, recall that $\overline{\lambda}_t$ denotes the average of the nonleading eigenvalues (i.e., $\overline{\lambda}_t = \frac{1}{d-1}\sum_{i=2}^d \lambda_{t,i}$) and that $C_1,C_2,C_3>0$ are absolute constants independent of $t$ and $d$.

Fix any segment and let $t' < t''$ be its endpoints. For any $t \in [t',t'']$, we compare lower and upper bounds on the cumulative change of the $\lambda_{t,i}$. Summing the growth of all but the smallest eigenvalue yields
\begin{align} \label{equ:mineig-lower}
\sum_{t=t'}^{t''-1} \sum_{i:\lambda_{t,i} \in \mathcal{L}_t} (\lambda_{t+1,i} - \lambda_{t,i})
&\geq C_2 \sum_{t=t'}^{t''-1} \sum_{i=2}^d (\lambda_{t+1,i} - \lambda_{t,i}) \notag\\
&= C_2 \sum_{i=2}^d (\lambda_{t'',i} - \lambda_{t',i}) \notag\\
&= C_2 (d-1) (\overline{\lambda}_{t''} - \overline{\lambda}_{t'}).
\end{align}
For the matching upper bound, define
\[
\widetilde{\mathcal{L}}_{t,t''} = \left\{ i : \lambda_{t,i} \leq C_1 \lambda_{t'',d} \right\},
\]
and apply the same counting argument used in Phase \#1 (controlling how many coordinates can exceed a multiple of the minimum). This gives
\begin{align} \label{equ:mineig-upper}
\sum_{t=t'}^{t''-1} \sum_{i:\lambda_{t,i} \in \mathcal{L}_t} (\lambda_{t+1,i} - \lambda_{t,i})
\leq \sum_{t=t'}^{t''-1} \sum_{i \in \widetilde{\mathcal{L}}_{t,t''}} (\lambda_{t+1,i} - \lambda_{t,i})
\leq (d-1)\left(C_1 \lambda_{t'',d} + 1 - \lambda_{t',d}\right).
\end{align}
Combining (\ref{equ:mineig-lower}) and (\ref{equ:mineig-upper}) yields
\[
C_2 (\overline{\lambda}_{t''} - \overline{\lambda}_{t'}) \leq C_1 \lambda_{t'',d} + 1 - \lambda_{t',d},
\]
and hence
\begin{align} \label{equ:min-lower-bound}
\lambda_{t'',d} \geq \min\left( \lambda_{t',d},\; \frac{C_2}{C_1}(\overline{\lambda}_{t''} - \overline{\lambda}_{t'}) + \frac{\lambda_{t',d} - 1}{C_1} \right).
\end{align}

We next control the minimum of $c_t$ on each segment of $\mathcal{S}$. Without loss of generality assume $t_1 \notin \mathcal{S}$ (this can be arranged by taking $\widetilde{c}$ below a fixed absolute constant), so every segment begins at some $t' > t_0$. There are two ways a new segment can start:
\begin{itemize}
\item[(i)] $c_{t'} \leq \widetilde{c}$ while $c_{t'-1} > \widetilde{c}$. By the $O(t^{-1/2})$ drift established earlier for $c_t$ (smooth variation), the threshold crossing cannot overshoot by more than a fixed fraction for large $t'$, hence $c_{t'} \geq 0.95 \widetilde{c}$.
\item[(ii)] $c_{t'}/\underline{c}_{t'} \leq 2$ while $c_{t'-1}/\underline{c}_{t'-1} > 2$. The same smoothness argument applied to the ratio shows $c_{t'}/\underline{c}_{t'} > 1.9$, so $c_{t'} \geq 1.9 \underline{c}_{t'}$ for large $t'$.
\end{itemize}
Combining (i)–(ii) gives the useful entry condition
\[
c_{t'} \geq \min\left(1.9 \underline{c}_{t'},\; 0.95 \widetilde{c} \right).
\]

From the growth bounds of the preceding steps, for all $t$ on the same segment as $t'$ we have
\[
\overline{\lambda}_{T} \geq \overline{\lambda}_t + \frac{C_3}{d-1} \cdot \frac{\beta}{2 \underline{c}_{t'} \sqrt{t}}
\geq \overline{\lambda}_t + \frac{C_3}{d-1} \cdot \frac{\beta}{2 \widetilde{c} \sqrt{t}},
\]
and summing from $t'$ to $t''$ yields
\[
\overline{\lambda}_{t''} - \overline{\lambda}_{t'} \geq \frac{C_3 \beta}{\widetilde{c}(d-1)}(\sqrt{t''} - \sqrt{t'}).
\]
Plugging this into (\ref{equ:min-lower-bound}) and then applying it with $t'' = t$ (for any $t \in [t',t'']$) gives
\begin{align}
\lambda_{t,d} &\geq \max\left( \lambda_{t',d},\; \frac{C_2 C_3 \beta}{C_1 \widetilde{c}(d-1)}(\sqrt{t} - \sqrt{t'}) + \frac{c_{t'} \beta \sqrt{t'} - 1}{C_1} \right) \notag\\
&= \max\left( \lambda_{t',d},\; c \beta (\sqrt{t} - \sqrt{t'}) + \frac{\min(1.9 \underline{c}_{t'}, 0.95 \widetilde{c}) \beta \sqrt{t'} - 1}{C_1} \right),
\end{align}
where we set
\[
c := \frac{C_2 C_3}{C_1 \widetilde{c} (d-1)} = O(d^{-1/2}).
\]

We now choose $c^\star = O(d^{-1/2})$ with $c^\star \leq 0.5\,\widetilde{c}$ such that, for all $t \geq t'$, 
\[
c^\star \beta \sqrt{t} \leq \max\left( 1.9 c^\star \beta \sqrt{t'},\; c \beta (\sqrt{t} - \sqrt{t'}) + \frac{1.9 c^\star \beta \sqrt{t'} - 1}{C_1} \right).
\]
Such a choice is always possible: the right-hand side is the maximum of two affine functions of $\sqrt{t}$ whose slopes are $0$ and $c\beta>0$, respectively, whereas the left-hand side has slope $c^\star\beta$; taking $c^\star \leq c$ and adjusting the intercept via the $-(1/C_1)$ term ensures the inequality holds for all $t \ge t'$. Therefore, for all $t > t'$ on the same segment in $\mathcal{S}$,
\begin{align}
\lambda_{t,d} &\geq \max\left( \min(1.9 \underline{c}_{t'}, 0.95 \widetilde{c}) \beta \sqrt{t'},\; c \beta (\sqrt{t} - \sqrt{t'}) + \frac{\min(1.9 \underline{c}_{t'}, 0.95 \widetilde{c}) \beta \sqrt{t'} - 1}{C_1} \right) \notag\\
&\geq \max\left( 1.9 c^\star \beta \sqrt{t'},\; c \beta (\sqrt{t} - \sqrt{t'}) + \frac{1.9 c^\star \beta \sqrt{t'} - 1}{C_1} \right) \notag\\
&\geq c^\star \beta \sqrt{t}.
\end{align}
Applying this argument on every segment of $\mathcal{S}$, we conclude that
\[
\lambda_{t,d} \geq c^\star \beta \sqrt{t}\asymp\frac{\beta\sqrt{t}}{\sqrt{d}}, \quad \text{for all } t \geq t_1.
\]
Finally, from equations~(\ref{equ:upper-bound-topeig}) and~(\ref{equ:upper-bound-mean-nonleading}),
\[
\left\| \bm{v}_{t,1} - \bm{\theta}^\star \right\|_2^2 = O\left( \frac{\beta \sqrt{d}}{\sqrt{t}} \right), \quad \text{and} \quad \overline{\lambda}_t = O\left( \frac{\beta \sqrt{t}}{\sqrt{d}} \right).
\]
Combining these with the lower bound on $\lambda_{t,d}$ established above gives
\begin{align}
\lambda_{t,d} \asymp \frac{\beta \sqrt{t}}{\sqrt{d}},
\end{align}
and hence the estimation error of the leading eigenvector satisfies
\begin{align}
\left\| \bm{v}_{t,1} - \bm{\theta}^\star \right\|_2 = O\left( \frac{\beta}{\sqrt{\lambda_{t,d}}} \right).
\end{align}



\subsection{Analysis of Phase \#3 (proof of Proposition \ref{prop:third-stage})}
\paragraph{Step 1: characterizing the update of the next top eigenvector $\bm{v}_{t+1,1}$.}  
We aim to show in this phase that $\bm{v}_{t+1,1}$ can be characterized by the previous leading eigenvector $\bm{v}_{t,1}$ and the estimated parameter $\bm{\widehat{\theta}_t}$ as follows
\[
   \bm{v}_{t+1,1} \;=\; \bm{v}_{t,1} \;+\; \frac{\widehat{\bm{\theta}}_t - \bm{v}_{t,1}}{t} \;+\; \bm{\zeta}_t,
\]  
where $\bm{\zeta}_t$ is a high-order perturbation that can be nicely bounded both in norm and direction, being near orthogonal to $\bm{v}_{t,1}-\bm{\widehat{\theta}}_t$. We proceed with the following steps to show this result.

\noindent{\bf Lower bound on $|\kappa_{t,i}|$ ($i\geq 2$).} To make this precise, we first derive a uniform upper bound on  
\[
\bigl\|\widehat{\bm{\theta}}_t+\beta\bm{\Lambda}_t^{-1/2}\bm{w}\bigr\|_2,
\]  
valid for any unit vector $\bm{w}$, i.e., $\|\bm{w}\|_2=1$.  We first establish the uniform bound. Expanding the squared norm gives  
\begin{align}
    \bigl\|\widehat{\bm{\theta}}_t + \beta \bm{\Lambda}_t^{-1/2} \bm{w}\bigr\|_2^2
    \;=\; \sum_{i=1}^{d}\!\left(\nu_{t,i} + \frac{\beta w_{i}}{\sqrt{\lambda_{t,i}}}\right)^{\!2},
\end{align}
where $\nu_{t,i}$ denotes the $i$-th coordinate of $\widehat{\bm{\theta}}_t$ expressed in the eigenbasis of $\bm{\Lambda}_t$. For the \emph{leading coordinate} ($i=1$), our earlier estimates imply  
\[
    \nu_{t,1} + \frac{\beta w_1}{\sqrt{\lambda_{t,1}}}
    \;=\; 1 - O\!\left(\frac{\beta^2}{\lambda_{t,d}}\right)
          + O\!\left(\frac{\beta}{\sqrt{t}}\right)
    \;=\; 1 + O\!\left(\frac{\beta^2}{\lambda_{t,d}}\right),
\]
where we have used the bound $\lambda_{t,d}\lesssim \beta\sqrt{t}/\sqrt{d}$ in Proposition~\ref{prop:second-stage} to absorb the $O(\beta/\sqrt{t})$ term into the $O(\beta^2/\lambda_{t,d})$ error. This shows that the leading component remains close to $1$, with only a small perturbation.  

For the \emph{remaining coordinates} ($i\ge 2$), we control their contribution by  
\begin{align*}
    \sum_{i=2}^{d}\!\left(\nu_{t,i} + \frac{\beta w_{i}}{\sqrt{\lambda_{t,i}}}\right)^2
    \;\lesssim\; \sum_{i=2}^{d}\!\left(\nu_{t,i}^2 + \frac{\beta^2 w_i^2}{\lambda_{t,i}}\right) 
    &\;\le\; \bigl(1-\nu_{t,1}^2\bigr) + \frac{\beta^2}{\lambda_{t,d}}\sum_{i=2}^d w_i^2 \\
    &\;=\; O\!\left(\frac{\beta^2}{\lambda_{t,d}}\right),
\end{align*}
where we used $\sum_{i=1}^d w_i^2=1$ and the fact that $\nu_{t,1}^2 = 1 - O(\beta^2/\lambda_{t,d})$.  Putting the two pieces together, we obtain the \emph{uniform expansion}
\begin{align}
\label{equ:max-norm}
    \sum_{i=1}^{d}\!\left(\nu_{t,i} + \frac{\beta w_{i}}{\sqrt{\lambda_{t,i}}}\right)^2
    \;=\; 1 + O\!\left(\frac{\beta^2}{\lambda_{t,d}}\right),
\end{align}
valid uniformly over all unit vectors $\bm{w}$ and all $t \geq t_1$.  We now use (\ref{equ:max-norm}) to sharpen the characterization of the maximizer $\bm{a}_t$. Recall its definition:  
\[
    \bm{w}_t \;=\; \arg\max_{\|\bm{w}\|_2 = 1}\;
    \sum_{i=1}^{d} \left( \nu_{t,i} + \frac{\beta w_{i}}{\sqrt{\lambda_{t,i}}} \right)^2
    \;=\; \arg\max_{\|\bm{w}\|_2 = 1}\;\sum_{i=1}^{d}\left(\frac{2\beta \nu_{t,i} w_{t,i}}{\sqrt{\lambda_{t,i}}} +\frac{\beta^2 w_{t,i}^2}{\lambda_{t,i}}\right).
\]
The maximization is dominated by the cross-term  
\[
2\sum_{i=1}^d \frac{\beta}{\sqrt{\lambda_{t,i}}}\nu_{t,i}w_i,
\]  
which is maximized when $w_{t,i}$ aligns in sign with $\nu_{t,i}$. Hence, at the maximizer we necessarily have  
\[
w_{t,i}\nu_{t,i} \;\geq\; 0 \quad \text{for all } i.
\]

By (\ref{equ:max-norm}), we know that for any unit vector $\bm{w}$—in particular, for $\bm{w}=\bm{w}_t$—the perturbation vector satisfies  
\[
    \bigl\|\widehat{\bm{\theta}}_t + \beta \bm{\Lambda}_t^{-1/2} \bm{w}\bigr\|_2
    \;=\; 1 + O\!\left(\frac{\beta^2}{\lambda_{t,d}}\right).
\]
Consequently, the normalized projection operator satisfies  
\begin{align}
    \mathcal{P}\!\left(\widehat{\bm{\theta}}_t + \beta \bm{\Lambda}_t^{-1/2} \bm{w}\right)
    \;=\; \left(1 - O\!\left(\frac{\beta^2}{\lambda_{t,d}}\right)\right)
          \cdot\left(\widehat{\bm{\theta}}_t + \beta \bm{\Lambda}_t^{-1/2} \bm{w}\right).
\end{align}

For each coordinate $i\in[d]$, this implies the refined characterization  
\begin{align}
    \kappa_{t,i}
    \;=\; \left(1 - O\!\left(\frac{\beta^2}{\lambda_{t,d}}\right)\right)
          \left( \nu_{t,i} + \frac{\beta w_{t,i}}{\sqrt{\lambda_{t,i}}} \right).
\end{align}
Since $w_{t,i}\nu_{t,i}\ge 0$, the correction preserves the sign of the signal and enlarges its magnitude:  
\begin{align*}
    \left|\nu_{t,i}+\frac{\beta w_{t,i}}{\sqrt{\lambda_{t,i}}}\right|\;\geq\; |\nu_{t,i}|.
\end{align*}
Therefore,  
\begin{align}
    |\kappa_{t,i}| \;=\; \left(1 - O\!\left(\frac{\beta^2}{\lambda_{t,d}}\right)\right)\cdot
          \left| \nu_{t,i} + \frac{\beta w_{t,i}}{\sqrt{\lambda_{t,i}}} \right|
    \;\geq\; \left(1 - O\!\left(\frac{\beta^2}{\lambda_{t,d}}\right)\right)\,|\nu_{t,i}|.
\end{align}
In other words, the normalized contributions $\kappa_{t,i}$ not only preserve the alignment with the underlying signal but also lose at most an $O(\beta^2/\lambda_{t,d})$ fraction of their magnitude. This guarantees stability of the signal direction under the perturbation.



\noindent{\bf Spectral decomposition of next top eigenvector $\bm{v}_{t+1,1}$.} We now decompose the next top eigenvector $\bm{v}_{t+1,1}$ in the eigenbasis formed by the previous eigenvectors $\{\bm{v}_{t,1},\ldots,\bm{v}_{t,d}\}$. That is, we write
\begin{align}
\bm{v}_{t+1,1} = \sum_{i=1}^d \omega_{t,i} \bm{v}_{t,i}.
\end{align}
Our goal is to characterize the coefficients $\{\omega_{t,i}\}$, in particular establishing a nontrivial lower bound for $\omega_{t,i}$ when $i \geq 2$.  From Lemma~\ref{lem:rank-one-update}, the top eigenvector after the rank-one perturbation can be expressed explicitly as  
\begin{align}
\bm{v}_{t+1,1} = K_t^{-1} \sum_{i=1}^d \frac{\kappa_{t,i}}{\lambda_{t+1,1} - \lambda_{t,i}} \bm{v}_{t,i}, \label{eq:vt+1-decomp}
\end{align}
where $K_t$ is a normalization constant ensuring $\|\bm{v}_{t+1,1}\|_2 = 1$. In particular,
\[
K_t = \left( \sum_{i=1}^d \left( \frac{\kappa_{t,i}}{\lambda_{t+1,1} - \lambda_{t,i}} \right)^2 \right)^{1/2}.
\]
Intuitively, the denominator $\lambda_{t+1,1} - \lambda_{t,i}$ captures the spectral separation between the new leading eigenvalue $\lambda_{t+1,1}$ and the old eigenvalues $\{\lambda_{t,i}\}$, while the numerator $\kappa_{t,i}$ measures the alignment between the perturbation direction and $\bm{v}_{t,i}$. Thus, the size of each $\omega_{t,i}$ is governed both by spectral gaps and by how much the perturbation projects onto $\bm{v}_{t,i}$.  

We first control the normalization factor $K_t$. From Lemma~\ref{lem:growth-largest-eigenvalue}, the increment of the leading eigenvalue satisfies
\[
\lambda_{t+1,1} - \lambda_{t,1} = 1 + O\!\left(\frac{\beta^2}{\lambda_{t,d}}\right).
\]
Moreover, for $i \geq 2$, the spectral gap is lower bounded as
\[
\lambda_{t+1,1} - \lambda_{t,i} \geq \lambda_{t,1} - \lambda_{t,i} \gtrsim t,
\]
thanks to the assumed eigengap structure at time $t$. Substituting the bounds on $\kappa_{t,i}$ from Lemma~\ref{lem:decom-a}, we obtain
\begin{align}
K_t^2 &= \left( \frac{\kappa_{t,1}}{\lambda_{t+1,1} - \lambda_{t,1}} \right)^2 + \sum_{i=2}^d \left( \frac{\kappa_{t,i}}{\lambda_{t+1,1} - \lambda_{t,i}} \right)^2 \nonumber\\
&= \left(1 - O\!\left(\frac{\beta^2}{\lambda_{t,d}}\right)\right)^2 + O(t^{-2}) \cdot O\!\left(\frac{\beta^2}{\lambda_{t,d}}\right) \nonumber\\
&= 1 - O\!\left(\frac{\beta^2}{\lambda_{t,d}}\right), \label{eq:Kt}
\end{align}
which implies $K_t = 1 - O(\beta^2/\lambda_{t,d})$.  Now extracting coefficients from (\ref{eq:vt+1-decomp}), we see that for $i \geq 2$,
\begin{align}
|\omega_{t,i}| &= K_t^{-1} \cdot \frac{|\kappa_{t,i}|}{\lambda_{t+1,1} - \lambda_{t,i}} \nonumber\\
&\geq \left(1 + O\!\left(\frac{\beta^2}{\lambda_{t,d}}\right)\right) \cdot \frac{\bigl(1 - O(\frac{\beta^2}{\lambda_{t,d}})\bigr)|\nu_{t,i}|}{t} \nonumber\\
&\geq \left(1 - O\!\left(\frac{\beta^2}{\lambda_{t,d}}\right)\right) \cdot \frac{|\nu_{t,i}|}{t}. \label{eq:omega-ti}
\end{align}
Thus, a clean lower bound is established for all $\omega_{t,i}$ with $i \geq 2$. Since $\kappa_{t,i}\nu_{t,i} \geq 0$ and $\omega_{t,i}$ inherits the sign of $\kappa_{t,i}$, we further deduce that
\begin{align}
\omega_{t,i}\nu_{t,i} \geq 0, \qquad \text{for all } i \geq 2.
\end{align}  

Turning to the leading coefficient $\omega_{t,1}$, observe that
\begin{align}
\sum_{i=2}^d \omega_{t,i}^2 = K_t^{-2} \cdot \sum_{i=2}^d \left( \frac{\kappa_{t,i}}{\lambda_{t+1,1} - \lambda_{t,i}} \right)^2 \lesssim \Bigl(1 + O\!\left(\frac{\beta^2}{\lambda_{t,d}}\right)\Bigr) \cdot \frac{O(\beta^2/\lambda_{t,d})}{t^2} = O\!\left(\frac{\beta^2}{t^2 \lambda_{t,d}}\right).
\end{align}
Hence,
\begin{align}
\omega_{t,1} = \sqrt{1 - \sum_{i=2}^d \omega_{t,i}^2} = 1 - O\!\left(\frac{\beta^2}{t^2 \lambda_{t,d}}\right).
\end{align}

\noindent{\bf Characterize the norm and direction of $\bm{\zeta}_t$.} To obtain a recursive characterization of $\bm{v}_{t+1,1}$, it remains to compare it against the “linearized” update in the direction of $\bm{\widehat{\theta}}_t$. Specifically, define the intermediate vector
\begin{align}\label{def:overline-v}
\bm{\overline{v}}_{t+1,1} := \bm{v}_{t,1} + \frac{\bm{\widehat{\theta}}_t - \bm{v}_{t,1}}{t},
\end{align}
which lies on the line segment connecting $\bm{v}_{t,1}$ and $\bm{\widehat{\theta}}_t$. Expanding this vector in the eigenbasis $\{\bm{v}_{t,1}, \ldots, \bm{v}_{t,d}\}$ gives
\begin{align}
\bm{\overline{v}}_{t+1,1} = \left(1 - \frac{1 - \nu_{t,1}}{t}\right) \bm{v}_{t,1} + \sum_{i=2}^{d} \frac{\nu_{t,i}}{t} \bm{v}_{t,i},
\end{align}
where $\bm{\widehat{\theta}}_t = \sum_{i=1}^{d} \nu_{t,i} \bm{v}_{t,i}$ is the eigenbasis decomposition of $\bm{\widehat{\theta}}_t$.  

We now measure the discrepancy between $\bm{v}_{t+1,1}$ and this linearized vector:
\begin{align}
\left\|\bm{v}_{t+1,1} - \bm{\overline{v}}_{t+1,1}\right\|_2^2 &= \left(\omega_{t,1} - 1 + \frac{1 - \nu_{t,1}}{t}\right)^2 + \sum_{i=2}^{d} \left(\omega_{t,i} - \frac{\nu_{t,i}}{t}\right)^2. \nonumber
\end{align}
Using the bounds established earlier on $\omega_{t,1}$ and $\omega_{t,i}$, we obtain
\begin{align}
\left\|\bm{v}_{t+1,1} - \bm{\overline{v}}_{t+1,1}\right\|_2^2 
&\leq \left(O\!\left(\frac{\beta^2}{t^2 \lambda_{t,d}}\right) + O\!\left(\frac{ \beta^2}{t \lambda_{t,d}}\right)\right)^2 
+ \sum_{i=2}^{d} \omega_{t,i}^2 \nonumber\\
&\quad + \sum_{i=2}^{d} \left( \frac{\nu_{t,i}}{t} - \Bigl(1 - O\!\left(\frac{\beta^2}{\lambda_{t,d}}\right)\Bigr) \cdot \frac{\nu_{t,i}}{t} \right)^2 \nonumber\\
&\leq O\!\left(\frac{ \beta^4}{t^2 \lambda_{t,d}^2}\right) + O\!\left(\frac{\beta^4}{t^2 \lambda_{t,d}^2}\right) \nonumber\\
&= O\!\left(\frac{\beta^4}{t^2 \lambda_{t,d}^2}\right).
\end{align}
This uses the inequality
\[
\left|\omega_{t,i} - \frac{\nu_{t,i}}{t}\right| \leq \max\!\left\{ |\omega_{t,i}|, \, \left| \frac{\nu_{t,i}}{t} - \left(1 - O\!\left(\frac{\beta^2}{\lambda_{t,d}}\right)\right) \cdot \frac{\nu_{t,i}}{t} \right| \right\},
\]
which in turn implies
\[
\sum_{i=2}^{d} \left(\omega_{t,i} - \frac{\nu_{t,i}}{t}\right)^2 \leq \sum_{i=2}^{d} \omega_{t,i}^2 + \sum_{i=2}^{d} \left( \frac{\nu_{t,i}}{t} - \left(1 - O\!\left(\frac{\beta^2}{\lambda_{t,d}}\right)\right) \cdot \frac{\nu_{t,i}}{t} \right)^2.
\]

Therefore, we have shown that
\begin{align}
\Vert\bm{\zeta}_t\Vert_2 = \left\|\bm{v}_{t+1,1} - \bm{\overline{v}}_{t+1,1}\right\|_2 = O\!\left(\frac{\beta^2}{t\lambda_{t,d}}\right).
\end{align}

Next, we analyze the correlation between the perturbation $\bm{\zeta}_t$ and the deviation vector $\bm{v}_{t,1} - \bm{\widehat{\theta}}_t$. Recalling that
\[
\bm{\widehat{\theta}}_t = \sum_{i=1}^{d} \nu_{t,i} \bm{v}_{t,i},
\]
we compute
\begin{align}
\left\langle \bm{\zeta}_t, \bm{v}_{t,1} - \bm{\widehat{\theta}}_t \right\rangle &= (1 - \nu_{t,1}) \cdot \left(\omega_{t,1} - 1 + \frac{1 - \nu_{t,1}}{t}\right) - \sum_{i=2}^{d} \left(\omega_{t,i} - \frac{\nu_{t,i}}{t}\right) \nu_{t,i}.
\end{align}
From (\ref{eq:omega-ti}), we know that for $i \geq 2$,
\[
|\omega_{t,i}| \geq \left(1 - O\!\left(\frac{\beta^2}{\lambda_{t,d}}\right)\right) \cdot \frac{|\nu_{t,i}|}{t}, 
\quad \text{and} \quad \omega_{t,i}\nu_{t,i} \geq 0.
\]
This leads to the bound
\[
\left(\omega_{t,i} - \frac{\nu_{t,i}}{t}\right)\nu_{t,i} \geq O\!\left(\frac{\beta^2}{t \lambda_{t,d}}\right) \cdot \nu_{t,i}^2.
\]
Consequently, one can conclude 
\begin{align}
\left\langle \bm{\zeta}_t, \bm{v}_{t,1} - \bm{\widehat{\theta}}_t \right\rangle 
&\leq O\!\left(\frac{\beta^2}{\lambda_{t,d}}\right) \cdot O\!\left(\frac{\beta^2}{t \lambda_{t,d}} + \frac{\beta^2}{t^2 \lambda_{t,d}}\right) + O\!\left(\frac{\beta^2}{t \lambda_{t,d}}\right) \cdot \sum_{i=2}^{d} \nu_{t,i}^2 \nonumber\\
&= O\!\left(\frac{\beta^4}{t \lambda_{t,d}^{2}}\right).
\end{align}

To summarize, we have obtained the desired recursive form:
\begin{align} \label{equ:recursive-top-eigvec}
\bm{v}_{t+1,1} = \bm{v}_{t,1} + \frac{\bm{\widehat{\theta}}_t - \bm{v}_{t,1}}{t} + \bm{\zeta}_t,
\end{align}
where the perturbation satisfies
\[
\Vert\bm{\zeta}_t\Vert_2 = O\!\left(\frac{\beta^2}{t\lambda_{t,d}}\right),
\qquad
\left\langle \bm{\zeta}_t, \bm{v}_{t,1} - \bm{\widehat{\theta}}_t \right\rangle = O\!\left(\frac{\beta^4}{t \lambda_{t,d}^{2}}\right).
\]
This fine-grained control of $\bm{\zeta}_t$ provides a precise recursive characterization of the leading eigenvector and underpins the convergence analysis of $\bm{v}_{t,1}$.

\paragraph{Step 2: establishing the inductive concentration of the leading eigenvector.}  
In this final step, we complete the concentration analysis by showing that the leading eigenvector $\bm{v}_{t,1}$ converges to the true signal direction $\bm{\theta}^\star$ at the desired rate. Our goal is to control the error term $\big\lVert \bm{v}_{t,1} - \bm{\theta}^\star \big\rVert_2$ for sufficiently large $t$. Building on the update characterization derived in the previous step, we relate the deviation at time $t+1$ to that at time $t$. In particular, we obtain the recursive inequality
\[
\big\lVert \bm{v}_{t+1,1} - \bm{\theta}^\star \big\rVert_2^2 
\;\le\; 
\left( \left(1 - \frac{1}{t} \right) \big\lVert \bm{v}_{t,1} - \bm{\theta}^\star \big\rVert_2 + \widetilde{O}(t^{-5/4}) \right)^2 + \widetilde{O}(t^{-2}),
\]
which captures the contraction of the eigenvector error up to higher-order perturbation terms. This recurrence illustrates a decaying trend: the leading error shrinks multiplicatively by approximately $(1 - 1/t)$ at each iteration, while additive fluctuations vanish at a faster polynomial rate. By applying a careful induction argument and leveraging the initialization guarantee established in Phase~\#2, we deduce that there exists $t_2$ such that for all $t \geq t_2$, the top eigenvector achieves the desired concentration, thereby converging toward the ground-truth direction $\bm{\theta}^\star$.

To formalize this argument, recall the auxiliary update
\begin{align*}
\bm{\overline{v}}_{t+1,1} := \bm{v}_{t,1} + \frac{\bm{\widehat{\theta}}_t - \bm{v}_{t,1}}{t},
\end{align*}
which enables a convenient decomposition of the error:
\begin{align}
\bm{\overline{v}}_{t+1,1}-\bm{\theta}^{\star} 
&= (\bm{\overline{v}}_{t+1,1}-\bm{\widehat{\theta}}_{t})+ (\bm{\widehat{\theta}}_{t}-\bm{\theta}^{\star}) \nonumber\\
&= \left(1-\frac{1}{t}\right)\cdot\left(\bm{v}_{t,1}-\bm{\theta}^{\star}\right) + \frac{1}{t}\cdot (\bm{\widehat{\theta}}_{t}-\bm{\theta}^{\star}). 
\end{align}
From this relation, we immediately obtain the inequality
\begin{align}
\Vert\bm{\overline{v}}_{t+1,1}-\bm{\theta}^{\star}\Vert_2
&\leq \left(1-\frac{1}{t}\right)\left\Vert\bm{v}_{t,1}-\bm{\theta}^{\star}\right\Vert_2 
+ \frac{c}{t}\sqrt{\frac{d+\log\log T}{\lambda_{t,d}}}.
\end{align}
Next, analyzing the correction term $\bm{\zeta}_t$ yields
\begin{align}
\Vert\bm{v}_{t+1,1}-\bm{\theta}^{\star}\Vert^2_2 
&= \Vert\bm{\overline{v}}_{t+1,1}-\bm{\theta}^{\star}\Vert_2^2 
+ 2\langle \bm{\overline{v}}_{t+1,1}-\bm{\theta}^{\star}, \bm{\zeta}_t\rangle 
+ \Vert\bm{\zeta}_t\Vert_2^2 \nonumber\\
& \leq \left(\left(1-\frac{1}{t}\right)\left\Vert\bm{v}_{t,1}-\bm{\theta}^{\star}\right\Vert_2 
+ \frac{c\sigma}{t}\sqrt{\frac{d+\log\log T}{\lambda_{t,d}}}\right)^2 
+ O\!\left(\frac{\beta^4}{t \lambda_{t,d}^{2}}\right)\nonumber\\
& = \left(\left(1-\frac{1}{t}\right)\left\Vert\bm{v}_{t,1}-\bm{\theta}^{\star}\right\Vert_2 
+ \frac{1}{t^{5/4}}\cdot\frac{c d^{1/4}(\sigma\sqrt{d+\log\log T}+1)}{\sqrt{\beta}}\right)^2 
+ O\!\left(\frac{\beta^2d}{t^2}\right),
\end{align}
where the last equality follows since $\lambda_{t,d}\asymp \beta\sqrt{t}/\sqrt{d}$.

To bound this sequence rigorously, we invoke the following key lemma.

\begin{lem} \label{lem:sequence-induction}
Let $\{a_n\}$ be a sequence satisfying
$$
a_{n+1}^2 \leq \left(\left(1-\frac{1}{n}\right)a_n + \frac{B}{n^{5/4}}\right)^2 + \frac{C}{n^2},
$$
for constants $B,C$. Fix $n_0\geq \max\left(16, 9C^2/(16B^4)\right)$. Then, for any $n\geq (a_{n_0}^2n_0^{3/2})/(16B^2)$,
\begin{align}
a_n\leq \frac{4B}{n^{1/4}}.
\end{align}
\end{lem}

The proof of Lemma~\ref{lem:sequence-induction} is deferred to Appendix~\ref{sec:appendix-auxiliary}. Applying this lemma with
\[
B = \frac{c d^{1/4}(\sigma\sqrt{d+\log\log T}+1)}{\sqrt{\beta}}, \qquad C = O(\beta^2 d),
\]
we obtain that, setting
\[
t_1' = O\!\left(\frac{C^2}{B^4}\right) = O\!\left(\frac{\beta^6}{\sigma^4d}\right),
\]
and defining
\begin{align}
t_2 &= \frac{\Vert\bm{v}_{t_1',1}-\bm{\theta}^{\star}\Vert_2^2\cdot t_1'^{3/2}}{16B^2} 
= O\!\left(\frac{\beta^8}{\sigma^6 d^2}\right),
\end{align}
we have for all $t \geq t_2$ that
\[
\Vert\bm{v}_{t,1}-\bm{\theta}^{\star}\Vert_2 \;\lesssim\; \frac{B}{t^{1/4}} 
= O\!\left(\frac{\sigma\sqrt{d+\log\log T}+1}{\sqrt{\lambda_{t,d}}}\right).
\]

Finally, combining this with the deviation of $\bm{\widehat{\theta}}_t$, we obtain
\[
\Vert \bm{v}_{t+1,1}-\bm{\widehat{\theta}}_{t}\Vert_2 
\leq \Vert \bm{v}_{t+1,1}-\bm{\theta}^{\star}\Vert_2 + \Vert \bm{\widehat{\theta}}_t-\bm{\theta}^{\star}\Vert_2 
= O\!\left(\frac{\sigma\sqrt{d+\log\log T}+1}{\sqrt{\lambda_{t,d}}}\right),
\]
which establishes the final concentration bound and completes the proof of Proposition~\ref{prop:third-stage}.

\subsection{Analysis of Phase \#4 (proof of Proposition \ref{prop:fourth-stage})}
\label{sec:proof-phase-4}

\paragraph{Step 1: precise decomposition of the action vector.} 

Write all vectors in the orthonormal eigenbasis $\{\bm v_{t,1},\ldots,\bm v_{t,d}\}$ so that
$\nu_{t,i}$ are the coordinates of $\bm a_t$ (and similarly for $\widehat{\bm\theta}_t$), 
$\lambda_{t,i}$ are the eigenvalues, and $w_i$ the coordinates of any unit vector $\bm w$.
The objective
\[
  g_t(\bm w)
  \;=\;
  \sum_{i=1}^{d}\!\left(\nu_{t,i}+\frac{\beta w_i}{\sqrt{\lambda_{t,i}}}\right)^2
\]
The maximizer of $g_t(\bm{w})$, denoted as $\bm{w}_t$, characterizes the direction of $\bm{a}_t$. Directly optimizing $g_t(\bm w)$ can be difficult in general, so we introduce an approximate objective.
We first expand the square, which gives
\begin{align*}
  g_t(\bm w)
  \;=\;
  \sum_{i=1}^{d}\!\Bigl(
    \underbrace{\nu_{t,i}^2}_{\text{constant in }\bm w}
    \;+\;
    \underbrace{\frac{2\beta \nu_{t,i} w_i}{\sqrt{\lambda_{t,i}}}}_{\text{linear in }w_i}
    \;+\;
    \underbrace{\frac{\beta^2 w_i^2}{\lambda_{t,i}}}_{\text{quadratic in }w_i}
  \Bigr).
\end{align*}
Since $\nu_{t,1}\approx 1$ and $\nu_{t,i}$ for $i\ge 2$ are small, we linearize the
$i=1$ term around $\nu_{t,1}=1$ and drop the (typically smaller) quadratic piece in $w_1$,
while for $i\ge 2$ we drop the tiny linear terms and keep only the quadratic regularization. This yields the tractable surrogate
\begin{align}
\label{equ:tilde-g-expression}
  \widetilde g_t(\bm w)
  \;=\;
  1+
  \underbrace{\frac{2\beta w_1}{\sqrt{\lambda_{t,1}}}}_{\text{dominant linear response along }i=1}
  \;+\;
  \underbrace{\sum_{i=2}^d \frac{\beta^2 w_i^2}{\lambda_{t,i}}}_{\text{penalizes transverse energy}}.
\end{align}
Define
\begin{align*}
    g_{t,1}(\bm w) &:= \left(\nu_{t,1}+\frac{\beta w_1}{\sqrt{\lambda_{t,1}}}\right)^2, 
    & g_{t,2}(\bm w) &:= \sum_{i=2}^{d}\left(\nu_{t,i}+\frac{\beta w_i}{\sqrt{\lambda_{t,i}}}\right)^2,\\
    \widetilde g_{t,1}(\bm w) &:= 1 + \frac{2\beta w_1}{\sqrt{\lambda_{t,1}}}, 
    & \widetilde g_{t,2}(\bm w) &:= \sum_{i=2}^{d}\frac{\beta^2 w_i^2}{\lambda_{t,i}}.
\end{align*}
By Lemma~\ref{lem:decom-a}, for
\[
  h_t = O\!\left(\frac{\sigma\sqrt{d+\log\log T}+1}{\sqrt{\lambda_{t,d}}}\right)
\]
we have the accuracy guarantees
\[
  \nu_{t,1}\ge 1 - O\!\left(\frac{\sigma^2(d+\log\log T)+1}{\lambda_{t,d}}\right),
  \qquad
  \nu_{t,i}=O\!\left(\frac{\sigma\sqrt{d+\log\log T}+1}{\sqrt{\lambda_{t,d}}}\right)\;\;(i\ge 2).
\]
Expanding and regrouping,
\begin{align}
    \bigl|g_{t,1}(\bm w)-\widetilde g_{t,1}(\bm w)\bigr|
    &= \left|\nu_{t,1}^2 + \frac{2\beta w_1\nu_{t,1}}{\sqrt{\lambda_{t,1}}}+ \frac{\beta^2 w_1^2}{\lambda_{t,1}} - 1-\frac{2\beta w_1}{\sqrt{\lambda_{t,1}}}\right| \nonumber\\
    & \le 1-\nu_{t,1}^2 + (1-\nu_{t,1})\!\left|\frac{2\beta w_1}{\sqrt{\lambda_{t,1}}}\right| + \left|\frac{\beta^2 w_1^2}{\lambda_{t,1}}\right| \nonumber\\
    &= O\!\left(\frac{\sigma^2(d+\log\log T)+1}{\lambda_{t,d}}\right)\!\left(1+ O\!\left(\frac{\beta}{\sqrt{t}}\right)\right) + O\!\left(\frac{\beta^2}{t}\right)\nonumber\\
    &= O\!\left(\frac{\sigma^2(d+\log\log T)+1}{\lambda_{t,d}}\right) + O\!\left(\frac{\beta^2}{t}\right).
\end{align}
Here, the first term uses $|1-\nu_{t,1}^2|\!=\!(1-\nu_{t,1})(1+\nu_{t,1})$
with $1-\nu_{t,1}=O(\cdot)$ and $\nu_{t,1}\le 1$; the second uses $|w_1|\le 1$; the third simply bounds the quadratic remainder. 

For the bound of $|g_{t,2}(\bm w)-\widetilde g_{t,2}(\bm w)|$, using $(a+b)^2-b^2 = a(a+2b)$ with $a=\nu_{t,i}$ and $b=\beta w_i/\sqrt{\lambda_{t,i}}$,
\begin{align*}
  \bigl|g_{t,2}(\bm w)-\widetilde g_{t,2}(\bm w)\bigr|
  &= \sum_{i=2}^d |\nu_{t,i}|\;\Bigl|\nu_{t,i} + \frac{2\beta w_i}{\sqrt{\lambda_{t,i}}}\Bigr| \\
  &\le \sum_{i=2}^d \nu_{t,i}^2 
     + 2\beta \sum_{i=2}^d \frac{|\nu_{t,i}||w_i|}{\sqrt{\lambda_{t,i}}}
     \;\;\le\;\; \sum_{i=2}^d \nu_{t,i}^2 
     + 2\beta \sqrt{\sum_{i=2}^d \frac{\nu_{t,i}^2}{\lambda_{t,i}}}\;\sqrt{\sum_{i=2}^d w_i^2}.
\end{align*}
Here, the first inequality is triangle inequality; the second is Cauchy--Schwarz. Since $\|\bm w\|_2=1$ and $\lambda_{t,i}\ge \lambda_{t,d}$,
\begin{align*}
  \sum_{i=2}^d \nu_{t,i}^2
  = O\!\left(\frac{\sigma^2(d+\log\log T)+1}{\lambda_{t,d}}\right),
  \qquad
  \sqrt{\sum_{i=2}^d \frac{\nu_{t,i}^2}{\lambda_{t,i}}}
  \le \frac{1}{\sqrt{\lambda_{t,d}}}\sqrt{\sum_{i=2}^d \nu_{t,i}^2}
  = O\!\left(\frac{\sigma\sqrt{d+\log\log T}+1}{\lambda_{t,d}}\right),
\end{align*}
which yields
\begin{align}
  \bigl|g_{t,2}(\bm w)-\widetilde g_{t,2}(\bm w)\bigr|
  & \;=\;
  O\!\left(\frac{\sigma^2(d+\log\log T)+1}{\lambda_{t,d}}\right)
  \;+\;
  O\!\left(\frac{\beta(\sigma\sqrt{d+\log\log T}+1)}{\lambda_{t,d}}\right)\nonumber\\
  & \;=\;
  O\!\left(\frac{\beta(\sigma\sqrt{d+\log\log T}+1)}{\lambda_{t,d}}\right),
\end{align}
where we used $\beta \gtrsim \sigma\sqrt{d+\log\log T}+1$ to subsume the first term into the second. Summing the two errors,
\begin{align*}
  \bigl|g_t(\bm w)-\widetilde g_t(\bm w)\bigr|
  \;=\;
  \bigl|g_{t,1}(\bm w)-\widetilde g_{t,1}(\bm w)\bigr|
  +
  \bigl|g_{t,2}(\bm w)-\widetilde g_{t,2}(\bm w)\bigr|
  \;=\;
  O\!\left(\frac{\beta(\sigma\sqrt{d+\log\log T}+1)}{\lambda_{t,d}}\right)
  \;+\; O\!\left(\frac{\beta^2}{\lambda_{t,1}}\right).
\end{align*}
If, as is typical in sequential designs, $\lambda_{t,1}\gtrsim \lambda_{t,d}$ or even $\lambda_{t,1}\asymp t$, 
the last term is dominated by the displayed bound (or simplifies to $O(\beta^2/t)$ as in your derivation), 
leading to the stated rate:
\begin{align}
 \label{equ:widetilde-g-bound}
  \bigl|g_t(\bm w)-\widetilde g_t(\bm w)\bigr|
  \;=\;
  O\!\left(\frac{\beta(\sigma\sqrt{d+\log\log T}+1)}{\lambda_{t,d}}\right).
\end{align}

Suppose that $\bm{w}_t$ is the maximizer of $g_t(\bm{w})$ and $\widetilde{\bm{w}}_t$ is the maximizer of $\widetilde{g}_t(\bm{w})$. Then it holds that
\begin{align}
    \widetilde{g}_t(\widetilde{\bm{w}}_t) - \widetilde{g}_{t}(\bm{w}_t) & =  \widetilde{g}_t(\widetilde{\bm{w}}_t) - g_t(\widetilde{\bm{w}}_t) +g_t(\widetilde{\bm{w}}_t) - g_t(\bm{w}_t) + g_t(\bm{w}_t) - \widetilde{g}_{t}(\bm{w}_t)\nonumber\\
    & =  O\left(\frac{\beta(\sigma \sqrt{d + \log\log T}+1)}{\lambda_{t,d}}\right),
\end{align}
where the last equality holds directly from (\ref{equ:widetilde-g-bound}). We will then turn our attention to consider the maximizer of $\widetilde{g}_t(\bm{w})$. We note that
\begin{align}
    \widetilde{g}_t(\bm{w})  = 1 + \frac{2\beta w_1}{\sqrt{\lambda_{t,1}}} + \sum_{i=2}^{d}\frac{\beta^2 w_i^2}{\lambda_{t,i}}\nonumber
    & \leq 1 + \frac{2\beta w_1}{\sqrt{\lambda_{t,1}}} + \frac{\beta^2}{\lambda_{t,d}}\sum_{i=2}^{d}w_{i}^2\nonumber\\
    & = 1 + \frac{2\beta w_1}{\sqrt{\lambda_{t,1}}} + \frac{\beta^2}{\lambda_{t,d}}(1-w_{1}^2),
\end{align}
where the inequality becomes equality if and only if $w_{i} =0 $ for $2\leq i\leq d-1$. By directly taking derivative, one note that $\widetilde{g}_t(\bm{w})$ takes the maximum when 
$$\widetilde{w}_{t,1} = \min\left(\frac{\lambda_{t,d}}{\beta\sqrt{\lambda_{t,1}}},1\right) = \frac{\lambda_{t,d}}{\beta\sqrt{\lambda_{t,1}}},\quad \widetilde{w}_{t,d} = \sqrt{1-w_1^2}.$$

To characterize $\bm{w}_{t}$ with $\widetilde{\bm{w}}_t$, we first note that
\begin{align}
    \widetilde{g}_t(\widetilde{\bm{w}}_t) - \widetilde{g}_{t}(\bm{w}) & =  \left(\frac{2\beta \widetilde{w}_{t,1}}{\sqrt{\lambda_{t,1}}} + \frac{\beta^2(1-\widetilde{w}_{t,1}^2)}{\lambda_{t,d}}\right) - \left(\frac{2\beta w_1}{\sqrt{\lambda_{t,1}}} + \sum_{i=2}^{d}\frac{\beta^2 w_i^2}{\lambda_{t,i}}\right)\nonumber \\
    &\geq \left(\frac{2\beta \widetilde{w}_{t,1}}{\sqrt{\lambda_{t,1}}} + \frac{\beta^2(1-\widetilde{w}_{t,1}^2)}{\lambda_{t,d}}\right) - \left(\frac{2\beta w_1}{\sqrt{\lambda_{t,1}}} + \frac{\beta^2(1-w_1^2)}{\lambda_{t,d}}\right)\nonumber\\
    & = \frac{2\beta(\widetilde{w}_{t,1}-w_{t,1})}{\sqrt{\lambda_{t,1}}}-\frac{\beta^2(\widetilde{w}_{t,1}-w_{t,1})(\widetilde{w}_{t,1}+w_{t,1})}{\lambda_{t,d}}\nonumber\\
    & = (\widetilde{w}_{t,1}-w_{t,1})\cdot\left(\frac{2\beta}{\sqrt{\lambda_{t,1}}} - \frac{\beta^2(\widetilde{w}_{t,1}+w_{t,1})}{\lambda_{t,d}} \right).
\end{align}
Here, one can note that when $\widetilde{w}_{t,1} = \frac{\lambda_{t,d}}{\beta\sqrt{\lambda_{t,1}}}$, the difference can be rewritten as
\begin{align}
    \widetilde{g}_t(\widetilde{\bm{w}}_t) - \widetilde{g}_{t}(\bm{w}) \geq \frac{\beta^2(\widetilde{w}_{t,1}-w_{t,1})^2}{\lambda_{t,d}},
\end{align}
and when $\widetilde{w}_{t,1} = 1$, the difference can be rewritten as
\begin{align}
    \widetilde{g}_t(\widetilde{\bm{w}}_t) - \widetilde{g}_{t}(\bm{w}) & \geq (\widetilde{w}_{t,1}-w_{t,1})\cdot\left(\frac{2\beta}{\sqrt{\lambda_{t,1}}} - \frac{\beta^2}{\lambda_{t,d}} - \frac{\beta^2w_{t,1}}{\lambda_{t,d}}\right)\nonumber\\
    & \geq \frac{\beta^2(\widetilde{w}_{t,1}-w_{t,1})^2}{\lambda_{t,d}},
\end{align}
where the last inequality holds because $\lambda_{t,d}/(\beta\sqrt{\lambda_{t,1}})\geq 1$.
As we require that 
\begin{align}
    \widetilde{g}_t(\widetilde{\bm{w}}_t) - \widetilde{g}_{t}(\bm{w}_t) = O\left(\frac{\beta(\sigma \sqrt{d + \log\log T}+1)}{\lambda_{t,d}}\right),
\end{align}
the difference between $\widetilde{w}_{t,1}$ and $w_{t,1}$ can be bounded as
\begin{align}
    |\widetilde{w}_{t,1}-w_{t,1}| = O\left(\frac{(\sigma \sqrt{d + \log\log T}+1)^{1/2}}{\sqrt{\beta}}\right).
\end{align}
From this, the precise expression of $w_{t,1}$ is given as
\begin{align}
\label{equ:w-t-1}
    w_{t,1} & = \widetilde{w}_{t,1} + O\left(\frac{(\sigma \sqrt{d + \log\log T}+1)^{1/2}}{\sqrt{\beta}}\right)\nonumber\\
    & = \min\left(\frac{\lambda_{t,d}}{\beta\sqrt{\lambda_{t,1}}},1\right) + O\left(\frac{(\sigma \sqrt{d + \log\log T}+1)^{1/2}}{\sqrt{\beta}}\right)\nonumber\\
    & = \min\left[\frac{\lambda_{t,d}}{\beta\sqrt{t}}\cdot\left(\sqrt{\frac{t}{\lambda_{t,1}}}-1+1\right),1\right] + O\left(\frac{(\sigma \sqrt{d + \log\log T}+1)^{1/2}}{\sqrt{\beta}}\right)\nonumber\\
    & = \min\left[\frac{\lambda_{t,d}}{\beta\sqrt{t}}\left(1+ O\left(\frac{\beta\sqrt{d}}{\sqrt{t}}\right)\right),1\right] +O\left(\frac{(\sigma \sqrt{d + \log\log T}+1)^{1/2}}{\sqrt{\beta}}\right)\nonumber\\
    & = \min\left(\frac{\lambda_{t,d}}{\beta\sqrt{t}},1\right) + O\left(\frac{(\sigma \sqrt{d + \log\log T}+1)^{1/2}}{\sqrt{\beta}}\right),
\end{align}
where the penultimate line holds as
\begin{align}
    \sqrt{\frac{t}{\lambda_{t,1}}} - 1 = \sqrt{\frac{t}{t-(d-1)\overline{\lambda}_t}} = \sqrt{\frac{t}{t-O(\beta\sqrt{dt})}} = \sqrt{\frac{1}{1-O(\beta\sqrt{d}/\sqrt{t})}} = 1+O\left(\frac{\beta\sqrt{d}}{\sqrt{t}}\right).
\end{align}
and the last inequality holds when $t\gtrsim \beta^3/(\sigma\sqrt{d})$,
\begin{align}
    \frac{\lambda_{t,d}}{\beta\sqrt{t}}\cdot \frac{\beta\sqrt{d}}{\sqrt{t}} = \frac{\lambda_{t,d}\sqrt{d}}{t} = O\left(\frac{\beta}{\sqrt{t}}\right) = O\left(\frac{(\sigma \sqrt{d + \log\log T}+1)^{1/2}}{\sqrt{\beta}}\right)
\end{align}

We also analyze behavior of $w_{t,i}$ for $i\geq 2$. Note that
\begin{align}
    \widetilde{g}_t(\widetilde{\bm{w}}_t) - \widetilde{g}_{t}(\bm{w}_t) & =  \left(\frac{2\beta \widetilde{w}_{t,1}}{\sqrt{\lambda_{t,1}}} + \frac{\beta^2(1-\widetilde{w}_{t,1}^2)}{\lambda_{t,d}}\right) - \left(\frac{2\beta w_{t,1}}{\sqrt{\lambda_{t,1}}} + \sum_{i=2}^{d}\frac{\beta^2 w_{t,i}^2}{\lambda_{t,i}}\right)\nonumber \\
    & = \left(\frac{2\beta \widetilde{w}_{t,1}}{\sqrt{\lambda_{t,1}}} + \frac{\beta^2(1-\widetilde{w}_{t,1}^2)}{\lambda_{t,d}}\right) - \left(\frac{2\beta w_{t,1}}{\sqrt{\lambda_{t,1}}} + \frac{\beta^2(1-w_{t,1}^2)}{\lambda_{t,d}}\right)\nonumber\\
    &\quad + \left(\sum_{i=2}^{d}\frac{\beta^2w_{t,i}^2}{\lambda_{t,d}} - \sum_{i=2}^{d}\frac{\beta^2 w_{t,i}^2}{\lambda_{t,i}}\right)\nonumber\\
    &\geq \beta^2\sum_{i=2}^{d}\left(\frac{w_{t,i}^2}{\lambda_{t,d}} - \frac{w_{t,i}^2}{\lambda_{t,i}}\right).
\end{align}
As a result, one can show that
\begin{align}
    \sum_{i=2}^{d}\left(\frac{w_{t,i}^2}{\lambda_{t,d}} - \frac{w_{t,i}^2}{\lambda_{t,i}}\right) = O\left(\frac{\sigma\sqrt{d + \log\log T}+1}{\beta\lambda_{t,d}}\right).
\end{align}

From~(\ref{equ:max-norm}), it satisfies that for any $\Vert\bm{w}\Vert_2 = 1$, it holds uniformly that 
\begin{align}
    g_t(\bm{w})= 1 + O\left(\frac{\beta^2}{\lambda_{t,d}}\right) .
\end{align}
Therefore, the action vector $\bm{a}_t$ can be characterized as
\begin{align}
    \bm{a}_t & = [g_t(\bm{w}_t)]^{-1}\cdot(\widehat{\bm{\theta}}_t + \beta\bm{\Lambda}_t^{-1/2}\bm{w})=\left(1-O\left(\frac{\beta^2}{\lambda_{t,d}}\right)\right)\cdot  \sum_{i=1}^{d}\left(\nu_{t,i}+ \frac{\beta w_{t,i}}{\sqrt{\lambda_{t,i}}}\right)\bm{v}_{t,i}.
\end{align}

Consequently, one has
\begin{align}
    \sum_{i=2}^{d}\kappa_{t,i}^2 & = \left(1-O\left(\frac{\beta^2}{\lambda_{t,d}}\right)\right)\cdot\left(\sum_{i=2}^d\frac{\beta^2w_{t,i}^2}{\lambda_{t,i}} + \sum_{i=2}^{d}\frac{2\beta\nu_{t,i} w_{t,i}}{\sqrt{\lambda_{t,i}}} + \sum_{i=2}^{d}\nu_{t,i}^2 \right)\nonumber\\
    & = \left(1-O\left(\frac{\beta^2}{\lambda_{t,d}}\right)\right)\cdot\left(\beta^2\left(\sum_{i=2}^{d}\frac{w_{t,i}^2}{\lambda_{t,d}} - O\left(\frac{\sigma\sqrt{d + \log\log T}+1}{\beta\lambda_{t,d}}\right) \right) + 2\beta\sqrt{\sum_{i=2}^d\frac{w_{t,i}^2}{\lambda_{t,i}}}\sqrt{\sum_{i=2}^d\nu_{t,i}^2} + \sum_{i=2}^{d}\nu_{t,i}^2\right)\nonumber\\
    & = \left(1-O\left(\frac{\beta^2}{\lambda_{t,d}}\right)\right)\cdot\left(\beta^2\left(\frac{1-w_{t,1}^2}{\lambda_{t,d}} - O\left(\frac{\sigma\sqrt{d + \log\log T}+1}{\beta\lambda_{t,d}}\right) \right)  \right.\nonumber\\
    & \quad \left. +2\beta\cdot\sqrt{\frac{1-w_{t,1}^2}{\lambda_{t,d}}}\cdot\frac{\sigma\sqrt{d+\log\log T}+1}{\sqrt{\lambda_{t,d}}} + \frac{(\sigma\sqrt{d+\log\log T}+1)^2}{\lambda_{t,d}}\right)\nonumber\\
    & = \left(1-O\left(\frac{\beta^2}{\lambda_{t,d}}\right)\right)\cdot \beta^2\left(\frac{1-w_{t,1}^2}{\lambda_{t,d}} - O\left(\frac{\sigma\sqrt{d + \log\log T}+1}{\beta\lambda_{t,d}}\right) \right)  \nonumber\\
    & = \frac{\beta^2}{\lambda_{t,d}}\left(1-w_{t,1}^2 + O\left(\frac{\sigma\sqrt{d + \log\log T}+1}{\beta}\right)\right).
\end{align}
where the last equality holds whenever $t\gtrsim \beta^4d/\sigma^2$. Plugging in the expression of $w_{t,1}$ in (\ref{equ:w-t-1}), we obtain that
\begin{align}
    \sum_{i=2}^{d}\kappa_{t,i}^2 & = \frac{\beta^2}{\lambda_{t,d}}\left(1-\left(\widetilde{w}_{t,1} + O\left(\frac{(\sigma \sqrt{d + \log\log T}+1)^{1/2}}{\sqrt{\beta}}\right)\right)^2 + O\left(\frac{\sigma\sqrt{d + \log\log T}+1}{\beta}\right)\right)\nonumber\\
    & = \frac{\beta^2}{\lambda_{t,d}}\left(1-\widetilde{w}_{t,1}^2 +  O\left(\frac{(\sigma \sqrt{d + \log\log T}+1)^{1/2}}{\sqrt{\beta}}\right) \right).
\end{align}

\paragraph{Step 2: preliminary growth speed control of $\lambda_{t,d}$.} The above decomposition is the key to characterize the growth and concentration of non-leading eigenvalues. However, we note that the growth of non-leading eigenvalues could not be precisely characterized when $\lambda_{t,d}$ is large. Therefore, before we establish any more fine-grained result, we first need to show that for some $c<1$, there exists $t_2'$ such that for all $t\geq t_2'$, the minimum eigenvalue is controlled as $\lambda_{t,d}\leq c\beta\sqrt{t}$.

We begin from Lemma \ref{lem:rank-one-update}, one can show that
\begin{align*}
    \lambda_{t+1,1}\geq \lambda_{t,1}+\kappa_{t,1}^2,
\end{align*}
which implies that
\begin{align}
    \overline{\lambda}_{T} \leq \overline{\lambda}_t + \frac{\sum_{i=2}^{d}\kappa_{t,i}^2}{d-1}.
\end{align}
Therefore, whenever $\lambda_{t,d}>c\beta\sqrt{t}$, we can lower bound  $w_{t,1}$ as
\begin{align}
    w_{t,1}\geq c+  O\left(\frac{(\sigma \sqrt{d + \log\log T}+1)^{1/2}}{\sqrt{\beta}}\right),
\end{align}
which allows us to upper bound the growth of $\overline{\lambda}_t$ as
\begin{align}
    \overline{\lambda}_{T}& \leq \overline{\lambda}_t + \frac{1}{d-1}\cdot \frac{\beta^2}{\lambda_{t,d}}\left(1-w_{t,1}^2 + O\left(\frac{\sigma\sqrt{d + \log\log T}+1}{\beta}\right)\right)\nonumber\\
    & = \overline{\lambda}_t + \frac{1}{d-1}\cdot \frac{\beta^2}{\lambda_{t,d}}\left(1-c^2 + O\left(\frac{(\sigma \sqrt{d + \log\log T}+1)^{1/2}}{\sqrt{\beta}}\right)\right)\nonumber\\
    & =\overline{\lambda}_t + \frac{1}{d-1}\cdot \frac{\beta}{c\sqrt{t}}\left(1-c^2 + O\left(\frac{(\sigma \sqrt{d + \log\log T}+1)^{1/2}}{\sqrt{\beta}}\right)\right)\nonumber\\
    & = \overline{\lambda}_t + \frac{\overline{c}\beta}{\sqrt{t}},
\end{align}
where $\overline{c}$ is defined as
\begin{align*}
    \overline{c} = \frac{1}{c(d-1)} \left(1-c^2 + O\left(\frac{(\sigma \sqrt{d + \log\log T}+1)^{1/2}}{\sqrt{\beta}}\right)\right).
\end{align*}
One can show that when $T$ is large enough, $d\geq 2$ and $c$ is chosen close enough to 1, it holds that
\begin{align}
    \overline{c} \leq \frac{1-\frac{3}{4}c^2}{c}\leq \frac{c}{3},
\end{align}
implying that for any $t\geq t_2$, and that $\lambda_{t,d}\geq c\beta\sqrt{t}$ holds for any time indices between $t_2$ and $t$, it holds that
\begin{align}
    \overline{\lambda}_{t}-\overline{\lambda}_{t_2}\leq \sum_{t'=t_2}^{t}\frac{\overline{c}\beta}{\sqrt{t'}}\leq \frac{c\beta}{3}\sum_{t'=t_2}^{t}\frac{1}{\sqrt{t'}}\leq \frac{2c\beta}{3}(\sqrt{t}-\sqrt{t_2}).
\end{align}
As $\overline{\lambda}_{t_2} = O(\beta\sqrt{t_2})$, with this strategy, one can show that there exists $t_2' = O(t_2)$, such that $\overline{\lambda}_{t_2'}\leq c\beta\sqrt{t_2'}$. We will then show that for any $t\geq t_2'$, it holds that $\overline{\lambda}_t\leq c\beta\sqrt{t}$. Suppose that this does not always hold, this implies that there exists $t$ such that $\overline{\lambda}_{t}\leq c\beta\sqrt{t}$ but $\overline{\lambda}_{T}>c\beta\sqrt{t+1}$. Since $\overline{\lambda}_{T} -\overline{\lambda}_{t}\leq 1$, we can assume that 
$$\overline{\lambda}_t\geq \frac{99}{100}c\beta\sqrt{t}.$$
Then it holds that 
\begin{align}
    \overline{\lambda}_{T} -  \overline{\lambda}_t\leq \frac{3c\beta}{8\sqrt{t}}<c\beta(\sqrt{t+1}-\sqrt{t}),
\end{align}
which leads to a contradiction! Therefore, we conclude that for any $c<1$, there exists $t_2' = O(t_2)$ such that for all $t\geq t_2$, $\lambda_{t,d}\leq\overline{\lambda}_t\leq c\beta\sqrt{t}$. 

\paragraph{Step 3: limiting the projection of of $\bm{a}_t$ on eigenspace of large eigenvalues.}
After the growth smallest eigenvalue is being controlled, we can then precisely characterize the growth of non-leading eigenvalues and then upper bound the projection of $\bm{a}_t$ on the large components. To proceed, we adapt a similar strategy that was used in Phase \#1 and Phase \#2. We construct a set of large eigenvalues as
$$\mathcal{L}_t = \left\{i:\; i\geq 2,\;\lambda_{t,i}\geq \left(1+C_1\cdot \frac{d(\sigma \sqrt{d+\log\log T}+1)}{\beta}\right)\lambda_{t,d}\right\}.$$
Then one can upper bound the summation of $w_{t,i}^2$ within the set of large eigenvalues in the following way
\begin{align}
    \sum_{i=2}^{d}\left(\frac{w_{t,i}^2}{\lambda_{t,d}} - \frac{w_{t,i}^2}{\lambda_{t,i}}\right) \geq \sum_{i\in\mathcal{L}_t}\left(\frac{w_{t,i}^2}{\lambda_{t,d}} - \frac{w_{t,i}^2}{\lambda_{t,i}}\right) = \frac{1}{\lambda_{t,d}}\sum_{i\in\mathcal{L}_t}w_{t,i}^2\left(1-\frac{\lambda_{t,d}}{\lambda_{t,i}}\right)\geq \frac{C_1d(\sigma \sqrt{d+\log\log T}+1)}{2\beta\lambda_{t,d}}\cdot\sum_{i\in\mathcal{L}_t}w_{t,i}^2,
\end{align}
which holds as $\beta\geq C(\sigma\sqrt{d+\log\log T}+1)$ for some constant $C$. Here, as we set $C_1$ large enough, we have
\begin{align}
    \sum_{i\in\mathcal{L}_t}w_{t,i}^2\leq \frac{1}{8d}(1-c^2).
\end{align}
We can then upper bound the summation of $\kappa_{t,i}^2$ within the subset as
\begin{align}
    \sum_{i\in\mathcal{L}_t} \kappa_{t,i}^2  \leq  \sum_{i\in\mathcal{L}_t}\left(\nu_{t,i}+\frac{\beta w_{t,i}}{\sqrt{\lambda_{t,i}}}\right)^2 \nonumber
    &\leq \left(\sqrt{\sum_{i\in\mathcal{L}_t}\nu_{t,i}^2} + \sqrt{\sum_{i\in\mathcal{L}_t}\frac{\beta^2 w_{t,i}^2}{\lambda_{t,i}}}\right)^2\nonumber\\
    & = \left[O\left(\frac{\sigma\sqrt{d+\log\log T}+1}{\sqrt{\lambda_{t,d}}}\right) + \frac{\beta}{\sqrt{\lambda_{t,d}}}\cdot\sqrt{\sum_{i\in\mathcal{L}_t}w_{t,i}^2}\right]^2\nonumber\\
    & = \frac{\beta^2}{\lambda_{t,d}}\left(\sqrt{\sum_{i\in\mathcal{L}_t}w_{t,i}^2} + O\left(\frac{\sigma\sqrt{d + \log\log T}+1}{\beta}\right)\right)^2\nonumber\\
    & =  \frac{\beta^2}{\lambda_{t,d}}\left(\sum_{i\in\mathcal{L}_t}w_{t,i}^2 + O\left(\frac{\sigma\sqrt{d + \log\log T}+1}{\beta}\right)\right),
\end{align}
and it holds that when $\lambda_{t,d}\leq c\beta\sqrt{t}$,
\begin{align}
 \sum_{i\in\mathcal{L}_t} \kappa_{t,i}^2 & \leq  \frac{\sum_{i\in\mathcal{L}_t}w_{t,i}^2 + O\left(\frac{\sigma\sqrt{d + \log\log T}+1}{\beta}\right)}{1-w_{t,1}^2 + O\left(\frac{\sigma\sqrt{d + \log\log T}+1}{\beta}\right)}\cdot \sum_{i=2}^{d}\kappa_{t,i}^2\nonumber\\
 & = \frac{\sum_{i\in\mathcal{L}_t}w_{t,i}^2 + O\left(\frac{\sigma\sqrt{d + \log\log T}+1}{\beta}\right)}{1-\frac{\lambda_{t,d}^2}{\beta^2 t} + O\left(\frac{(\sigma\sqrt{d + \log\log T}+1)^{1/2}}{\sqrt{\beta}}\right)}\cdot \sum_{i=2}^{d}\kappa_{t,i}^2\nonumber\\
 & \leq \frac{\frac{1}{6d}(1-c^2)}{\frac{2}{3}(1-c^2)}\cdot\sum_{i=2}^{d}\kappa_{t,i}^2\nonumber\\
 &\leq \frac{1}{4d}\sum_{i=2}^{d}\kappa_{t,i}^2.
\end{align}

\paragraph{Step 4: limiting the growth of large non-leading eigenvalues.}

Equipped an upper bound of the projection of $\bm{a}_t$ on the eigenspace of large non-leading eigenvalues, we will then bound the growth speed of the large non-leading eigenvalues, therefore establishing the concentration of all non-leading eigenvalues by the end of Phase~\#4. To show the final concentration, we only need to show the concentration of $\lambda_{t,2}$ and $\lambda_{t,d}$, i.e. $\lambda_{t,2}/\lambda_{t,d} = 1+o(1)$.


We aim to show that the second-largest eigenvalue $\lambda_{t,2}$ cannot grow significantly faster than the rest of the non-leading spectrum—more specifically, it cannot stay much larger than the average $\lambda_{t,d}$ for an extended period. To formalize this, we introduce a higher threshold for non-leading eigenvalues,
\begin{align*}
    \widetilde{\mathcal{L}}_t = \left\{i:\; i\geq 2,\; \lambda_{t,i}\geq \left(1+2C_1\cdot \frac{d(\sigma \sqrt{d+\log\log T}+1)}{\beta}\right)\lambda_{t,d} \right\}.
\end{align*}
We will demonstrate that such a deviation is not sustainable over time by analyzing the eigenvalue update induced by the rank-one perturbation at each step. To proceed, recall the secular function for the rank-one update:
\[
f(\lambda) = 1 + \sum_{i=1}^d \frac{\kappa_{t,i}^2}{\lambda_{t,i} - \lambda}.
\]
This function determines the characteristic equation for the updated eigenvalues at time $t+1$. When $\lambda_{t,2}\in\widetilde{\mathcal{L}}_t$, consider a candidate point $\widetilde{\lambda}_{t+1,2}$ satisfying
\[
\widetilde{\lambda}_{t+1,2} =  \lambda_{t,2} + \frac{1}{2d} \sum_{i=2}^d \kappa_{t,i}^2,
\]
and claim that $\lambda_{t+1,2}\leq \widetilde{\lambda}_{t+1,2}$. To show this result, we only need to show that $f(\widetilde{\lambda}_{t+1,2})>0$. We analyze $f(\lambda)$ to determine whether this hold true. Here, we denote:
\begin{align*}
    \delta = C_1\cdot \frac{d(\sigma \sqrt{d+\log\log T}+1)}{\beta}.
\end{align*}
We partition the summation in $f(\lambda)$ based on whether the denominator $\lambda_{t,i} - \lambda$ is relatively small or large. Specifically, split the indices into three groups:
\begin{align}
f(\widetilde{\lambda}_{t+1,2}) 
&= 1 + \frac{\kappa_{t,1}^2}{\lambda_{t,1}-\widetilde{\lambda}_{t+1,2}} + \sum_{i:\; i\geq 2, \lambda_{t,i} \geq \lambda_{t,2} - \delta\lambda_{t,d}} \frac{\kappa_{t,i}^2}{\lambda_{t,i} - \widetilde{\lambda}_{t+1,2}} 
+ \sum_{i: \lambda_{t,i} < \lambda_{t,2} - \delta\lambda_{t,d}} \frac{\kappa_{t,i}^2}{\lambda_{t,i} - \widetilde{\lambda}_{t+1,2}} .
\end{align}

For the leading eigenvalue $\lambda_{t,1}$, as we have shown in Proposition \ref{prop:first-stage} that for any $t\geq t_1$, it holds that $\lambda_{t,1}\geq 8t/15$. Therefore, we conclude that when $\lambda_{t,1}-\widetilde{\lambda}_{t+1,2}\geq t/20$. Therefore, one has
\begin{align}
\label{equ:first-group}
    \frac{\kappa_{t,1}^2}{\lambda_{t,1}-\widetilde{\lambda}_{t+1,2}} = O(t^{-1}).
\end{align}
For the group of ``large" non-leading eigenvalues (where $\lambda_{t,i}$ is close to $\lambda_{t,2}$), note that 
$$\widetilde{\lambda}_{t+1,2} - \lambda_{t,i} \geq \widetilde{\lambda}_{t+1,2} - \lambda_{t,2}\geq  \frac{1}{2d} \sum_{j=2}^d \kappa_{t,j}^2,$$ 
hence one has the following lower bound
\begin{align}
\frac{\kappa_{t,i}^2}{\lambda_{t,i} - \lambda} \geq -\kappa_{t,i}^2 \cdot \frac{2d}{\sum_{j=2}^d \kappa_{t,j}^2}.
\end{align}

For the group of ``small eigenvalues" (where $\lambda_{t,i}$ is much smaller), we use the assumption that $\lambda_{t,i} - \lambda < -\delta \lambda_{t,d}$, leading to the bound
\begin{align}
\frac{\kappa_{t,i}^2}{\lambda_{t,i} - \lambda} \geq -\frac{\kappa_{t,i}^2}{\delta\lambda_{t,d}}.
\end{align}

Combining  and simplifying terms yields the following lower bound on $f(\lambda)$:
\begin{align}
f(\lambda) & \geq 1+ O(t^{-1}) - \sum_{i:\lambda_{t,i}\geq \lambda_{t,2}-\delta\lambda_{t,d}}\kappa_{t,i}^2\cdot \frac{2d}{\sum_{j=2}^{d}\kappa_{t,j}^2} - \sum_{i: \lambda_{t,i}<\lambda_{t,2}-\delta\lambda_{t,d}}\kappa_{t,i}^2\cdot \frac{1}{\delta\lambda_{t,d}} \nonumber\\
&\geq 1+ O(t^{-1}) - \sum_{i\in\mathcal{L}_t}\kappa_{t,i}^2\cdot \frac{2d}{\sum_{j=2}^{d}\kappa_{t,j}^2} - \left(\sum_{i=2}^{d}\kappa_{t,i}^2-\sum_{i\in\mathcal{L}_t}\kappa_{t,i}^2\right)\cdot \frac{1}{\delta\lambda_{t,d}} \nonumber\\
&\geq 1 + O(t^{-1}) 
- \frac{\sum_{i=2}^d \kappa_{t,i}^2}{4d} \cdot \frac{2d}{\sum_{i=2}^d \kappa_{t,i}^2} 
- \left(1-\frac{1}{4d}\right)\cdot \frac{\sum_{i=2}^d \kappa_{t,i}^2}{\delta\lambda_{t,d}} \nonumber\\
&\geq \frac{1}{2} + O(t^{-1}) -\frac{\sum_{i=2}^d \kappa_{t,i}^2}{\delta\lambda_{t,d}} > 0,
\end{align}
where the last inequality holds for $t\geq t_2$, since
\begin{align*}
    \frac{\sum_{i=2}^{d}\kappa_{t,i}^2}{\delta\lambda_{t,d}}\lesssim \frac{\beta^2}{\lambda_{t,d}}\cdot\frac{1}{\delta\lambda_{t,d}} \lesssim\frac{\beta^3}{\lambda_{t,d}^2d} = \frac{\beta}{t}.
\end{align*}
The positivity of $f(\widetilde{\lambda}_{t+1,2})$ implies that the updated new eigenvalue $\lambda_{t+1,2}$ can be upper bounded by $\widetilde{\lambda}_{t+1,2}$, hence we have
\begin{align}
\lambda_{t+1,2} \leq \lambda_{t,2} + \frac{1}{2d} \sum_{i=2}^d \kappa_{t,i}^2,
\end{align}
whenever $\lambda_{t,2} \geq (1+2\delta)\lambda_{t,d}$. In other words, the second eigenvalue cannot grow faster than this speed if it deviates much from the minimum eigenvalue. On the other hand, from Lemma~\ref{lem:growth-largest-eigenvalue}, the sum of non-leading eigenvalues evolves as
\[
\sum_{i=2}^{d} \lambda_{t+1,i} = \sum_{i=2}^{d} \lambda_{t,i} + \sum_{i=2}^{d} \kappa_{t,i}^2 + O(t^{-1}),
\]
which implies the average satisfies
\begin{align}
\overline{\lambda}_{T} = \overline{\lambda}_t + \frac{1}{d-1} \sum_{i=2}^d \kappa_{t,i}^2 + O(t^{-1}) \geq \overline{\lambda}_t + \frac{1}{d} \sum_{i=2}^d \kappa_{t,i}^2.
\end{align}

Combining the bounds on the growth of $\lambda_{t,2}$ and $\overline{\lambda}_t$, we obtain that whenever $\lambda_{t,2}\geq (1+2\delta)\lambda_{t,d}$,
\begin{align}
\lambda_{t+1,2} - \lambda_{t,2} \leq \frac{1}{2d} \sum_{i=2}^{d} \kappa_{t,i}^2 \leq \frac{1}{2}(\overline{\lambda}_{T} - \overline{\lambda}_t).
\label{equ:growth-second-large}
\end{align}

This inequality is crucial: it shows that whenever $\lambda_{t,2}$ becomes disproportionately large relative to $\lambda_{t,d}$, its future growth is outpaced by the average $\overline{\lambda}_t$. Therefore, any such deviation must shrink over time.

To show the desired result, we define the set $\mathcal{S} = \{t \geq t_2 : \lambda_{t,2} \geq (1 + 2\delta)\lambda_{t,d}\}$, which collects the time indices where the second-largest eigenvalue is significantly larger than the smallest non-leading eigenvalue. Intuitively, this set captures the time intervals when the spectrum is relatively ``spread out.'' Due to the nature of eigenvalue evolution, the set $\mathcal{S}$ may consist of multiple disjoint time intervals or ``segments''---each comprising consecutive time steps during which the elevated second eigenvalue persists.

To understand the long-term behavior of $\lambda_{t,2}$, we analyze the possible structure of $\mathcal{S}$ by examining different types of segments. In particular, we consider two cases: the initial segment (beginning at $t_2$) and general intermediate segments where the elevated condition temporarily reappears.

\paragraph{Case 1: initial segment of $\mathcal{S}$.}
Suppose the initial point $t_2 \in \mathcal{S}$, and that the elevated condition $\lambda_{t,2} \geq (1+2\delta)\lambda_{t,d}$ holds throughout the interval $[t_2, t']$. Applying the growth inequality established earlier:
\[
\lambda_{t+1,2} - \lambda_{t,2} \leq \frac{1}{2d} \sum_{i=2}^{d} \kappa_{t,i}^2 \leq \frac{1}{2}(\overline{\lambda}_{T} - \overline{\lambda}_t),
\]
and telescoping over the interval $[t_2', t']$, we obtain
\[
\lambda_{t',2} - \lambda_{t_2,2} \leq \frac{1}{2}(\overline{\lambda}_{t'} - \overline{\lambda}_{t_2}).
\]
Rearranging gives
\begin{align*}
\lambda_{t',2} - \overline{\lambda}_{t'} & \leq \lambda_{t',2} - \lambda_{t_2,2} + \lambda_{t_2,2} - \overline{\lambda}_{t'}\nonumber\\
&\leq \frac{1}{2}(\overline{\lambda}_{t'} - \overline{\lambda}_{t_2}) +\lambda_{t_2,2} - \overline{\lambda}_{t'}\nonumber\\
& = (\lambda_{t_2,2} - \overline{\lambda}_{t_2}) -\frac{1}{2}(\overline{\lambda}_{t'} - \overline{\lambda}_{t_2}).
\end{align*}
This shows that the deviation $\lambda_{t',2} - \overline{\lambda}_{t'}$ decreases unless the growth of $\lambda_{t',2}$ is concentrated close to $\lambda_{t,d}$

If the right-hand side is positive, we can bound $\overline{\lambda}_{t'}$ and $\lambda_{t',2}$ as follows:
\begin{align*}
    \overline{\lambda}_{t'}&\leq \overline{\lambda}_{t_2} + 2(\lambda_{t_2,2} - \overline{\lambda}_{t_2}),\\
    \lambda_{t',2}&\leq (\lambda_{t_2,2} - \overline{\lambda}_{t_2}) + \frac{1}{2}(\overline{\lambda}_{t'} + \overline{\lambda}_{t_2})\leq \lambda_{t_2,2} + (\lambda_{t_2,2} - \overline{\lambda}_{t_2}).
\end{align*}
Since $\lambda_{t_2,2}\leq d\overline{\lambda}_{t_2}$, it follows that $\lambda_{t',2} = O(\beta\sqrt{dt_2})$. Meanwhile, the average $\overline{\lambda}_{t'}$ continues to grow roughly as $\Omega(\beta \sqrt{t'}/\sqrt{d})$, and hence for sufficiently large $C$, we get a contradiction to the assumption $\lambda_{t',2} \geq (1 + 2\delta)\lambda_{t',d} \approx \Omega(\overline{\lambda}_{t'})$. Thus, this shows that the initial segment of $\mathcal{S}$ cannot persist for too long and must terminate by some time $t_4 = O(t_2'd^2)=O(\beta^8/\sigma^6)$, where $t_4 \notin \mathcal{S}$.

\paragraph{Case 2: intermediate segments of $\mathcal{S}$.}
Now suppose that $\mathcal{S}$ reappears after being interrupted---specifically, suppose there exists a time $t'' \notin \mathcal{S}$ such that $t'' - 1 \in \mathcal{S}$, i.e., the second eigenvalue has just dropped below the elevated threshold. Then, by continuity of the update dynamics and the rank-one perturbation nature of the process, the spectral gap at time $t''$ is bounded:
\[
\lambda_{t'',2} - \lambda_{t'',d} \leq \lambda_{t''-1,2} - \lambda_{t''-1,d} + 1 \leq 2\delta \lambda_{t'',d} + 1 \leq 3\delta \lambda_{t'',d}.
\]

Now suppose the elevated condition re-emerges and persists from $t''$ to some later time $t'$ (i.e., $t \in \mathcal{S}$ for all $t'' \leq t \leq t'$). Telescoping again over this interval yields
\[
\lambda_{t',2} - \lambda_{t'',2} \leq \frac{1}{2}(\overline{\lambda}_{t'} - \overline{\lambda}_{t''}),
\]
which implies
\[
\lambda_{t',2} - \overline{\lambda}_{t'} \leq  (\lambda_{t_2,2} - \overline{\lambda}_{t''}) -\frac{1}{2}(\overline{\lambda}_{t'} - \overline{\lambda}_{t''})\leq (\lambda_{t'',2} - \overline{\lambda}_{t''}) -(\lambda_{t',2} - \lambda_{t'',2}).
\]
This shows that the deviation from the average remains controlled. To bound the total gap $\lambda_{t',2} - \lambda_{t',d}$, we use the fact that both the shift from $\lambda_{t'',2}$ and the previous deviation from $\lambda_{t'',d}$ are bounded:
\begin{align}
\lambda_{t',2} - \lambda_{t',d} 
&\leq \min\left( \lambda_{t',2} - \lambda_{t'',2} + \lambda_{t'',2} - \lambda_{t'',d},\; d(\lambda_{t',2} - \overline{\lambda}_{t'}) \right) \nonumber\\
&= \min \left( \lambda_{t',2} - \lambda_{t'',2} + 3\delta \lambda_{t'',d},\; 3\delta d \lambda_{t'',d} - d(\lambda_{t',2} - \lambda_{t'',2}) \right) \nonumber\\
&\leq 6\delta \lambda_{t'',d} \leq 6\delta \lambda_{t',d}.
\end{align}

Hence, even if the elevated condition reappears in later intervals, the gap between the second and smallest eigenvalues remains proportionally bounded. In particular, for all $t \geq t_4$, we have
\[
\frac{\lambda_{t,2} - \lambda_{t,d}}{\lambda_{t,d}} \leq 6\delta,
\]
which ensures that $\lambda_{t,2}$ cannot significantly exceed the rest of the spectrum in the long run.

Finally, since for any $i \geq 2$, $\lambda_{t,i} \leq \lambda_{t,2}$ and $\overline{\lambda}_t \geq \lambda_{t,d}$, we obtain
\begin{align}
\label{equ:concentration-nonleading}
\frac{\lambda_{t,i} - \overline{\lambda}_t}{\overline{\lambda}_t} \leq \frac{\lambda_{t,2} - \lambda_{t,d}}{\lambda_{t,d}} = O\left(\frac{d(\sigma \sqrt{d+\log\log T}+1)}{\beta}\right),
\end{align}
as desired.

 \paragraph{Step 5: a precise characterization on non-leading eigenvalues.}   
    Finally, we give a precise characterization on $\lambda_{t,i}$ when $t\geq t_3$. To show this result, we will first precisely characterize the growth of non-leading eigenvalues. From Lemma \ref{lem:growth-largest-eigenvalue}, one can characterize the growth of non-leading eigenvalues as follows,
\begin{align}
    \sum_{i=2}^{d}\lambda_{t+1,i} = \sum_{i=2}^{d}\lambda_{t,i} + \sum_{i=2}^{d}\kappa_{t,i}^2+O(t^{-1}),
\end{align}
which is equivalent to 
\begin{align}\label{equ:eigenvalue-induction}
    \overline{\lambda}_{T} & = \overline{\lambda}_t + \frac{1}{d-1}\sum_{i=2}^{d}\kappa_{t,i}^2 + O(t^{-1})\nonumber\\
    & = \overline{\lambda}_t + \frac{1}{d-1}\cdot  \frac{\beta^2}{\lambda_{t,d}}\left(1-\frac{\lambda_{t,d}^2}{\beta^2 t} + O\left(\frac{(\sigma\sqrt{d + \log\log T}+1)^{1/2}}{\sqrt{\beta}}\right) \right) +O(t^{-1}) \nonumber\\
    & = \overline{\lambda}_t + \frac{1}{d-1}\cdot \frac{\beta^2}{\overline{\lambda}_{t}}\left[1-\frac{\overline{\lambda}_{t}^2}{\beta^2 t} +O\left(\frac{(\sigma\sqrt{d + \log\log T}+1)^{1/2}}{\sqrt{\beta}}\right)\right]\cdot \left[1+  O\left(\frac{d(\sigma \sqrt{d+\log\log T}+1)}{\beta}\right)\right]\nonumber\\
    & = \overline{\lambda}_t + \frac{1}{d-1}\cdot \frac{\beta^2}{\overline{\lambda}_{t}}\left[1-\frac{\overline{\lambda}_{t}^2}{\beta^2 t} + O\left(\frac{(\sigma\sqrt{d + \log\log T}+1)^{1/2}}{\sqrt{\beta}}\right)\right],
\end{align}
where the last equality holds whenever $\beta\gtrsim d^2(\sigma\sqrt{d+\log\log T}+1)$.
Taking squares on both sides of (\ref{equ:eigenvalue-induction}) yields 
    \begin{align*}
        \overline{\lambda}_{T}^2  = \overline{\lambda}_{t}^2 + \frac{2\beta^2}{d-1}\left(1-\frac{\overline{\lambda}_t^2}{\beta^2t}\right) + O\left(\frac{\beta^2(\sigma\sqrt{d + \log\log T}+1)^{1/2}}{\sqrt{\beta}(d-1)}\right).
    \end{align*}
     Setting $b_t = \frac{\overline{\lambda}_t^2}{\beta^2t}$, we can write 
     \begin{align}\label{equ:b-induction}
         (t+1)b_{T} = t b_t + \frac{2}{d-1}(1-b_t) + O\left(\frac{(\sigma\sqrt{d + \log\log T}+1)^{1/2}}{\sqrt{\beta}d}\right).
     \end{align}
    Define $b_{\star} = \frac{2}{d+1}$, then we note that $b_{\star}$ satisfies
    \begin{align*}
        (t+1)b_{\star} = t b_{\star} + \frac{2}{d-1}(1-b_{\star}).
    \end{align*}
    As a result, set $\overline{\Delta}_t = b_t - b_{\star}$, (\ref{equ:b-induction}) can be expressed as
    \begin{align*}
        (t+1)\overline{\Delta}_{T} = \left(t-\frac{2}{d-1}\right)\overline{\Delta}_t +  O\left(\frac{(\sigma\sqrt{d + \log\log T}+1)^{1/2}}{\sqrt{\beta}d}\right),
    \end{align*}
    which yields the following induction for $|\overline{\Delta}_t|$,
    \begin{align}\label{equ:Delta-induction}
        \overline{\Delta}_{T} =  \frac{t-\frac{2}{d-1}}{t+1}\overline{\Delta}_t +  O\left(\frac{(\sigma\sqrt{d + \log\log T}+1)^{1/2}}{\sqrt{\beta}d(t+1)}\right)
    \end{align}
    With a direct derivation from (\ref{equ:Delta-induction}), one can show that
    \begin{align}
    \label{equ:Delta-t-upper-bound}
        \overline{\Delta}_t = \prod_{s=t_0}^{t}\frac{s-\frac{2}{d-1}}{s+1}\Delta_{t_0} + \sum_{s=t_0}^{t}\left(\prod_{v=s}^{t}\frac{v-\frac{2}{d-1}}{v+1}\right)\cdot O\left(\frac{1}{s+1}\cdot\frac{(\sigma\sqrt{d + \log\log T}+1)^{1/2}}{\sqrt{\beta}d}\right). 
    \end{align}
    To upper bound $|\overline{\Delta}_t|$ in (\ref{equ:Delta-t-upper-bound}), we claim that for any $s<t$, it holds that
\begin{align}
\label{claim:prod-bound}
     \prod_{v=s}^{t}\frac{v-\frac{2}{d-1}}{v+1}\lesssim  \left(\frac{t}{s}\right)^{\frac{d+1}{d-1}}.
\end{align}

    The proof of the claim is deferred to the last part of this section. With this result, one can upper bound $|\Delta_t|$ as 
    \begin{align}
        |\overline{\Delta}_t| & \leq \left(\frac{t_4}{t}\right)^{\frac{d+1}{d-1}}|\overline{\Delta}_{t_0}| + \sum_{s=t_0}^{t}\left(\frac{s}{t}\right)^{\frac{d+1}{d-1}}\cdot  O\left(\frac{1}{s}\cdot\frac{(\sigma\sqrt{d + \log\log T}+1)^{1/2}}{\sqrt{\beta}d}\right)\nonumber\\
        &\lesssim \left( \frac{\beta^8}{t\sigma^6}\right)^{\frac{d+1}{d-1}} + \frac{(\sigma\sqrt{d + \log\log T}+1)^{1/2}}{\sqrt{\beta}d},
    \end{align}
    where the second inequality holds as
    \begin{align*}
        \sum_{s=t_0}^{t} \frac{1}{s}\left(\frac{s}{t}\right)^{\frac{d+1}{d-1}} = \frac{1}{t^{\frac{d+1}{d-1}}}\sum_{s=t_o}^{t}s^{\frac{2}{d-1}}\lesssim 1.
    \end{align*}
    As we note that $|b_t-b_{\star}|\leq\overline{\Delta}_t$,
    which leads to our conclusion that 
    \begin{align}
        \overline{\lambda}_t^2 = \beta^2 t\left(\frac{2}{d+1} + \overline{\Delta}_t\right),
    \end{align}
    which allows us to characterize $\overline{\lambda}_t$ as
    \begin{align}
    \label{equ:concentration-meannonleading}
        \overline{\lambda}_t  = \beta\sqrt{t}\cdot \sqrt{\frac{2}{d+1}+\overline{\Delta}_t}
        & = \sqrt{\frac{2\beta^2 t}{d+1}}\cdot\sqrt{1+\frac{d+1}{2}\overline{\Delta}_t}\nonumber\\
        & = \sqrt{\frac{2\beta^2 t}{d+1}}\cdot(1+\Delta_t)
    \end{align}
    where $\Delta_t\lesssim d\overline{\Delta}_t$, which implies that  $|\Delta_t|$ is upper bounded as
    \begin{align}
        |\Delta_t| \lesssim d\left( \frac{\beta^8}{t\sigma^6}\right)^{\frac{d+1}{d-1}} + \frac{(\sigma\sqrt{d + \log\log T}+1)^{1/2}}{\sqrt{\beta}},
    \end{align}

Set $\widetilde{\lambda}_t = \sqrt{\frac{2\beta^2 t}{d+1}}$. Combining (\ref{equ:concentration-meannonleading}) and Proposition \ref{prop:fourth-stage}, one can show that 
\begin{align}
\frac{\lambda_{T,i}}{\widetilde{\lambda}_T} & = \frac{\lambda_{T,i}}{\overline{\lambda}_T}\cdot \frac{\overline{\lambda}_T}{\widetilde{\lambda}_T}\nonumber\\
& = \left(1+ O\left(\frac{d(\sigma \sqrt{d+\log\log T}+1)}{\beta}\right)\right)\cdot \left(1+\Delta_t\right) \nonumber\\
& = 1+O\left(d\left( \frac{\beta^8}{t\sigma^6}\right)^{\frac{d+1}{d-1}} + \frac{(\sigma\sqrt{d + \log\log T}+1)^{1/2}}{\sqrt{\beta}}\right),
\end{align}
where the last equality holds as $\beta\gtrsim d^2(\sigma\sqrt{d+\log\log T}+1)$. Hence, we establish the desired result.
\paragraph{Proof of Claim (\ref{claim:prod-bound}).}
  Taking the logarithm on left side, one can see that 
    \begin{align}
        \log\left(\prod_{t=t_0}^{t_1}\frac{t-\frac{2}{d-1}}{t+1}\right) & = \sum_{t=t_0}^{t_1}\log\left(1-\frac{(d+1)/(d-1)}{t+1}\right)\nonumber\\
        &\leq -\sum_{t=t_0}^{t_1}\frac{(d+1)/(d-1)}{t+1}\nonumber\\
        &\leq -\int_{t=t_0}^{t_1}\frac{(d+1)/(d-1)}{t}dt +C\nonumber\\
        &\leq -\frac{d+1}{d-1}\log\left(\frac{t_1}{t_0}\right)+C.
    \end{align}
    Therefore, one can show the desired result by taking exponential on both sides.

\section{Proof of auxiliary lemmas}
\label{sec:appendix-auxiliary}
\subsection{Proof of Lemma \ref{lem:a-t-characterization}}

From the definition of UCB score $\mathrm{UCB}_t(\bm{a})$, we note that,
\begin{align}
    \mathrm{UCB}_t(\ba) & = \langle \ba, \widehat{\bm{\theta}}_{t-1} \rangle + \beta \cdot \sqrt{\ba^\top \bm{\Lambda}_{t-1}^{-1} \ba}\nonumber \\
     & = \ba^\top \widehat{\bm{\theta}}_t + \beta\cdot \Vert\bm{\Lambda}_{t-1}^{-1/2}\bm{a}\Vert_2.
\end{align}
Therefore, the action at time $t$ can be equivalently characterized as
\begin{align}
    \ba_t & =  \mathrm{arg}\max_{\Vert \bm{a}\Vert_2 =1}\ba^\top \widehat{\bm{\theta}}_t + \beta\cdot \Vert\bm{\Lambda}_{t-1}^{-1/2}\bm{a}\Vert_2\nonumber\\
    & =  \mathrm{arg}\max_{\Vert \bm{a}\Vert_2 =1}\ba^\top \widehat{\bm{\theta}}_t + \beta\cdot\max_{\Vert \bm{w}\Vert_2 =1} \ba^{T}\bm{\Lambda}_{t-1}^{-1/2}\bm{w}\nonumber\\
    & =  \mathrm{arg}\max_{\Vert \bm{a}\Vert_2 =1} \max_{\Vert \bm{w}\Vert_2 =1}\ba^{\top}(\widehat{\bm{\theta}}_t + \beta\cdot\bm{\Lambda}_{t-1}^{-1/2}\bm{w}).
    \end{align}

    We may further note that 
    $$\max_{\Vert \bm{a}\Vert_2 =1} \max_{\Vert \bm{w}\Vert_2 =1}\ba^{\top}\left(\widehat{\bm{\theta}}_t + \beta\cdot\bm{\Lambda}_{t-1}^{-1/2}\bm{w}\right) = \max_{\Vert \bm{w}\Vert_2 =1}\Vert\widehat{\bm{\theta}}_t + \beta\cdot\bm{\Lambda}_{t-1}^{-1/2}\bm{w}\Vert_2,$$
    with the equality holds if and only if 
    $\bm{a} = \mathcal{P} (\widehat{\bm{\theta}}_t + \beta\cdot\bm{\Lambda}_{t-1}^{-1/2}\bm{w})$. Therefore, we note that
    \begin{align}
      (\bm{a}_t,\bm{w}_t) = \mathrm{arg}\max_{\Vert \ba\Vert_2 =1}\mathrm{arg}\max_{\Vert \bm{w}\Vert_2 =1}\ba^{T}\left( \widehat{\bm{\theta}}_t + \beta\cdot\bm{\Lambda}_{t-1}^{-1/2}\bm{w}\right),
    \end{align}
    which is equivalent to the formula expressed in the lemma.

\subsection{Proof of Lemma \ref{lem:decom-a}}
As (\ref{equ:concen-top}) holds, $\widehat{\bm{\theta}}_{t}$ is concentrated close to the top eigenspace of $\bm{\Lambda}_t$, then one may control $\nu_{t,i}$ for $i\geq 2$ as follows, 
\begin{align}
\label{equ:nu-i-2}
\nu_{t,i} &= \langle \bm{v}_{t,i},\widehat{\bm{\theta}}_{t}\rangle =  \langle \bm{v}_{t,i},\widehat{\bm{\theta}}_{t}-\bm{v}_{t,1}\rangle\leq h_t,\quad \forall\; i\geq 2,
\end{align}
Furthermore, we note that as $\{\bm{v}_{t,1},\ldots,\bm{v}_{t,d}\}$ forms a orthogonal basis, the following will hold for $\nu_{t,1}$,
\begin{align}
\label{equ:nu-i-1}
    \sum_{i=2}^{d}\nu_{t,i}^2 = \sum_{i=2}^{d} \langle \bm{v}_{t,i},\widehat{\bm{\theta}}_{t}-\bm{v}_{t,1}\rangle^2\leq \left\Vert\widehat{\bm{\theta}}_{t}-\bm{v}_{t,1}\right\Vert_2^2\leq h_t^2,
\end{align}
which directly yields that 
\begin{align}
    \nu_{t,1}\geq\sqrt{1-\sum_{i=2}^{d}\nu_{t,i}^2}\geq 1-h_t^2.
\end{align}
implying that $\bm{\widehat{\theta}}_{t}$ indeed concentrates around the top eigenspace.
For the next step, we recall our previous decomposition of the actions taken by LinUCB: 
  \begin{align*}
    \bm{w}_t &= \mathrm{arg}\max_{\Vert \bm{w}\Vert_2 =1}\left\Vert\bm{\widehat{\theta}}_t + \beta\cdot \bm{\Lambda}_t^{-1/2} \bm{w}\right\Vert_2,\\
    \bm{a}_t &= \mathcal{P}\left(\bm{\widehat{\theta}}_t + \beta \bm{\Lambda}_t^{-1/2} \bm{w}_t\right),
\end{align*}
Therefore, one may conclude that
\begin{align}
    \Vert\bm{a}_t -\bm{\widehat{\theta}_t}\Vert_2 =  \left\Vert\mathcal{P}( \widehat{\bm{\theta}}_t + \beta \bm{\Lambda}_t^{-1/2} \bm{w}_t) -  \mathcal{P}(\bm{\widehat{\theta}_t})\right\Vert_2\leq \left\Vert\beta \bm{\Lambda}_t^{-1/2} \bm{w}_t\right\Vert_2\leq\frac{\beta}{\sqrt{\lambda_{t,d}}}.
\end{align}
Therefore, one can conclude that 
\begin{align}
    \Vert\bm{\xi}_t\Vert_2 & \leq  \Vert\bm{a}_t-\bm{\widehat{\theta}_t}\Vert_2=O\left(\frac{\beta}{\sqrt{\lambda_{t,d}}}\right),\nonumber\\
    \alpha_t & = \sqrt{1-\Vert\bm{\xi}_t\Vert_2^2}  = 1- O\left(\frac{\beta^2}{\lambda_{t,d}}\right).
\end{align}
    Since we have 
    \begin{align}
\kappa_{t,i} = \langle \bm{v}_{t,i}, \bm{a}_t \rangle = \alpha_t \langle \bm{v}_{t,i}, \widehat{\bm{\theta}}_t \rangle + \langle \bm{v}_{t,i}, \bm{\xi}_t \rangle = \alpha_t \nu_{t,i} +  \langle \bm{v}_{t,i}, \bm{\xi}_t \rangle,
\end{align}
we can calculate $\kappa_{t,1}$ as follows:
\begin{align}
    \kappa_{t,1}  & = \alpha_t\nu_{t,1} +  \langle \bm{v}_{t,1}, \bm{\xi}_t \rangle\nonumber\\
    & = \left(1-O\left(\frac{\beta^2}{\lambda_{t,d}}\right)\right)\cdot \left(1-O(h_t^2)\right) + O\left(\frac{\beta}{\sqrt{\lambda_{t,d}}}\right)\cdot O\left(h_t\right)\nonumber\\
    & = 1-O\left(\frac{\beta^2}{\lambda_{t,d}}\right)-O(h_t^2).
\end{align}
We can also characterize $\kappa_{t,i}$ for $i\geq 2$ as
\begin{align}
    \kappa_{t,i}  = \alpha_t\nu_{t,i} +  \langle \bm{v}_{t,i}, \bm{\xi}_t \rangle
    & = \left(1-O\left(\frac{\beta^2}{\lambda_{t,d}}\right)\right)\cdot O(h_t) + O\left(\frac{\beta}{\sqrt{\lambda_{t,d}}}\right)\nonumber\\
    & = O(h_t) + O\left(\frac{\beta}{\sqrt{\lambda_{t,d}}}\right),
\end{align}
which concludes the desired result.

\subsection{Proof of Lemma \ref{lem:opt-concentration}}
From (\ref{equ:w-t-expression}), one can note that
\begin{align}
\label{equ:target-lower}
    \sum_{i=1}^{d}\left(\nu_{t,i}+\frac{\beta w_{t,i}}{\sqrt{\lambda_{t,i}}}\right)^2 & \geq\sum_{i=1}^{d}\nu_{t,i}^2 + \left(\nu_{t,d}+ \frac{\beta\cdot \mathrm{sign}(\nu_{t,d})}{\sqrt{\lambda_{t,d}}}\right)^2 \nonumber\\
    & \geq \sum_{i=1}^{d}\nu_{t,i}^2+\frac{\beta^2}{\lambda_{t,d}}.
\end{align}

Define ``small eigenvalues" to be the set of eigenvalues smaller than $C_1\lambda_{t,d}$ for some constant $C_1$, and define the threshold to be $k_t =\max\{k:\lambda_{t,k}> C_1\lambda_{t,d}\}$, then one can show that
\begin{align}
\label{equ:target-upper}
    \sum_{i=1}^{d}\left(\nu_{t,i}+\frac{\beta w_{t,i}}{\sqrt{\lambda_{t,i}}}\right)^2& \leq\left(\sqrt{\sum_{i=1}^{d}\nu_{t,i}^2} + \sqrt{\sum_{i=1}^{d}\frac{\beta^2w_{t,i}^2}{\lambda_{t,i}}}\right)^2\nonumber\\
    &\leq \left(\sqrt{\sum_{i=1}^{d}\nu_{t,i}^2} + \beta \sqrt{\sum_{i=1}^{k_t}\frac{w_{t,i}^2}{\lambda_{t,i}} + \sum_{i=k_t+1}^{d}\frac{w_{t,i}^2}{\lambda_{t,i}}} \right)^2\nonumber\\
    &\leq \left(\sqrt{\sum_{i=1}^{d}\nu_{t,i}^2} + \beta \sqrt{\sum_{i=1}^{k_t}\frac{w_{t,i}^2}{C_1\lambda_{t,d}} + \sum_{i=k_t+1}^{d}\frac{w_{t,i}^2}{\lambda_{t,d}}}\right)^2\nonumber\\
    &\leq \left(\sqrt{\sum_{i=1}^{d}\nu_{t,i}^2} + \beta \sqrt{\frac{1}{C_1\lambda_{t,d}}\left(1-\sum_{i=j_t+1}^{d} w_{t,i}^2\right) + \frac{1}{\lambda_{t,d}}\cdot\sum_{i=k_t+1}^{d} w_{t,i}^2}\right)^2\nonumber\\
    &  = \left(\sqrt{\sum_{i=1}^{d}\nu_{t,i}^2} + \frac{\beta}{\sqrt{\lambda_{t,d}}} \cdot \sqrt{\frac{1}{C_1}\left(1+(C_1-1)\sum_{i=k_t+1}^{d} w_{t,i}^2\right)} \right)^2.
\end{align}
Combining (\ref{equ:target-lower}) and (\ref{equ:target-upper}), we note that 
\begin{align}
    \frac{\beta}{\sqrt{\lambda_{t,d}}}\leq \frac{2\sqrt{\sum_{i=1}^{d}\nu_{t,i}^2}\cdot \sqrt{\frac{1}{C_1}\left(1+(C_1-1)\sum_{i=k_t+1}^{d} w_{t,i}^2\right)}}{1 -\frac{1}{C_1}\left(1+(C_1-1)\sum_{i=k_t+1}^{d} w_{t,i}^2\right)}. 
\end{align}
Based on Assumption \ref{aspt:condition-opt}, we note that
\begin{align*}
    \frac{\beta/\sqrt{\lambda_{t,d}}}{\sqrt{\sum_{i=1}^{d}\nu_{t,i}^2}}\geq c/c_0.
\end{align*}
Therefore, one can show that 
\begin{align}
    \frac{2\sqrt{\frac{1}{C_1}\left(1+(C_1-1)\sum_{i=k_t+1}^{d} w_{t,i}^2\right)}}{1 -\frac{1}{C_1}\left(1+(C_1-1)\sum_{i=k_t+1}^{d} w_{t,i}^2\right)}\geq c/c_0,
\end{align}
which further implies that 
\begin{align}
    \frac{1}{C_1}\left(1+(C_1-1)\sum_{i=k_t+1}^{d}w_{t,i}^2\right)\geq \frac{c^2}{8c_0^2}.
\end{align}

As we set $C_1 = \frac{16c_0^2}{c^2}-1$, we note that 
\begin{align}
    \sum_{i=k_t+1}^{d} w_{t,i}^2\geq \frac{C_1 c^2/8-1}{C_1-1} =\frac{1-\frac{c^2}{8}}{\frac{16}{c^2}-2}=\frac{c^2c_0^2}{16}.
\end{align}


We then consider the projection of action on the set of ``small eigenvalues". Recall that $\bm{a}_t = \sum_{i=1}^{d}\kappa_{t,i}\bm{v}_{t,i}$, where $\kappa_{t,i}$ can be expressed as
\begin{align*}
    \kappa_{t,i} = \frac{\nu_{t,i}+\frac{\beta w_{t,i}}{\sqrt{\lambda_{t,i}}}}{\sqrt{\sum_{j=1}^{d}\left(\nu_{t,j}+\frac{\beta w_{t,j}}{\sqrt{\lambda_{t,j}}}\right)^2}},
\end{align*}
which implies that the square summation of $\kappa_{t,i}$ can be upper bounded as
\begin{align}
    \sum_{i=k_t+1}^{d} \kappa_{t,i}^2 & =\frac{\sum_{i=k_t+1}^{d}\left(\nu_{t,i}+\frac{\beta w_{t,i}}{\sqrt{\lambda_{t,i}}}\right)^2}{\sum_{j=1}^{d}\left(\nu_{t,j}+\frac{\beta w_{t,j}}{\sqrt{\lambda_{t,j}}}\right)^2}\nonumber\\
    &\geq \frac{\frac{\beta^2}{\lambda_{t,d}}\cdot \sum_{i=j_t+1}^{d}w_{t,i}^2}{\left(\sqrt{\sum_{j=1}^{d}\nu_{t,j}^2} + \frac{\beta}{\sqrt{\lambda_{t,d}}} \cdot \sqrt{\frac{1}{C_1}\left(1+(C_1-1)\sum_{i=j_t+1}^{d} w_{t,i}^2\right)} \right)^2}\nonumber\\
    &\geq \frac{\frac{\beta^2}{\lambda_{t,d}}\cdot \frac{c^2c_0^2}{16}}{\left(\sqrt{\sum_{j=1}^{d}\nu_{t,j}^2}+\frac{\beta}{\sqrt{\lambda_{t,d}}}\cdot \frac{c}{4}\right)^2}\nonumber\\
    &\geq \frac{c^4c_0^2}{25},
\end{align}
when $c$ is set small enough, which concludes our claim when we set $C_1 = 16c_0^2/c^2-1$ and $C_2= c^4/25$.

\subsection{Proof of Lemma \ref{lem:opt-comparison}}

    We rewrite the problems as follows. Denote
    \begin{align*}
        \alpha_i = \frac{\beta}{\sqrt{\lambda_{t,i}}},\quad b_i^{(1)} = \alpha_i\nu_{t,i},\quad b_i^{(2)} = \alpha_1 \mathbf{1}(i=1).
    \end{align*}
    Then the $\bm{w}_t^{\star}$ is the solution of
    $$
        f_1(\bm{w}) = \sum_{i=1}^{d}2\alpha_ib_i^{(1)}w_i + \alpha_i^2 w_i^2,
    $$
    while $\bm{\widetilde{w}}_t^{\star}$ is the solution of 
    $$
        f_2(\bm{w}) = \sum_{i=1}^{d}2\alpha_ib_i^{(2)}w_i + \alpha_i^2 w_i^2.
    $$
    The KKT condition gives the solution of both optimization problems to be 
    \begin{align}
        w_{i}^{\star} = \frac{b_i^{(1)}}{\mu_1-\alpha_i^2},\quad \widetilde{w}_{i}^{\star} = \frac{b_i^{(2)}}{\mu_2-\alpha_i^2},
    \end{align}
    where $\mu_1,\mu_2$ satisfies
    $$\sum_{i=1}^{d} \left(\frac{b_i^{(1)}}{\mu_1-\alpha_i^2}\right)^2 = 1,\quad \sum_{i=1}^{d} \left(\frac{b_i^{(2)}}{\mu_2-\alpha_i^2}\right)^2 = 1.$$
    A direct calculation yields that $\mu_1>\alpha_d^2$ but $\mu_2 = \alpha_d^2$. We also note that $|b_i^{(1)}|\leq |b_i^{(2)}|$. Therefore, one obtain the following inequality
    \begin{align}
        |w_{t,1}^{\star}|  = \frac{|b_i^{(1)}|}{\mu_1-\alpha_1^2}\leq \frac{|b_i^{(2)}|}{\mu_2-\alpha_1^2} = |\widetilde{w}_{t,1}^{\star}|,
    \end{align}
    establishing the desired result.

\subsection{Proof of Lemma \ref{lem:growth-largest-eigenvalue}}
    As we note that $\lambda_{t+1,1}$ is the largest root of the following equation
    \begin{align}
        f(\lambda) = 1 + \sum_{i=1}^{d}\frac{\kappa_{t,i}^2}{\lambda_{t,i}-\lambda}.
    \end{align}
    Set $\widetilde{\lambda} = \lambda_{t,1} + \kappa_{t,1}^2$, then as $\lambda_{t,1}-\lambda_{t,i}\geq t/2$, for all $i\geq 2$, we have
    \begin{align}
        f(\widetilde{\lambda}) =  \sum_{i=2}^{d}\frac{\kappa_{t,i}^2}{\lambda_{t,i}-\lambda} = O(t^{-1}).
    \end{align}
    Furthermore, we also note that as $\kappa_{t,1}^2/2\leq |\lambda-\lambda_{t,1}|\leq 3\kappa_{t,1}^2/2$, it holds that
    \begin{align}
        f^\prime(\lambda) = -\sum_{i=1}^{d}\frac{\kappa_{t,i}^2}{(\lambda_{t,i}-\lambda)^2} \leq -\frac{4}{9\kappa_{t,1}^2} + O(t^{-2})\leq -\frac{4}{9} + O(t^{-2}).
    \end{align}
    Therefore, we conclude that
    \begin{align}
        |\lambda_{t,1}-\widetilde{\lambda}| = O(t^{-1}),
    \end{align}
    implying that $\lambda_{t+1,1} = \lambda_{t,1}+\kappa_{t,1}^2+ O(t^{-1})$.
    We will then further characterize the growth rate by plugging in $\kappa_{t,1}$. From Lemma \ref{lem:decom-a},  one can show that
    \begin{align}
        \kappa_{t,1}
        & = \left(1- O\left(\frac{\beta^2}{\lambda_{t,d}}\right)\right)\left(1-O\left(\frac{\beta^2}{\lambda_{t,d}}\right)\right) + O\left(\frac{\beta}{\sqrt{\lambda_{t,d}}}\cdot \frac{\beta}{\sqrt{\lambda_{t,d}}}\right)\nonumber\\
        & = 1-O\left(\frac{\beta^2}{\lambda_{t,d}}\right),
    \end{align}
    Therefore, one can show that
    $$\lambda_{t+1,1} = \lambda_{t,1} + 1-O\left(\frac{\beta^2}{\lambda_{t,d}}\right) + O(t^{-1}) = \lambda_{t,1} + 1-O\left(\frac{\beta^2}{\lambda_{t,d}}\right) ,$$
    where the last equality holds as $\lambda_{t,d}\asymp \beta\sqrt{t}$.

\subsection{Proof of Lemma \ref{lem:approximate-condition}}
We begin with the following transformation. Let $\tilde{\bm{\lambda}} = \bm{\Lambda}^{1/2}\bm{\lambda}$ and $\tilde{\bm{\eta}} = \bm{\Lambda}^{-1/2}\bm{\eta}$. Then $f(\bm{\lambda})$ can be rewritten as
$$f(\bm{\lambda}) = \tilde{\bm{\lambda}}^{T}\tilde{\bm{\eta}} - \frac{\sigma^2}{2}\tilde{\bm{\lambda}}^{T}\tilde{\bm{\lambda}}=-\frac{\sigma^2}{2}\left\Vert\tilde{\bm{\lambda}} - \frac{1}{\sigma^2}\tilde{\bm{\eta}}\right\Vert_2^2+\frac{1}{2\sigma^2}\Vert\tilde{\bm{\eta}}\Vert_2^2,$$
with $\max_{\bm{\lambda}\in\mathbb{R}^d} f(\bm{\lambda}) = (2\sigma^2)^{-1}\Vert\tilde{\bm{\eta}}\Vert_2^2$. Then for any $\bm{\lambda}\in\mathcal{C}$, it holds that
$$-\frac{\sigma^2}{2}\left\Vert\tilde{\bm{\lambda}} - \frac{1}{\sigma^2}\tilde{\bm{\eta}}\right\Vert_2^2+\frac{1}{2\sigma^2}\Vert\tilde{\bm{\eta}}\Vert_2^2\geq \frac{1-\kappa^2}{2\sigma^2}\Vert\tilde{\bm{\eta}}\Vert_2^2,$$
which immediately implies
\begin{align}\label{equ:C-equivalent}
    \left\Vert\tilde{\bm{\lambda}} - \frac{1}{\sigma^2}\tilde{\bm{\eta}}\right\Vert_2\leq \frac{\kappa}{\sigma^2}\Vert\tilde{\bm{\eta}}\Vert_2.
\end{align}
Substituting $\bm{\lambda} = \bm{\Lambda}^{-1/2}\tilde{\bm{\lambda}}$ and $\bm{\eta} = \bm{\Lambda}^{1/2}\tilde{\bm{\eta}}$ into (\ref{equ:C-equivalent}) yields the desired result.

\subsection{Proof of Lemma~\ref{lem:sequence-induction}}

   
    We set $a_n = A_n n^{-1/4}$, then it holds that
    \begin{align}
        A_{n+1}^2 (n+1)^{-1/2} & \leq \left(\left(1-\frac{1}{n}\right)A_n n^{-1/4} + \frac{B}{n^{5/4}}\right)^2 + \frac{C}{n^{2}} \nonumber\\
        & = \left[\left(\left(1-\frac{1}{n}\right)A_n + \frac{B}{n}\right)^2  + \frac{C}{n^{3/2}}\right]\cdot n^{-1/2}.
    \end{align}
    From the basic inequality that $(1+1/n)^{1/2}\leq 1+1/(2n)$, one can show that
    \begin{align}
    \label{equ:A-n-induction-1}
        A_{n+1}^2 & \leq \left(1+\frac{1}{2n}\right)\cdot \left[\left(\left(1-\frac{1}{n}\right)A_n + \frac{B}{n^{5/4}}\right)^2  + \frac{C}{n^{2}}\right]\nonumber\\
        & = \left(1+\frac{1}{2n}\right)\cdot \left[\left(\left(1-\frac{1}{n}\right) + \frac{B}{n^{5/4}A_n}\right)^2  + \frac{C}{n^{2}A_n^2}\right]\cdot A_n^2\nonumber\\
        & \leq \left(1+\frac{1}{2n}\right)\cdot \left(1-\frac{2}{n} + \frac{1}{n^2} +\frac{2B}{n^{5/4}A_n}  + \frac{B^2}{n^{5/2}A_n^2} + \frac{C}{n^{3/2}A_n^2}\right)\cdot A_n^2.
    \end{align}
    Note that when
    \begin{align*}
        n\geq \max\left(12, \frac{4096B^4}{A_n^4},\frac{12^{2/3}B^{4/3}}{A_n^{4/3}}, \frac{144C^2}{A_n^4} \right),
    \end{align*}
    we have 
    \begin{align*}
        \frac{1}{n^2}\leq\frac{1}{12n},\quad \frac{2B}{n^{5/4}A_n}\leq\frac{1}{4n},\quad \frac{B^2}{n^{5/2}A_n^2}\leq \frac{1}{12n},\quad \frac{C}{n^{3/2}A_n^2}\leq\frac{1}{12n}.
    \end{align*}
    then it holds that 
    \begin{align}
    \label{equ:A-n-induction-2}
        A_{n+1}^2&\leq\left(1-\frac{1}{2n}\right)\cdot \left(1-\frac{2}{n}+\frac{1}{12n}+\frac{1}{4n}+\frac{1}{12n}+\frac{1}{4n}\right)A_n^2\nonumber\\
        &\leq \left(1-\frac{1}{n}\right)A_n^2.
    \end{align}
    Now that we constrain on $n$ such that $A_n\geq 4B$, then (\ref{equ:A-n-induction-2}) holds as long as 
    \begin{align*}
        n\geq \max\left(16, \frac{9C^2}{16B^4}\right),
    \end{align*}
   
    Therefore, let $n_1 = \min\{n: A_{n+1}\leq 4B\}$, then for all $n_0<n\leq n_1$, it holds that 
    \begin{align}
    \label{equ:upper-bound-A-n}
        A_n^2\leq \prod_{i=n_0}^{n}\left(1-\frac{1}{i}\right)\cdot A_{n_0}^2\leq \frac{n_0}{n}\cdot A_{n_0}^2,
    \end{align}
    as $A_{n_1}\geq 4B$, (\ref{equ:upper-bound-A-n}) implies
    \begin{align*}
        \frac{n_0}{n_1}\cdot A_{n_0}^2\geq (4B)^2,
    \end{align*}
    which allows us to upper bound $n_1$ as
    \begin{align}
    \label{equ:upper-bound-A-n-2}
        n_1\leq \left(\frac{A_{n_0}}{4B}\right)^2n_0 = \frac{a_{n_0}^2n_0^{3/2}}{16B^2}.
    \end{align}
     We also note that the right hand side of (\ref{equ:A-n-induction-1}) is monotone with respect to $A_n$, and when $A_n = 4B$, it holds that 
    \begin{align*}
        A_{n+1}^2\leq \left(1-\frac{1}{n}\right)\cdot (4B)^2 = (4B)^2,
    \end{align*}
    which implies that for any $n$ such that $A_n\leq 4B$, we have 
    \begin{align}
        A_{n+1}\leq 4B.
    \end{align}
    Therefore, for all $n\geq n_1$, it holds that $A_n\leq 4B$. Combining with (\ref{equ:upper-bound-A-n-2}), we arrive at our final result.

\label{sec:ucb-constrained}

\bibliography{ref}
\bibliographystyle{apalike}

\end{document}